\documentstyle{amsart}

\newtheorem{Theorem}{Theorem}[section]
\newtheorem{Lemma}[Theorem]{Lemma}

\newtheorem{Remark}[Theorem]{Remark}

\newtheorem{Definition}[Theorem]{Definition}

\begin{document}

\title{Local Monomialization of Transcendental Extensions}

\author{Steven Dale Cutkosky}
\thanks{  partially supported by NSF}

\maketitle
\section{Introduction}

Suppose that $R\subset S$ are local domains such that $S$ dominates $R$. We will say that $R\subset S$ is monomial
if $R$ and $S$ are regular and there are regular system of parameters $(x_1,\ldots,x_m)$
in $R$ and $(y_1,\ldots, y_n)$ in $S$,   there are units $\delta_1,\ldots,\delta_m$
in $S$ and an $m\times n$ matrix $A$ of natural numbers such that $A$ has maximal rank $m$ and
$$
\begin{array}{ll}
x_1&=y_1^{a_{11}}\cdots y_n^{a_{1n}}\delta_1\\
&\vdots\\
x_m&=y_1^{a_{m1}}\cdots y_n^{a_{mn}}\delta_m.
\end{array}
$$

If $R$ is of equicharacteristic zero, and $R\subset S$ is monomial then there is a finite \'etale extension $S\rightarrow \overline S$, obtained by extracting roots of the $\delta_i$,
such that $R\subset\overline S$ is monomial, with all $\delta_i=1$ in the above system of
equations.

The condition of  being monomial is very special.
If $R$ and $S$ are localizations of polynomial rings over a field, the variables of $R$
will be related to those of $S$ by  polynomials, which are of course in general definitely not monomial if $R$ and $S$ have dimension greater than 1.

In this paper, we prove a very general theorem (Theorem \ref{Theorem2}) showing that
when $R\subset S$ are essentially of finite type over a field of characteristic zero,
they are in fact very close to being monomial. Before we precisely state our result,
we will explain the philosophy behind it.

\subsection{Monoidal transforms, valuations rings and birational local algebra}
A particular example of a monomial extension is a monoidal transform of a regular local ring. Suppose that $R$ is a regular local 
ring, and that $P\subset R$ is a regular prime; that is, $R/P$ is a regular local ring. Let 
$f$ be a nonzero element of $P$, and suppose that $m$ is a prime ideal of $R[\frac{P}{f}]$
that contracts to the maximal ideal of $R$. Let $R_1=R[\frac{P}{f}]_m$.
$R\rightarrow R_1$ is called a monoidal transform. A monoidal transform takes on a particularly
simple form if $R$ contains an algebraically closed  field isomorphic to its residue field
 and $m$ is a maximal ideal of $R[\frac{P}{f}]$. Then there are systems
of regular parameters $(x_1,\ldots,x_m)$ in $R$ and $(y_1,\ldots,y_m)$ in $R_1$,  such that
$$
x_i=\left\{\begin{array}{ll}
y_iy_r&\text{ for }1\le i\le r-1\\
y_i&\text{ for }r<i\le \text{ dim }R
\end{array}\right.
$$
where  $r=\text{ height }(P)$.

Of course, if we iterate monoidal transforms, the composition will in general
not be monomial.

A monoidal transform is  just a local ring of the blowup  
$X=\text{proj}(\oplus_{n\ge 0}P^n)$ of $P$,  which is proper over $\text{spec}(R)$. It has many local rings which dominate $R$.
The whole blowup scheme contains all information about $R$, but it focuses attention on
some infinitesimal properties of $R$ at its maximal ideal.  If we want to interpret this in terms of
local algebra, we immediately reach a difficulty. Some local rings in  $X$ may be
more interesting than others, but there is in general no well defined way to pick out
a most important local ring.  We will state the problem precisely.

Suppose that  $K$ is a field.  
If $R$ is a Noetherian local domain with quotient field $K$, and $I$ is an ideal of $R$,
we want to determine a unique local ring $R_1$ 
 of the proper $R$ scheme $X=\text{proj}(\oplus_{n\ge 0}I^n)$ which dominates $R$.

There is a  solution to this problem  through valuation theory. We fix a valuation
$\nu$ of $K$ such that its valuation ring $V$ dominates $R$.  We choose $f\in R$
such that $\nu(f)=\text{ min }\{\nu(f)\mid f\in I\}$.   $R[\frac{I}{f}]$ is
contained in $V$, and if we let $m=m_V\cap R[\frac{I}{f}]$, where $m_V$ is the maximal ideal of $V$, then
$R_1=R[\frac{I}{f}]_m$ is the unique local ring of the blowup $X$ of $I$ which is
dominated by the valuation $\nu$.  If $I$ is a regular prime, we say that $R\rightarrow R_1$ is a
monoidal transform along $\nu$.

For this approach to be meaningful, we must consider all types of valuations which dominate
$R$.  As an example, suppose that $R$ is an excellent, normal two dimensional local ring.     
If we blow up the maximal ideal of $R$ and normalize the resulting scheme, there will only
be isolated singularities.  If the scheme is not regular, we can take $R_1$ to be the local
ring of one of the singular points.
Now we blow up the maximal ideal of $R_1$ and normalize. If this scheme is not regular,
we can take $R_2$ to be the local ring of a singular point.   We can iterate this procedure,
and obtain a birational sequence of  normal local rings
$$
R\rightarrow R_1\rightarrow R_2\rightarrow\cdots
$$
which continues as long as the normalization of the blowup of the maximal ideal of
$R_i$ is singular.

This problem can be easily interpreted in valuation theory.  If the sequence is infinite,
then we can consider the  subring $\cup_{i=0}^{\infty}R_i$ of 
the quotient field $K$ of $R$.  There exists a valuation ring $V$ of $K$ which dominates $R$.
In fact, since we are in dimension two, $\cup_{i=0}^{\infty}R_i$ is a valuation ring.
The valuation ring dominates each $R_i$ in the above sequence.
In general, $\cup_{i=0}^{\infty}R_i$ is not Noetherian, so we are forced
to consider  general valuations, whose value groups are much larger than
${\bf Z}$.    For instance, even on an algebraic function field of dimension two, there are valuations with value group ${\bf Q}$ (c.f. Example 3, Section 15, Chapter VI \cite{ZS}).  

The problem of resolution of singularities is formulated in
local birational algebra as local uniformization.  Suppose that
$R$ is a local domain with quotient field $K$, and $V$ is a valuation 
ring of $K$ which dominates $R$.  The problem of local uniformization
is to find a sequence of monoidal transforms  
$R\rightarrow R'$ of $R$ along $V$  such  that $R'$ is  a regular local
ring.

Local uniformization was originally proven in all dimensions for local rings essentially
of finite type over a field of characteristic 0 by Zariski  \cite{Z1}, \cite{Z2}.
Of course this is now also a consequence of Hironaka's general theorems on resolution
of singularities in characteristic 0 \cite{H1}.  Local uniformization is known to be true
for excellent local rings of dimension two,  (\cite{Z1}, \cite{Ab1}, \cite{H2}, \cite{L3}), and 
resolution of singularities is known for excellent surfaces of dimension two. Local uniformization is known  for local rings which are essentially of finite type over an
algebraically closed field of characteristic $p\ne 2,3$ or 5 \cite{Ab5}.
Recently there has been meaningful progress in local uniformization in higher dimensions
over fields of positive characteristic (\cite{Co}, \cite{HRW},\cite{Ku},\cite{Sp},\cite{T}).

For varieties of  dimension $\le 3$,  Local Uniformization implies resolution of singularities
(\cite{Z5},  \cite{Ab5}), but the validity of this implication is unknown in higher dimensions.

 \subsection{Local Monomialization} 

 We  now  state the problem of monomialization of our extension
of local rings $R\subset S$ in birational local algebra.

Suppose that 
$R\subset S$
is a local homomorphism of local rings essentially of finite type over a field $k$
and that $V$ is a valuation ring of the quotient field $K$ of $S$, such that $V$ dominates $S$. 
Then we can ask if there are sequences of monoidal transforms $R\rightarrow R'$ and $S\rightarrow S'$ along $V$ such that
$V$ dominates $S'$, $S'$ dominates $R'$, and $R\rightarrow R'$ is a monomial mapping. 
$$
\begin{array}{lll}
R'&\rightarrow&S'\subset V\\
\uparrow&&\uparrow\\
R&\rightarrow &S
\end{array}
$$

We completely answer this question in the affirmative when $k$ has characteristic 0 in Theorem \ref{Theorem2}. Notations are as in Section \ref{Not}

\begin{Theorem}\label{Theorem2} Suppose that $k$ is a field of
characteristic zero, $K\rightarrow K^*$ is a (possibly transcendental) extension
of algebraic function fields over $k$.  Suppose that $\nu^*$ is a valuation of
$K^*$ which is trivial on $k$. Suppose that $R$ is an algebraic local ring of $K$,
$S$ is an algebraic local ring of $K^*$ such that $S$ dominates $R$ and $\nu^*$ dominates
$S$.  
Then there exist sequences of monoidal transforms $R\rightarrow R'$ and $S\rightarrow S'$
along $\nu^*$ such that $R'$ and $S'$ are regular local rings, $S'$ dominates $R'$,
there exist regular parameters
$(y_1,\ldots,y_n)$ in $S'$, $(x_1,\ldots,x_m)$ in $R'$, 
units $\delta_1,\ldots,\delta_m\in S'$ and an $m\times n$
matrix $(c_{ij})$ of nonnegative integers such that $(c_{ij})$ has rank $m$, and 
\begin{equation}\label{eqI1}
x_i=\prod_{j=1}^ny_j^{c_{ij}}\delta_i
\end{equation}
for $1\le i\le m$.
\end{Theorem}

The case when the field extension $K\rightarrow K^*$ is finite is solved in Theorem 1.1
of our paper \cite{C1}.

The standard theorems on resolution of singularities allow one to easily find $R'$ and $S'$ such that (\ref{eqI1}) holds,
but, in general,  the essential condition that $(a_{ij})$ has maximal rank $m$ will not hold.
It is for this reason that we must construct the sequence of monoidal transforms
$R\rightarrow R'$, even if $R$ is regular. 
The difficulty of the
proof of the Theorem is to achieve this condition.

It is an interesting open problem to prove Theorem \ref{Theorem2} in  positive
characteristic, even in dimension 2 (\cite{CP}, \cite{CP2}).  

In \cite{C3} and \cite{CP2} we use Theorem 1.1 \cite{C1}, which is the finite field
extension $K\rightarrow K^*$ case of Theorem \ref{Theorem2} to prove very strong results
in the ramification theory of general valuations on characteristic zero algebraic
function fields, such as  Abhyankar's `` Weak local simultaneous resolution conjecture''. 
It is expected that Theorem \ref{Theorem2} can be used to extend this ramification theory
to arbitrary extensions of characteristic zero algebraic function fields.

We will make a few comments here about the proof of Theorem \ref{Theorem2}. Our
starting point is the  proof for finite extensions $K\rightarrow K^*$ of our paper
\cite{C1}. 
An overview of the proof (in the finite field extension case) can be found in  Section 1.3 of \cite{C1}. 

 Some parts of this proof generalize readily to the case when $K^*$ is transcendental over $K$.    For these parts, we give here the modified statements, and indicate the changes
which must be made in the original proofs.  However, there are some parts of the proof which
are quite different. The really new ingrediants in the proof are given in the critical sections 
\ref{UTS} and \ref{CUTS} of this paper.   
As in the proof for the case when $K\rightarrow K^*$ is finite, we reduce to the
case when $V^*$ has rank 1. Since $V= V^*\cap K$ then has rank $\le 1$, and we can assume
that $V$ is nontrivial, we are reduced to the case when $V$ has rank 1 also.
Two new complexities arise in the case when $K^*$ is transcendental over $K$. The rational
rank of a valuation $\nu$ is the dimension of the ${\bf Q}$ vector space $\Gamma_{\nu}\otimes{\bf Q}$, where $\Gamma_{\nu}$ is the valuation group of $\nu$.  
We have an inequality $\overline r=\text{ ratrank }(\nu)\le \overline s=\text{ ratrank }(\nu^*)$. If $K^*$ is finite over $K$ this is an equality.  The case when 
$\overline r=\text{ ratrank }(\nu)< \overline s=\text{ ratrank }(\nu^*)$ is significantly
more difficult. It is addressed in Section \ref{CUTS}. The second major new complexity 
lies in the extension of residue fields of valuations.  If $k(V^*)$ is the residue
field of $V^*$ and $k(V)$ is the residue field of $V$, then we have
$\text{ trdeg }_{k(V)}k(V^*)\le \text{ trdeg }_KK^*$. Thus $k(V^*)$ is algebraic
(though not generally finite) over $k(V)$ if $K^*$ is finite over $K$. The new arguments
which are required to handle the case when $k(V^*)$ is transcendental over $k(V)$
 are in Sections \ref{UTS} and \ref{CUTS}.

\subsection{Monomialization of morphisms of varieties}\label{section3}

We discuss an application of our Theorem \ref{Theorem2} to proper morphisms of
varieties, and the problem of monomialization of morphisms of varieties.

\begin{Definition} Suppose that $\Phi:X\rightarrow Y$ is a dominant morphism
of nonsingular integral finite type $k$ schemes. $\Phi$ is  monomial if for every
$p\in X$ there exist regular parameters $(y_1,\ldots,y_m)$ in 
${\cal O}_{Y,\Phi(p)}$, and an \'etale cover $U$ of an affine neighborhood of $p$, uniformizing parameters $(x_1,\ldots,x_n)$ on $U$ and a
matrix $a_{ij}$ such that
$$
\begin{array}{ll}
y_1=&x_1^{a_{11}}\cdots x_n^{a_{1n}}\\
&\vdots\\
y_m=&x_1^{a_{m1}}\cdots x_n^{a_{mn}}
\end{array}
$$
\end{Definition}

We do not assume that $X$ and $Y$ are separated in the above definition.
Since $\Phi$ is dominant, the matrix $(a_{ij})$ must have maximal rank $m$.

A quasi-complete variety over a field $k$ is an integral finite type $k$-scheme
which satisfies the existence part of the valuative criterion for properness
(Hironaka, Chapter 0, Section 6 of \cite{H1} and Chapter 8 of \cite{C1}).

The construction of a monomialization by  quasi-complete varieties follows from Theorem \ref{Theorem2}.

\begin{Theorem}\label{TheoremB} Let $k$ be a field of characteristic zero, $\Phi:X\rightarrow Y$ a dominant 
morphism of   proper  $k$-varieties.   Then there are
  birational morphisms of nonsingular quasi-complete  $k$-varieties $\alpha:X_{1}\rightarrow X$
and $\beta: Y_{1}\rightarrow Y$, and a  monomial morphism $\Psi:X_1\rightarrow Y_1$ such that  the diagram
\[
\begin{array}{lll}
X_1&{\Psi}\atop{\rightarrow}&Y_1\\
\downarrow&&\downarrow\\
X&{\Phi}\atop{\rightarrow}&Y
\end{array}
\]
commutes and $\alpha$ and $\beta$ are
 locally  products of blow ups of nonsingular subvarieties. 
That is, for every $z\in X_{1}$, there
exist affine neighborhoods $V_1$ of $z$, $V$ of $x=\alpha(z)$, such that $\alpha:V_1\rightarrow V$ 
is a  finite product of monoidal transforms, and  
there
exist affine neighborhoods $W_1$ of $\Psi(z)$, $W$ of $y=\beta(\Psi(z))$, such that $\beta:W_1\rightarrow W$ 
is a  finite product of monoidal transforms.
\end{Theorem}

A monoidal transform of a nonsingular $k$-scheme $S$ is the map $T\rightarrow S$ induced by an open subset $T$ of
$\text{Proj}(\oplus {\cal I}^n)$, where ${\cal I}$ is the ideal sheaf of a nonsingular subvariety of $S$.

The case of Theorem \ref{TheoremB} when $X\rightarrow Y$ is generically finite is proven in Theorem 1.2 of our paper \cite{C1}.

The proof of Theorem \ref{TheoremB} in general follows from Theorem \ref{Theorem2}, by
patching a finite number of local solutions, as in the proof of Theorem 1.2 \cite{C1}. 
The resulting schemes may not be separated.

It is an extremely interesting question to determine if the conclusions of Theorem \ref{TheoremB} hold, but with the stronger conditions that $\alpha$ and $\beta$ are products
of monoidal transforms on proper varieties $X_1$ and $Y_1$.

The strongest known result on monomialization is  our theorem below.

\begin{Theorem}\label{Mon}(Theorem 18.21, \cite{C5})
Suppose that $\Phi:X\rightarrow S$ is a dominant morphism from a 3 fold $X$ to
a surface $S$ (over an algebraically closed field $k$ of characteristic zero). Then there exist sequences of
blow ups of nonsingular subvarieties $X_1\rightarrow X$ and $S_1\rightarrow S$ such
that the induced map $\Phi_1:X_1\rightarrow S_1$ is a monomial morphism.
\end{Theorem}

A generalization of this result to prove monomialization of 
strongly prepared morphisms from $N$-folds to surfaces appears in the
paper \cite{CK} with Olga Kashcheyeva.

The author would like to thank the Mathematical Sciences Research Institute for
its hospitality while this manuscript was being prepared.

\section{Notations}\label{Not} 
We will denote the
maximal ideal of a local ring $R$ by $m(R)$. If $R$ contains a field $k$, we will denote
its residue field by $k(R)$. 
We will denote the  quotient field of a domain $R$ by $Q(R)$.
Suppose that $R\subset S$ is an inclusion of local rings. We will say that $R$ dominates
$S$ if $m(S)\cap R=m(R)$.  
Suppose that $K$ is an algebraic function field over a field $k$. We will say that
a local ring  $R$ with quotient  field $K$ is an algebraic local ring of $K$ if $R$ is essentially of finite type over $k$.
If $R$ is a local ring, $\hat R$ will denote the completion of $R$ at its maximal ideal.
If $L_1$ and $L_2$ are 2 subfields of a field $M$, then $L_1*L_2$ will denote the subfield
of $M$ generated by $L_1$ and $L_2$.

Good introductions to the valuation theory which we require in this paper can be found  in Chapter VI of \cite{ZS} and in \cite{Ab4}.  
A valuation $\nu$ of $K$ will be called a $k$-valuation if $\nu(k)=0$. 
We will denote by $V_{\nu}$ the associated valuation ring, which necessarily contains $k$. 
A valuation ring
$V$ of $K$ will be called a $k$-valuation ring if $k\subset V$.  The value group of a valuation $\nu$ with valuation ring $V$ will be denoted by $\Gamma_{\nu}$ or $\Gamma_V$. 
We will abuse notation by denoting the valuation $\nu$, which is a homorphism of the
group of units of $K$, as a function on $K$.
If $R$ is a subring of $V_{\nu}$ then the center of $\nu$ (the center of $V_{\nu}$)
on $R$ is the prime ideal $R\cap m(V_{\nu})$.  If $R$ is a Noetherian subring of $V_{\nu}$
and $I\subset R$ is an ideal, we will write $\nu(I)=\rho$ if $\rho=\text{ min }\{\nu(f)\mid f\in I\}$.

We will review the concept of composite valuations. For details, we refer to Section 10 of
Chapter II of \cite{Ab4} and Section 10, Chapter VI \cite{ZS}.
If $\nu$ is a valuation of rank greater than 1, then $\nu$ is a composite valuation.
That is, there are valuations $w$ and $\overline \nu$ where $w$ is a valuation of
$K$ and $\overline \nu$ is a valuation of the residue field of $V_w$ such that
if $\pi:V_w\rightarrow k(V_w)$ is the residue map, then $V_v=\pi^{-1}(V_{\overline \nu})$.
For $f\in V_w$ such that $\pi(f)\ne0$ we have $\nu(f)=\overline \nu(\pi(f))$. 
This gives us an inclusion of value groups $\Gamma_{\overline\nu}\subset \Gamma_{\nu}$.
 $\Gamma_{\overline\nu}$ is an isolated subgroup of $\Gamma_{\nu}$.
There exists a prime ideal $p$ in $V_{\nu}$ such that $V_{w}=(V_{\nu})_p$. For $f\in K$, $w(f)$ is the residue of $\nu(f)$   in $\Gamma_w=\Gamma_{\nu}/\Gamma_{\overline\nu}$.
We say that $\nu$ is the composite of $w$ and $\overline\nu$ and write
$\nu=w\circ\overline \nu$.

Suppose that $R$ is a local domain. A monoidal transform $R\rightarrow R_1$ is a 
birational extension of local domains such that $R_1=R[\frac{P}{x}]_m$ where $P$ is
a regular prime ideal of $R$, $0\ne x\in P$ and $m$ is a prime ideal of $R[\frac{P}{x}]$
such that $m\cap R=m(R)$. $R\rightarrow R_1$ is called a quadratic transform  if $P=m(R)$.

If $R$ is regular, and $R\rightarrow R_1$ is a monodial transform, then there exists a regular sustem of parameters $(x_1,\ldots, x_n)$ in
$R$ and $r\le n$ such that
$$
R_1=R\left[\frac{x_2}{x_1},\ldots,\frac{x_r}{x_1}\right]_m.
$$

Suppose that $\nu$ is a valuation of the quotient field $R$ with valuation ring $V_{\nu}$
which dominates $R$. Then $R\rightarrow R_1$ is a monoidal transform along $\nu$
(along $V_{\nu}$) if $\nu$ dominates $R_1$.

\section{valuations}
\begin{Lemma}\label{Lemma4}
Suppose that $K$ is a field containing a subfield $k$, $t_1,\ldots, t_{\alpha}$ are 
algebraically independent over $K$ and $\nu$ is a $k$-valuation of $K$
with  valuation ring $V$. Then there
exists a unique extension $\overline\nu$ of $\nu$ to $K(t_1,\ldots,t_{\alpha})$,
such that $\overline\nu(f)=\nu(f)$ for $f\in K$, $\overline\nu(t_i)=0$,
 and if $\overline V$ is the valuation ring of $\overline \nu$, then the images of $t_1,\ldots,t_{\alpha}$ in $k(\overline V)$ are 
algebraically independent over $k(V)$.
\end{Lemma}

\begin{pf} For $I=(i_1,\ldots,i_{\alpha})\in {\bf N}^{\alpha}$, let $t^I=t_1^{i_1}\ldots t_{\alpha}^{i_{\alpha}}$.  If
$$
0\ne h=\sum_{I\in{\bf N}^{\alpha}} f_It^I\in K[t_1,\ldots,t_{\alpha}]
$$
 with $f_I\in K$,
 define
$$
\overline \nu(f)=\text{min}\{\nu(f_I)\mid f_I\ne 0\}.
$$
This induces an extension of $\nu$ as desired.
We will verify that $\overline\nu(fg)=\overline\nu(f)+\overline \nu(g)$
for 
$$
f=\sum_If_It^I, g=\sum_J g_jt^J\in K[t_1,\ldots,t_{\alpha}].
$$  
$$
fg=\sum_A\left(\sum_{I+J=A}f_Ih_J\right)t^A.
$$
For each $A$,
$$
\nu\left(\sum_{I+J=A}f_Ig_J\right)\ge\text{ min }\{\nu(f_I)+\nu(g_J)\mid I+J=A\}
\ge \overline \nu(f)+\overline\nu(g).
$$
Let $I_0$ be such that $\nu(f_{I_0})>\nu(f_I)$ if $I<I_0$
(in the Lex order) and $\nu(f_I)\ge \nu(f_{I_0)})$ if $I>I_0$.
Similarily, let $J_0$ be such that $\nu(g_{J_0})>\nu(g_J)$ if $J<J_0$
 and $\nu(g_J)\ge \nu(g_{J_0)})$ if $J>J_0$. Let $A_0=I_0+J_0$. Then
$$
\nu(\sum_{I+J=A_0}f_Ig_J)=\nu(f_{I_0}g_{J_0})
=\overline\nu(f)+\overline\nu(g).
$$
Thus $\overline\nu(fg)=\overline\nu(f)+\overline\nu(g)$.

Suppose that $\overline\nu$ is an extension of $\nu$ with the desired properties. If $h=\sum f_It^I\in K[t_1,\ldots,t_{\alpha}]$
with $f_I\in K$,
let $f_J$ be such that $\nu(f_J)=\text{min}\{\nu(f_I)\}$.
If $\overline\nu(h)>\nu(f_J)$, then $\overline \nu(\sum _I\frac{f_I}{f_J}t^I)>0$
so that 
$$
[\sum\frac{f_I}{f_J}t^J]=0\text{ in }k(\overline V),
$$
where $[\beta]$ denotes the class of $\beta\in \overline V$ in $k(\overline V)$.
But by assumption, $[t_1],\ldots,[t_{\alpha}]$ are algebraically independent in $k(\overline V)$
over $k(V)$. This is a contradiction.
\end{pf}

\begin{Lemma}\label{Lemma25}  Suppose that $K$ is a field containing a subfield $k$,
$t_1,\ldots, t_{\alpha}$ are analytically independent over $K$ and $\nu$ is a $k$-valuation
of $K$ with valuation ring $V$. Suppose that $R$ is a noetherian local domain with quotient
field $K$ such that $V$ dominates $R$. Then there exists a unique extension $\overline\nu$
of $\nu$ to $Q(R[[t_1\ldots,t_{\alpha}]])$ such that $\overline\nu(f)=\nu(f)$ for $f\in K$,
$\overline \nu(t_i)=0$ for $1\le i\le\alpha$, and if $\overline V$ is the valuation ring of $\overline\nu$,
then the images of $t_1,\ldots,t_{\alpha}\in k(\overline V)$ are analytically independent
over $k(V)$.
\end{Lemma}

\begin{pf} For a series
$$
f=\sum a_It^I\in R[[t_1,\ldots, t_{\alpha}]]
$$
with $a_I\in R$, we define
$$
\overline\nu(f)=\text{ min }\{\nu(f_I)\mid f_I\in R\}.
$$
We first verify that $\overline\nu$ is well defined. Suppose that 
$$
\text{ min }\{\nu(f)\mid f_I\in R\}
$$
does not exist. Then there exists an infinite descending chain of values
$$
a_1>a_2>a_3>\cdots>0
$$
 and $f_i\in R$ such that $\nu(f_i)=a_i$ for all positive
integers $i$.  Let $I_i$ be the $R$ ideal
$$
I_i=\{g\in R\mid \nu(g)\ge a_i\}.
$$
Then we have an infinite strictly ascending chain of ideals in $R$,
$$
I_1\subset I_2\subset\cdots,
$$
a contradiction to the assumption that $R$ is Noetherian.

As in the proof of Lemma \ref{Lemma4}, $\overline\nu$ induces an extension of $\nu$
as desired.
As in the proof of Lemma \ref{Lemma4}, $\overline\nu$ is unique.
\end{pf}



\section{Rational Rank 1 Valuations}\label{RR1V}
Suppose that $k$ is a field of characteristic 0 and $K\rightarrow K^*$ is
an extension of algebraic function fields over $k$. Suppose that $\nu^*$
is a rank 1 $k$-valuation of $K^*$ with valuation ring $V^*$. Let
$\nu=\nu^*\mid K$ with valuation ring $V=V^*\cap K$. Necessarily, $\nu$ has rank $\le 1$
(c.f. Lemma \ref{Lemma10} and the
discussion following Lemma \ref{Lemma10}). 
Assume that 
$\nu$ has rank 1 and that $k(V)$ is algebraic over $k$.  

Let $\overline r=\text{ratrank }\nu$, $\overline s=\text{ratrank }\nu^*$
be the respective rational ranks.
Let 
$m=\text{trdeg}_k K$, $n=\text{trdeg}_k K^*
-\text{trdeg}_k k(V^*)$.
We necessarily have  $m\le n$, $\overline r\le m$, $\overline s\le n$ and $\overline r\le \overline s$.

Suppose that $R$ is an algebraic local ring of $K$, $S$ is an algebraic
local ring of $K^*$ such that $R$ and $S$ are regular, $S$ dominates $R$, and $V^*$ dominates 
$S$ (so that $V$ dominates $R$).  

Suppose that $\text{trdeg}_k k(S)=\text{trdeg}_k k(V^*)$.  Further suppose that $t_1,\ldots,t_{\overline \alpha}\in S$ are such that their residues in $k(V^*)$
are a transcendence basis of $k(V^*)$ over $k(V)$. We then have that the residues of $t_1,\ldots
t_{\overline \alpha}$ in $k(S)$ are also a transcendence basis of $k(S)$ over $k(R)$.

We define a monoidal transform sequence (MTS) as in Definition 3.1 of \cite{C1} and
define a uniformizing transform sequence (UTS), a rational uniformizing transform sequence (RUTS)
and a UTS along a valuation as in Definition 3.2 \cite{C1}.

We also define, for our $R$ with quotient field $K$ and extension ring $S$ with quotient field $K^*$ a compatable UTS (CUTS), a compatable RUTS (CRUTS) and a CUTS along $\nu^*$ as on page 29 of
\cite{C1}. Of course, in a CUTS $(R,\overline R''_n,\overline T_n)$ and
$(S,\overline U''_n,\overline U_n)$, we now have that the quotient field of $\overline U_i''$ is a finitely generated extension field of the quotient field
of $\overline T_i''$  for all $i$, as opposed to the much stronger condition of being a finite
extension, which holds in \cite{C1}.

Lemma 3.3 of page 29 of \cite{C1} on the compatability of a CRUTS and its associated MTS is valid in our
extended setting. The same proof holds.

Suppose that $(R,\overline T'',\overline T)$
and $(S,\overline U'',\overline U)$ is a CUTS along $\nu^*$. When there is no
danger of confusion, we will denote by $\nu^*$
our extension of $\nu^*$ to the quotient field of $\overline U''$ which dominates $\overline U''$, 
$\nu$
our extension of $\nu$ to the quotient field of $\overline T''$ which dominates $\overline T''$,
$\tilde \nu^*$
our extension of $\nu^*$ to the quotient field of $\overline U$ which dominates $\overline U$,
and $\tilde \nu$
our extension of $\nu$ to the quotient field of $\overline T$ which dominates $\overline T$.


 For $f\in\overline U$, we will write $\nu^*(f)<\infty$
to mean $\tilde\nu^*(f)\in\Gamma_{\nu^*}$. For $f\in\overline T$, $\nu(f)<\infty$ will mean
$\tilde \nu^*(f)\in\Gamma_{\nu}$.

Let $p_{\overline U}=\{f\in \overline U\mid \nu^*(f)=\infty\}$,
$p_{\overline T}=\{f\in\overline T\mid \nu(f)=\infty\}$.
Our extension of $\nu^*$ to $Q(\overline U/p_{\overline U})$ and of $\nu$ to 
$Q(\overline T/p_{\overline T})$ are canonical and have value groups $\Gamma_{\nu^*}$
and $\Gamma_{\nu}$ respectively. Note that we have  natural embeddings
$\overline T''\subset \overline T/p_{\overline T}$ and $\overline U''\subset \overline U/p_{\overline U}$. We will in general not be concerned with precise values of elements
in $Q(\overline U)$ and $Q(\overline T)$ which have infinite value.

\section{Perron Transforms}\label{PT}
In this section,  assumptions and notations will be as in Section \ref{RR1V}. 

We define a UTS $\overline T\rightarrow\overline T(1)$ of type $I$ and a UTS $\overline T\rightarrow \overline T(1)$ of type $II_r$ along $\nu$, using the ``Algorithm of Perron''
\cite{Z2} as in section 4.1 of \cite{C1}. Since our notations are a little different,
we summarize the final forms of the transformations here.  
We assume (as in section 4.1 of \cite{C1}) that $\overline T''$ has regular parameters
$(\overline x_1,\ldots,\overline x_m)$ such that $\nu(\overline x_1),\ldots,\nu(\overline x_{\overline r})$ is a rational basis of $\Gamma_{\nu}\otimes{\bf Q}$. 

We first state the equations defining a UTS $\overline T\rightarrow \overline T(1)$ of 
type $I$. $\overline T''(1)=\overline T'(1)$ has regular parameters
$(\overline x_1(1),\ldots,\overline x_m(1))$ such that
$$
\begin{array}{ll}
\overline x_1=&\overline x_1(1)^{a_{11}}\cdots\overline x_{\overline r}(1)^{a_{1\overline r}}\\
&\vdots\\
\overline x_{\overline r}=&\overline x_1(1)^{a_{\overline r1}}\cdots\overline x_{\overline r}(1)^{a_{\overline r\overline r}}
\end{array}
$$
and $\overline x_i=\overline x_i(1)$ for $\overline r<i\le m$. The matrix $A=(a_{ij})$
of natural numbers is computed using Perron's algorithm. We have
$\text{Det}(A)=\pm 1$, and $\nu(\overline x_1(1)),\ldots,\nu(\overline x_{\overline r}(1))$
are a rational basis of $\Gamma_{\nu}\otimes{\bf Q}$.

We now state the equations defining a UTS  $\overline T\rightarrow \overline T(1)$ of
type $II_r$ with $0<r\le m-\overline r$. $\overline T''(1)$ has regular parameters $(\overline x_1(1),\ldots,
\overline x_m(1))$ such that 
$$
\begin{array}{ll}
\overline x_1=&\overline x_1(1)^{a_{11}}\cdots\overline x_{\overline r}(1)^{a_{1\overline r}}
c^{a_{1,\overline r+1}}\\
&\vdots\\
\overline x_{\overline r}=&\overline x_1(1)^{a_{\overline r1}}\cdots\overline x_{\overline r}(1)^{a_{\overline r\overline r}}
c^{a_{\overline r,\overline r+1}}\\
\overline x_{\overline r+r}&=\overline x_1(1)^{a_{\overline r+1,1}}\cdots
\overline x_{\overline r}(1)^{a_{\overline r+1,\overline r}}(\overline x_{\overline r+r}(1)+1)
c^{a_{\overline r+1,\overline r+1}}
\end{array}
$$
and $\overline x_i=\overline x_i(1)$ for $\overline r<i\le m$ and $i\ne \overline r+r$.
We have that $c\in k(\overline T(1))$ and $A=(a_{ij})$ is a matrix of natural numbers such that $\text{Det}(A)=\pm 1$. $\nu(\overline x_1(1)),\ldots,\nu(\overline x_{\overline r}(1))$
are a rational basis of $\Gamma_{\nu}\otimes{\bf Q}$. 

We define UTSs $\overline U\rightarrow
\overline U(1)$ along $\nu^*$ in a similar way. Starting with regular parameters $(\overline
y_1,\ldots,\overline y_n)$ in $\overline U''$ such that $\nu^*(\overline y_1)\ldots,
\nu^*(\overline y_{\overline s})$ are a rational basis of $\Gamma_{\nu^*}\otimes{\bf Q}$,
we define a UTS $\overline U\rightarrow \overline U(1)$ of type $I$ so that
$\overline U''(1)=\overline U'(1)$ has regular parameters $(\overline y_1(1),\ldots,
\overline y_n(1))$ such that
$$
\begin{array}{ll}
\overline y_1&=\overline y_1(1)^{b_{11}}\cdots\overline y_{\overline s}(1)^{b_{1\overline s}}\\
&\vdots\\
\overline y_{\overline s}&=\overline y_1(1)^{b_{\overline s1}}\cdots\overline y_{\overline s}(1)^{b_{\overline s\overline s}}
\end{array}
$$
and $\overline y_i=\overline y_i(1)$ for $\overline s<i\le n$. We have that
$B=(b_{ij})$ is a matrix of natural numbers such that $\text{Det}(B)=\pm 1$ and
$\nu^*(\overline y_1(1)),\ldots,\nu^*(\overline y_{\overline s}(1))$ are a rational
basis of $\Gamma_{\nu^*}\otimes{\bf Q}$.

We define a UTS $\overline U\rightarrow \overline U(1)$ of type $II_r$ with
$0<r\le n-\overline s$ so that
$\overline U''(1)$ has regular parameters $(\overline y_1(1),\ldots,\overline y_n(1))$
such that
$$
\begin{array}{ll}
\overline y_1=&\overline y_1(1)^{b_{11}}\cdots\overline y_{\overline s}(1)^{b_{1\overline s}}
d^{b_{1,\overline s+1}}\\
&\vdots\\
\overline y_{\overline s}=&\overline y_1(1)^{b_{\overline s1}}\cdots\overline y_{\overline s}(1)^{b_{\overline s\overline s}}
d^{b_{\overline s,\overline s+1}}\\
\overline y_{\overline s+r}&=\overline y_1(1)^{b_{\overline s+1,1}}\cdots\overline y_{\overline s}(1)^{b_{\overline s+1,\overline s}}(\overline y_{\overline s+r}(1)+d)d^{b_{\overline s+1,
\overline s+1}}
\end{array}
$$
and $\overline y_i=\overline y_i(1)$ for $\overline s<i\le n$ and $i\ne \overline s+r$.
We have $d\in k(\overline U(1))$ and $B=(b_{ij})$ is a matrix of natural numbers such 
that $\text{Det}(B)=\pm1$ and  $\nu^*(\overline y_1(1)),\ldots,\nu^*(\overline y_{\overline s}(1))$ are a rational basis of $\Gamma_{\nu^*}\otimes{\bf Q}$.

\begin{Lemma}\label{Lemma4.3} Suppose that $(R,\overline T'',\overline T)$
and $(S,\overline U'',\overline U)$ is a CUTS along $\nu^*$, $\overline T''$ has regular 
parameters $(\overline x_1,\ldots,\overline x_m)$ and $\overline U''$ has regular
parameters $(\overline y_1,\ldots,\overline y_n)$, related by
$$
\begin{array}{ll}
\overline x_1&=\overline y_1^{c_{11}}\cdots \overline y_{\overline s}^{c_{1\overline s}}\alpha_1\\
&\vdots\\
\overline x_{\overline r}&=\overline y_1^{c_{\overline r1}}\cdots  \overline z_{\overline s}^{c_{\overline r\overline s}}\alpha_{\overline r}
\end{array}
$$
such that $\alpha_1,\ldots,\alpha_{\overline r}\in k(\overline U)$, $\nu(\overline x_1),\ldots,\nu(\overline x_{\overline r})$ are rationally independent, 
$\nu^*(\overline y_1),\ldots,\nu^*(\overline y_{\overline s})$ are rationally independent
and $(c_{ij})$ has rank $\overline r$. Suppose that $\overline T\rightarrow \overline T(1)$ is a 
UTS of type I along $\nu$, such that $\overline T'(1)=\overline T''(1)$ has regular
parameters $(\overline x_1(1),\ldots,\overline x_m(1))$ with

$$
\begin{array}{ll}
\overline x_1&=\overline x_1(1)^{a_{11}}\cdots\overline x_{\overline r}(1)^{a_{1\overline r}}\\
&\vdots\\
\overline x_{\overline r}&=\overline x_1(1)^{a_{\overline r1}}\cdots\overline x_{\overline r}(1)
^{a_{\overline r,\overline r}}.
\end{array}
$$

Then there exists a UTS of type I along $\nu^*$, $\overline U\rightarrow \overline U(1)$
such that $(R,\overline T''(1),\overline T(1))$ and $(S,\overline U''(1),\overline U(1))$
is a CUTS along $\nu^*$ and $\overline U'(1)=\overline U''(1)$ has regular parameters
$(\overline y_1(1),\ldots,\overline y_n(1))$ with
$$
\begin{array}{ll}
\overline y_1&=\overline y_1(1)^{b_{11}}\cdots\overline y_{\overline s}(1)^{b_{1\overline s}}\\
&\vdots\\
\overline y_{\overline s}&=\overline y_1(1)^{b_{\overline s}1}\cdots\overline y_{\overline s}(1)^{b_{\overline s\overline s}}
\end{array}
$$
and
$$
\begin{array}{ll}
\overline x_1(1)&=\overline y_1(1)^{c_{11}(1)}\cdots\overline y_{\overline s}(1)^{c_{1\overline s}(1)}\alpha_1(1)\\
&\vdots\\
\overline x_{\overline r}(1)&=\overline y_1(1)^{c_{\overline r1}(1)}\cdots\overline y_{\overline s}(1)^{c_{\overline r\overline s}(1)}\alpha_{\overline r}(1)
\end{array}
$$
where $\alpha_1(1),\ldots,\alpha_{\overline r}(1)\in k(\overline U(1))$
are products of integral powers of the $\alpha_i$,
$\nu(\overline x_1(1)),\ldots,\nu(\overline x_{\overline r}(1))$ are rationally independent,
$\nu^*(\overline y_1(1)),\ldots,\nu^*(\overline y_{\overline s}(1))$ are rationally
independent and $(c_{ij}(1))$ has rank $\overline r$.
\end{Lemma}

Lemma \ref{Lemma4.3} is a minor extension of Lemma 4.3 \cite{C1}. The same proof is valid,
after replacing $s$ in Lemma 4.3 with $\overline r$ and $\overline s$ as necessary.

\begin{Lemma}\label{Lemma4.4} Suppose that $(R,\overline T'',\overline T)$
and $(S,\overline U'',\overline U)$ is a CUTS along $\nu^*$, $\overline T''$ has regular 
parameters $(\overline x_1,\ldots,\overline x_m)$ and $\overline U''$ has regular
parameters $(\overline y_1,\ldots,\overline y_n)$ with
 
\begin{equation}\label{eq29}
\begin{array}{ll}
\overline x_1&=\overline y_1^{c_{11}}\cdots \overline y_{\overline s}^{c_{1\overline s}}\alpha_1\\
&\vdots\\
\overline x_{\overline r}&=\overline y_1^{c_{\overline r1}}\cdots  \overline z_{\overline s}^{c_{\overline r\overline s}}\alpha_{\overline r}\\
\overline x_{\overline r+1}&=\overline y_{\overline s+1}\\
&\vdots\\
\overline x_{\overline r+l}&=\overline y_{\overline s+l}
\end{array}
\end{equation} 
such that $\alpha_1,\ldots,\alpha_{\overline r}\in k(\overline U)$, $\nu(\overline x_1),\ldots,\nu(\overline x_{\overline r})$ are rationally independent, 
$\nu^*(\overline y_1),\ldots,\nu^*(\overline y_{\overline s})$ are rationally independent
and $(c_{ij})$ has rank $\overline r$. 

Suppose that $\overline T\rightarrow \overline T(1)$ is a 
UTS of type $II_r$ along $\nu$, with $r\le l$, such that $\overline T''(1)$ has regular
parameters $(\overline x_1(1),\ldots,\overline x_m(1))$ with
$$
\begin{array}{ll}
\overline x_1&=\overline x_1(1)^{a_{11}}\cdots\overline x_{\overline r}(1)^{a_{1\overline r}}
c^{a_{1,\overline r+1}}\\
&\vdots\\
\overline x_{\overline r}&=\overline x_1(1)^{a_{\overline r1}}\cdots\overline x_{\overline r}(1)
^{a_{\overline r,\overline r}}c^{a_{\overline r,\overline r+1}}\\
\overline x_{\overline r+r}&=\overline x_1(1)^{a_{\overline r+1,1}}\cdots\overline x_{\overline r}(1)
^{a_{\overline r+1,\overline r}}(\overline x_{\overline r+r}(1)+1)c^{a_{\overline r+1,\overline r+1}}.
\end{array}
$$
Then there exists a UTS of type $II_r$ (followed by a UTS of type $I$)
along $\nu^*$, $\overline U\rightarrow \overline U(1)$,
such that $(R,\overline T''(1),\overline T(1))$ and $(S,\overline U''(1),\overline U(1))$
is a CUTS along $\nu^*$ and $\overline U''(1)$ has regular parameters
$(\overline y_1(1),\ldots,\overline y_n(1))$ with
$$
\begin{array}{ll}
\overline y_1&=\overline y_1(1)^{b_{11}}\cdots\overline y_{\overline s}(1)^{b_{1\overline s}}
d^{b_{1\overline s+1}}\\
&\vdots\\
\overline y_{\overline s}&=\overline y_1(1)^{b_{\overline s}1}\cdots\overline y_{\overline s}(1)^{b_{\overline s\overline s}}d^{b_{\overline s,\overline s+1}}\\
\overline y_{\overline s+r}&=\overline y_1(1)^{b_{\overline s+1,1}}\cdots\overline y_{\overline s}(1)^{b_{\overline s+1,\overline s}}(\overline y_{\overline s+r}(1)+1)d^{b_{\overline s+1,\overline s+1}}
\end{array}
$$
and
$$
\begin{array}{ll}
\overline x_1(1)&=\overline y_1(1)^{c_{11}(1)}\cdots\overline y_{\overline s}(1)^{c_{1\overline s}(1)}\alpha_1(1)\\
&\vdots\\
\overline x_{\overline r}(1)&=\overline y_1(1)^{c_{\overline r1}(1)}\cdots\overline y_{\overline s}(1)^{c_{\overline r\overline s}(1)}\alpha_{\overline r}(1)\\
\overline x_{\overline r+1}(1)&=\overline y_{\overline s+1}\\
&\vdots\\
\overline x_{\overline r+l}&=\overline y_{\overline s+l}
\end{array}
$$
where $\alpha_1(1),\ldots,\alpha_{\overline r}(1)\in k(\overline U(1))$,
$\nu(\overline x_1(1)),\ldots,\nu(\overline x_{\overline r}(1))$ are rationally independent,
$\nu^*(\overline y_1(1)),\ldots,\nu^*(\overline y_{\overline s}(1))$ are rationally
independent and $(c_{ij}(1))$ has rank $\overline r$.
\end{Lemma}

Lemma \ref{Lemma4.4} is a simple variation of  Lemma 4.4 \cite{C1}.

\begin{Lemma}\label{Lemma4.5} Suppose that $(R,\overline T'',\overline T)$
and $(S,\overline U'',\overline U)$ is a CUTS along $\nu^*$, $\overline T''$ has regular 
parameters $(\overline x_1,\ldots,\overline x_m)$ and $\overline U''$ has regular
parameters $(\overline y_1,\ldots,\overline y_n)$ such that
 
\begin{equation}\label{eq38}
\begin{array}{ll}
\overline x_1&=\overline y_1^{c_{11}}\cdots \overline y_{\overline s}^{c_{1\overline s}}\alpha_1\\
&\vdots\\
\overline x_{\overline r}&=\overline y_1^{c_{\overline r1}}\cdots  \overline z_{\overline s}^{c_{\overline r\overline s}}\alpha_{\overline r}\\
\overline x_{\overline r+1}&=\overline y_{\overline s+1}\\
&\vdots\\
\overline x_{\overline r+l}&=\overline y_{\overline s+l}\\
\overline x_{\overline r+l+1}&=\overline y_1^{d_1}\cdots\overline y_{\overline s}^{d_{\overline s}}\overline y_{\overline s+l+1}
\end{array}
\end{equation} 
such that $\alpha_1,\ldots,\alpha_{\overline r}\in k(\overline U)$, $\nu(\overline x_1),\ldots,\nu(\overline x_{\overline r})$ are rationally independent, 
$\nu^*(\overline y_1),\ldots,\nu^*(\overline y_{\overline s})$ are rationally independent
and $(c_{ij})$ has rank $\overline r$. 

Suppose that $\overline T\rightarrow \overline T(1)$ is a 
UTS of type $II_{l+1}$ along $\nu$, such that $\overline T''(1)$ has regular
parameters $(\overline x_1(1),\ldots,\overline x_m(1))$ with
$$
\begin{array}{ll}
\overline x_1&=\overline x_1(1)^{a_{11}}\cdots\overline x_{\overline r}(1)^{a_{1\overline r}}
c^{a_{1,\overline r+1}}\\
&\vdots\\
\overline x_{\overline r}&=\overline x_1(1)^{a_{\overline r1}}\cdots\overline x_{\overline r}(1)
^{a_{\overline r,\overline r}}c^{a_{\overline r,\overline r+1}}\\
\overline x_{\overline r+l+1}&=\overline x_1(1)^{a_{\overline r+1,1}}\cdots\overline x_{\overline r}(1)
^{a_{\overline r,\overline r+1}}(\overline x_{\overline r+l+1}(1)+1)c^{a_{\overline r+1,\overline r+1}}
\end{array}
$$
Then there exists a UTS of type $II_{l+1}$ (followed by a UTS of type $I$)
along $\nu^*$ $\overline U\rightarrow \overline U(1)$
such that $(R,\overline T''(1),\overline T(1))$ and $(S,\overline U''(1),\overline U(1))$
is a CUTS along $\nu^*$ and $\overline U''(1)$ has regular parameters
$(\overline y_1(1),\ldots,\overline y_n(1))$ with
$$
\begin{array}{ll}
\overline y_1&=\overline y_1(1)^{b_{11}}\cdots\overline y_{\overline s}(1)^{b_{1\overline s}}
d^{b_{1\overline s+1}}\\
&\vdots\\
\overline y_{\overline s}&=\overline y_1(1)^{b_{\overline s}1}\cdots\overline y_{\overline s}(1)^{b_{\overline s\overline s}}d^{b_{\overline s,\overline s+1}}\\
\overline y_{\overline s+l+1}&=\overline y_1(1)^{b_{\overline s+1,1}}\cdots\overline y_{\overline s}(1)^{b_{\overline s+1,\overline s}}(\overline y_{\overline s+l+1}(1)+1)d^{b_{\overline s+1,\overline s+1}}
\end{array}
$$
and
$$
\begin{array}{ll}
\overline x_1(1)&=\overline y_1(1)^{c_{11}(1)}\cdots\overline y_{\overline s}(1)^{c_{1\overline s}(1)}\alpha_1(1)\\
&\vdots\\
\overline x_{\overline r}(1)&=\overline y_1(1)^{c_{\overline r1}(1)}\cdots\overline y_{\overline s}(1)^{c_{\overline r\overline s}(1)}\alpha_{\overline r}(1)\\
\overline x_{\overline r+1}(1)&=\overline y_{\overline s+1}\\
&\vdots\\
\overline x_{\overline r+l}&=\overline y_{\overline s+l}\\
\overline x_{\overline r+l+1}&=\overline y_{\overline s+l+1}
\end{array}
$$
where $\alpha_1(1),\ldots,\alpha_{\overline r}(1)\in k(\overline U(1))$,
$\nu(\overline x_1(1)),\ldots,\nu(\overline x_{\overline r}(1))$ are rationally independent,
$\nu^*(\overline y_1(1)),\ldots,\nu^*(\overline y_{\overline s}(1))$ are rationally
independent and $(c_{ij}(1))$ has rank $\overline r$.
\end{Lemma}

Lemma \ref{Lemma4.5} is a simple variation of  Lemma 4.5 \cite{C1}.

\begin{Lemma}\label{Lemma4.6} Suppose that $(R,\overline T'',\overline T)$
and $(S,\overline U'',\overline U)$ is a CUTS along $\nu^*$, $\overline T''$ has regular 
parameters $(\overline x_1,\ldots,\overline x_m)$ and $\overline U''$ has regular
parameters $(\overline y_1,\ldots,\overline y_n)$ such that
$$
\begin{array}{ll}
\overline x_1&=\overline y_1^{c_{11}}\cdots \overline y_{\overline s}^{c_{1\overline s}}\alpha_1\\
&\vdots\\
\overline x_{\overline r}&=\overline y_1^{c_{\overline r1}}\cdots  \overline y_{\overline s}^{c_{\overline r\overline s}}\alpha_{\overline r}\\
\overline x_{\overline r+1}&=\overline y_{\overline s+1}\\
&\vdots\\
\overline x_{\overline r+l}&=\overline y_{\overline s+l}\\
\overline x_{\overline r+l+1}&=\overline y_1^{d_1}\cdots\overline y_{\overline s}^{d_{\overline s}}\delta
\end{array}
$$ 
such that $\delta\in \overline U''$ is a unit, $\alpha_1,\ldots,\alpha_{\overline r}\in k(\overline U)$, $\nu(\overline x_1),\ldots,\nu(\overline x_{\overline r})$ are rationally independent, 
$\nu^*(\overline y_1),\ldots,\nu^*(\overline y_{\overline s})$ are rationally independent
and $(c_{ij})$ has rank $\overline r$. 

Suppose that $\overline T\rightarrow \overline T(1)$ is a 
UTS of type $II_{l+1}$ along $\nu$, such that $\overline T''(1)$ has regular
parameters $(\overline x_1(1),\ldots,\overline x_m(1))$ with
$$
\begin{array}{ll}
\overline x_1&=\overline x_1(1)^{a_{11}}\cdots\overline x_{\overline r}(1)^{a_{1\overline r}}
c^{a_{1,\overline r+1}}\\
&\vdots\\
\overline x_{\overline r}&=\overline x_1(1)^{a_{\overline r1}}\cdots\overline x_{\overline r}(1)
^{a_{\overline r,\overline r}}c^{a_{\overline r,\overline r+1}}\\
\overline x_{\overline r+l+1}&=\overline x_1(1)^{a_{\overline r+1,1}}\cdots\overline x_{\overline r}(1)
^{a_{\overline r,\overline r+1}}(\overline x_{\overline r+l+1}(1)+1)c^{a_{\overline r+1,\overline r+1}}.
\end{array}
$$
Then there exists a UTS of type  $I$
along $\nu^*$ $\overline U\rightarrow \overline U(1)$
such that $(R,\overline T''(1),\overline T(1))$ and $(S,\overline U''(1),\overline U(1))$
is a CUTS along $\nu^*$ and $\overline U'(1)$ has regular parameters
$(\hat y_1(1),\ldots,\hat y_n(1))$ with
$$
\begin{array}{ll}
\overline y_1&=\hat y_1(1)^{b_{11}}\cdots\hat y_{\overline s}(1)^{b_{1\overline s}}\\
&\vdots\\
\overline y_{\overline s}&=\hat y_1(1)^{b_{\overline s}1}\cdots\hat y_{\overline s}(1)^{b_{\overline s\overline s}}\\
\end{array}
$$
and $\overline U''(1)$ has regular parameters $(\overline y_1(1),\ldots,\overline y_n(1))$
such that $\overline y_i(1)=\epsilon_i\hat y_i(1)$ for $1\le i\le \overline s$ for some
units $\epsilon_i\in\overline U''(1)$, 
$$
\begin{array}{ll}
\overline x_1(1)&=\overline y_1(1)^{c_{11}(1)}\cdots\overline y_{\overline s}(1)^{c_{1\overline s}(1)}\alpha_1(1)\\
&\vdots\\
\overline x_{\overline r}(1)&=\overline y_1(1)^{c_{\overline r1}(1)}\cdots\overline y_{\overline s}(1)^{c_{\overline r\overline s}(1)}\alpha_{\overline r}(1)\\
\overline x_{\overline r+1}(1)&=\overline y_{\overline s+1}\\
&\vdots\\
\overline x_{\overline r+l}&=\overline y_{\overline s+l}\\
\end{array}
$$
where $\alpha_1(1),\ldots,\alpha_{\overline r}(1)\in k(\overline U(1))$,
$\nu(\overline x_1(1)),\ldots,\nu(\overline x_{\overline r}(1))$ are rationally independent,
$\nu^*(\overline y_1(1)),\ldots,\nu^*(\overline y_{\overline s}(1))$ are rationally
independent and $(c_{ij}(1))$ has rank $\overline r$.
\end{Lemma}

Lemma \ref{Lemma4.6} is a simple variation of  Lemma 4.6 \cite{C1}.
The last two lines of page 45  \cite{C1} must be replaced with the following lines:

After possibly interchanging $\hat y_1(1),\ldots, \hat y_{\overline s}(1)$, we may
assume that if
$$
\tilde C=\left(\begin{array}{lll} c_{11}(1)&\cdots&c_{1\overline r}(1)\\
\vdots&&\vdots\\
c_{\overline r1}(1)&\cdots&c_{\overline r\overline r}(1)
\end{array}\right)
$$
then $\text{Det}(\tilde C)\ne 0$. Set $(e_{ij})=\tilde C^{-1}$, and set
$$
\epsilon_i=\left\{
\begin{array}{ll}
\delta_1^{e_{i1}}\cdots\delta_{\overline r}^{e_{i\overline r}}
(\frac{\delta_{l+1}}{c})
^{-\gamma_1e_{i1}-\cdots-\gamma_{\overline r}e_{i\overline r}}
&
\text{for }1\le i\le \overline r,\\
1& \text{ for }\overline r <i\le\overline s
\end{array}\right.
$$

\section{UTSs of form $\overline m$}\label{UTS}

In this section, assumptions and notations will be as in Section \ref{RR1V}.

\begin{Definition}\label{Definition1} 
Suppose that $(R,\overline T'',\overline T)$
 is a UTS along $\nu$, such that $\overline T''$
contains a subfield isomorphic to $k(c_0)$ for some $c_0\in \overline T''$
and 
$\overline T''$ has regular 
parameters $(\overline z_1,\ldots,\overline z_m)$

Suppose that $0\le \overline m\le m-\overline r$.
 We will say that a UTS along $\nu$

\begin{equation}\label{eq52}
\overline T\rightarrow \overline T(1)\rightarrow\cdots\rightarrow\overline T(t)
\end{equation}
is of form $\overline m$ if for $0\le\alpha\le t$,
$\overline T''(\alpha)$ has regular parameters
$$
(\overline z_1(\alpha),\ldots,\overline z_m(\alpha))\text{ and }(\tilde{\overline z}'_1(\alpha),
\ldots,\tilde{\overline z}'_m(\alpha)),
$$
where $\overline z_i(0)=\overline z_i$ for $1\le i\le m$, 
$\overline T''(\alpha)$ contains a subfield isomorphic to $k(c_0,\ldots,c_{\alpha})$
and there are polynomials 
$$
P_{i,\alpha}\in k(c_0,\ldots,c_{\alpha})[\overline z_1(\alpha),\ldots,\overline z_{i-1}(\alpha)]\text{ for }\overline r+1\le i\le\overline r+\overline m
$$
 such that
$$
\tilde{\overline z}'_i(\alpha)=\left\{
\begin{array}{ll}
\overline z_i(\alpha)-P_{i,\alpha}&\text{if }\overline r+1\le i\le \overline r+\overline m\\
\overline z_i(\alpha)&\text{otherwise}
\end{array}\right.
$$

Each $\overline T(\alpha)\rightarrow \overline T(\alpha+1)$ will either be of type $I$ or
of type $II_r$ with $1\le r\le \overline r+\overline m$.

In a transformation $\overline T(\alpha)\rightarrow\overline T(\alpha+1)$ of type $I$,
$\overline T''(\alpha+1)$ will have regular parameters $(\overline z_1(\alpha+1),\ldots,
\overline z_m(\alpha+1))$ defined by 
$$
\begin{array}{ll}
\tilde{\overline z}'_1(\alpha)&=\overline z_1(\alpha+1)^{a_{11}(\alpha+1)}\cdots
\overline z_{\overline r}(\alpha+1)^{a_{1\overline r}(\alpha+1)}\\
&\vdots\\
\tilde{\overline z}'_{\overline r}(\alpha)&=\overline z_1(\alpha+1)^{a_{\overline r1}(\alpha+1)}
\cdots\overline z_{\overline r}(\alpha+1)^{a_{\overline r\overline r}(\alpha+1)}.
\end{array}
$$
and $c_{\alpha+1}$ is defined to be 1. In a transformation $\overline T(\alpha)\rightarrow
\overline T(\alpha+1)$ of type $II_r$ ($1\le r\le \overline r+\overline m$) $\overline T''(\alpha+1)$
will have regular parameters 
$(\overline z_1(\alpha+1),\ldots,
\overline z_m(\alpha+1))$ defined by 
$$
\begin{array}{ll}
\tilde{\overline z}'_1(\alpha)&=\overline z_1(\alpha+1)^{a_{11}(\alpha+1)}\cdots
\overline z_{\overline r}(\alpha+1)^{a_{1\overline r}(\alpha+1)}
c_{\alpha+1}^{a_{1\overline r+1}(\alpha+1)}\\
&\vdots\\
\tilde{\overline z}'_{\overline r}(\alpha)&=\overline z_1(\alpha+1)^{a_{\overline r1}(\alpha+1)}
\cdots\overline z_{\overline r}(\alpha+1)^{a_{\overline r\overline r}(\alpha+1)}
c_{\alpha+1}^{a_{\overline r,\overline r+1}(\alpha+1)}\\
\tilde{\overline z}'_{\overline s+r}(\alpha)&=
\overline z_1(\alpha+1)^{a_{\overline r+1,1}(\alpha+1)}\cdots\overline z_{\overline r}(\alpha+1)
^{a_{\overline r+1,\overline r}(\alpha+1)}\\
&\,\,\,\,\,\,\,\,\,\,\cdot
(\overline z_{\overline r+r}(\alpha+1)+1)c_{\alpha+1}^{a_{\overline r+1,
\overline r+1}(\alpha+1)}
\end{array}
$$
\end{Definition}

\begin{Theorem}\label{Theorem3}
Suppose that $(R,\overline T'',\overline T)$
 is a UTS along $\nu$, such that $\overline T''$
contains a subfield isomorphic to $k(c_0)$ for some $c_0\in \overline T''$
and 
$\overline T''$ has regular 
parameters $(\overline z_1,\ldots,\overline z_m)$.
\begin{description} 
\item[(A1)]
Suppose that $0\le\overline m\le m-\overline r$. Then there exists a UTS (\ref{eq52})
of form $\overline m$
 such that 
$$
p_{\overline m}(i)=\{f\in k(\overline T(i))[[\overline z_1(i),\ldots,\overline z_{\overline r+\overline m}(i)]]\mid \nu(f)=\infty\}
$$
has the form 
\begin{equation}\label{eq13}
p_{\overline m}(t)=
(\overline z_{r(1)}(t)-Q_{r(1)}(\overline z_1(t),\ldots,\overline z_{r(1)-1}),\ldots,
\overline z_{r(\tilde m)}(t)-Q_{r(\tilde m)}(\overline z_1(t),\ldots,\overline z_{r(\tilde m)-1}))
\end{equation}
for some $0\le\tilde m\le\overline m$ and $\overline r<r(1)<r(2)<\cdots<r(\tilde m)\le\overline r+\overline m$,
where $Q_{r(i)}$ are power series with coefficients in $k(c_0,\ldots,c_t)$
\item[(A2)] Suppose that $L$ is a finite extension field of $k(c_0)(t_1,\ldots,t_{\beta})$,
with $0\le\beta\le\overline\alpha$ (For the definition of $t_1,\ldots t_{\overline\alpha}$
see section \ref{RR1V}). Suppose that $\nu'$ is an extension of $\tilde \nu$ to 
$Q(L[[\overline z_1,\ldots,\overline z_{\overline r+\overline m}]])$ such that $\nu'$
dominates $L[[\overline z_1,\ldots,\overline z_{\overline r+\overline m}]]$. $h\in L[[\overline z_1,\ldots,\overline z_{\overline r+\overline m}]]$
for some $0\le \overline m$ with $\overline r+\overline m\le m$ and $\nu(h)<\infty$.
Then there exists a UTS (\ref{eq52}), of form $\overline m$, such that (\ref{eq13}) holds in $\overline T(t)$ and
$$
h=\overline z_1(t)^{d_1}\cdots\overline z_{\overline r}(t)^{d_{\overline r}}u
$$
where $u\in L(c_1,\ldots,c_t)[[\overline z_1(t),\ldots,\overline z_{\overline r+\overline m}(t)]]$ is a unit power series. 

If $h\in L[\overline z_1,\ldots,\overline z_{\overline r+\overline m}]$, then
$u\in L(c_1,\ldots,c_t)[\overline z_1(t),\ldots,\overline z_{\overline r+\overline m}(t)]$.
\item[(A3)]
Suppose that $L$ is a finite extension field of 
$k(c_0)(t_1,\ldots,t_{\beta})$ with $0\le\beta\le\overline\alpha$
(For the definition of $t_1,\ldots t_{\overline\alpha}$
see section \ref{RR1V}). Suppose that $\nu'$ is an extension of $\tilde \nu$ to 
$Q(L[[\overline z_1,\ldots,\overline z_{\overline r+\overline m}]])$ such that
$\nu'$ dominates $L[[\overline z_1,\ldots,\overline z_{\overline r+\overline m}]]$.
$h\in L[[\overline z_1,\ldots,\overline z_{\overline r+\overline m}]]$
for some $0<\overline m$ with $\overline r+\overline m\le m$, $\nu(h)=\infty$
and $A>0$ is given.
Then there exists a UTS (\ref{eq52}), of form $\overline m$, such that (\ref{eq13}) holds in $\overline T(t)$ and
$$
h=\overline z_1(t)^{d_1}\cdots\overline z_{\overline r}(t)^{d_{\overline r}}\Sigma
$$
where 
$\nu(\overline z_1(t)^{d_1}\cdots\overline z_{\overline r}(t)^{d_{\overline r}})>A$ and
$\Sigma\in L(c_1,\ldots,c_t)[[\overline z_1(t),\ldots,\overline z_{\overline r+\overline m}]]$.

If $h\in L[\overline z_1,\ldots,\overline z_{\overline r+\overline m}]$, then
$\Sigma\in L(c_1,\ldots,c_t)[\overline z_1(t),\ldots,\overline z_{\overline r+\overline m}(t)]$.

\end{description}
\end{Theorem}
\begin{pf}
(A1), (A2), (A3) replace (53), (54) and (55) of the proof of Theorem 4.7 \cite{C1}.
We observe that (A1) is trivial for $\overline m=0$ since $p_{0}=(0)$.
The proof of (A2) when $\overline m=0$ follows from the ``Proof of (54) for $m=s$''
on page 50 of the proof of Theorem 4.7 \cite{C1}.

We will now establish (A1), (A2) and (A3) by proving the following inductive 
statements.

$A(m')$: (A1), (A2) and (A3) for $\overline m\le m'$ imply (A1) for $\overline m=m'$.

$B(m')$: (A2), (A3) for $\overline m<m'$ and (A1) for $\overline m=m'$ imply
(A2) and (A3) for $\overline m=m'$.

The ``Proof of $A(\overline m)$ $(s<\overline m)$'' on pages 51-55 of \cite{C1}
proves $A(m')$ for $0<m'$.

We now give the proof of $B(m')$ for $m'>0$.
By assumption, there exists a UTS $\overline T\rightarrow \overline T(t)$
satisfying (A1) for $\overline m=m'$. After replacing $\overline T''$
with
$\overline T''(t)$
and replacing $c_0$ with a primitive element of $k(c_0,\ldots,c_t)$ over $k$, we may assume that 
$$
p_{m'}=
(\overline z_{r(1)}-Q_{r(1)}(\overline z_1,\ldots,\overline z_{r(1)-1}),\ldots,
\overline z_{r(\tilde m')}-Q_{r(\tilde m')}(\overline z_1,\ldots,\overline z_{r(\tilde m')-1}))
$$
where the $Q_{r(i)}$ are power series with coefficients in $k(c_0)$.

Let $\overline R=L[[\overline z_1,\ldots,\overline z_{\overline r+m'}]]$. Let 
$$
\overline p = \{f\in \overline R\mid \nu(f)=\infty\}.
$$
Let
$$
\overline p_{m'}=p_{m'}\cap k(c_0)[[\overline z_1,\ldots,\overline z_{\overline r+m'}]].
$$
We first establish the following formula:

\begin{equation}\label{eq6} \overline p=\overline p_{m'}\overline R.
\end{equation}

We first prove the identity (\ref{eq6}) when $L=k(c_0)(t_1,\ldots,t_{\beta})$
for $\beta$ with $0\le\beta\le\overline\alpha$.
Let
$$
\tilde R=k(c_0)[[\overline z_1,\ldots,\overline z_{\overline r+ m'}]]/\overline p_{ m'}.
$$

Let $a\in \overline p$.
Suppose that $N>0$ is given. Chevalley's Theorem (Theorem 13, Section 5, Chapter VIII \cite{ZS}) implies there exists $M$ such that
$g\in\tilde R$ and $\nu(g)>M$ implies $g\in m(\tilde R)^N$. There also exists
$N_0\ge N$ such that $\Omega\in m(\overline R)^{N_0}$ implies $\nu(\Omega)>M$.
Write $a=H+\Omega$ with $H\in k(c_0)(t_1,\ldots,t_{\beta})[\overline z_1,\ldots,
\overline z_{\overline r+ m'}]$ and $\Omega\in m(\overline R)^{N_0}$. Then $\nu(H)=\nu(\Omega)>M$.
There exists $0\ne h\in k(c_0)[t_1,\ldots,t_{\beta}]$ such that
$$
hH=\sum_{I=(i_1,\ldots,i_{\beta})}\alpha_It_1^{i_1}\cdots t_{\beta}^{i_{\beta}}
$$
is a polynomial with all $\alpha_I\in k(c_0)[\overline z_1,\ldots,\overline z_{\overline r+m'}]$. Thus $\nu(\alpha_I)>M$ for all $I$ (by Lemma
\ref{Lemma4}). Let $\tilde \alpha_I$ be the residue of $\alpha_I$ in $\tilde R$. $\tilde\alpha_I\in
m(\tilde R)^N$ for all $I$ implies 
$$
\alpha_I\in \overline p_{m'}+m(k(c_0)[[\overline z_1,\ldots,\overline z_{\overline r_m'}]])^N
$$
 for
all $I$ so that $a\in \overline p_{m'}\overline R+m(\overline R)^N$. Since this is true for all
$N$,
$$
a\in \cap_{N>0}\left(\overline p_{m'}\overline R+m(\overline R)^N\right)=\overline p_{m'}\overline R.
$$
(\ref{eq6}) is thus established when $L=k(c_0)(t_1,\ldots,t_{\beta})$.

Now suppose that $L$ is a finite extension of $k(c_0)(t_1,\ldots,t_{\beta})$. Let $\overline M=k(c_0)(t_1,\ldots,t_{\beta})$
and let $M'$ be a Galois closure of $L$ over $\overline M$. Let $G$ be the Galois group of
$M'$ over $\overline M$.  Suppose that $a\in \overline p$. Set
$$
g=\prod_{\sigma\in G}\sigma(a)\in \overline M[[\overline z_1,\ldots,\overline z_{\overline r+ m'}]].
$$
$\nu(g)=\infty$ since $a\mid g$ in $L[[\overline z_1,\ldots,\overline z_{\overline r+m'}]]$. Thus $g\in \overline p_{m'}M'[[\overline z_1,\ldots,\overline z_{\overline r+ m'}]]$ which is a prime ideal invariant under $G$. Thus there exists
$\sigma\in G$ such that $\sigma(a)\in \overline p_{m'}M'[[\overline z_1,\ldots,\overline z_{\overline r+m'}]]$. Necessarily we then have that 
$$
a\in \left(p_{m'}M'[[\overline z_1,\ldots,\overline z_{\overline r+m'}]]\right)\cap
L[[\overline z_1,\ldots,\overline z_{\overline r+m'}]]=\overline p_{m'}L[[\overline z_1,\ldots,\overline z_{\overline r+m'}]].
$$
 (\ref{eq6}) is now established.

Now suppose that $h\in L[[\overline z_1,\ldots,\overline z_{\overline r+m'}]]$ and $\nu(h)<\infty$. Let $\overline M=k(c_0)(t_1,\ldots,t_{\beta})$
and let $M'$ be a Galois closure of $L$ over $\overline M$. Let $G$ be the Galois group of
$M'$ over $\overline M$.   Set
$$
g=\prod_{\sigma\in G}\sigma(h)\in \overline M[[\overline z_1,\ldots,\overline z_{\overline r+ m'}]].
$$
It follows from (\ref{eq6})  that $\nu(g)<\infty$. We will construct a UTS (\ref{eq52}) so that
$$
g=u\overline z_1(t)^{e_1}\cdots\overline z_{\overline r}^{e_{\overline r}}(t)
$$
where $u$ is a unit power series in $k(c_0,\ldots,c_t,t_1,\ldots,t_{\beta})[[\overline z_1(t),
\ldots,\overline z_{\overline r+m'}(t)]]$ and $h\in L(c_1,\ldots,c_t)[[\overline z_1(t),
\ldots,\overline z_{\overline r+m'}(t)]]$. Since $h\mid g$ in $L(c_1,\ldots,c_t)[[\overline z_1(t),
\ldots,\overline z_{\overline r+m'}(t)]]$, we will have that $h$ has the desired form in
$L(c_1,\ldots,c_t)[[\overline z_1(t),
\ldots,\overline z_{\overline r+m'}(t)]]$.

We now follow the argument of the proof of Theorem 4.7 \cite{C1} 
as modified to fit our notation,
from 
``Set $g=\overline z_1^{d_1}\cdots\overline z_{\overline r}^{d_{\overline r}}g_0$''
on the 17th line of the ``Proof of $B(\overline m)$'' of page 55 of \cite{C1} to
the fourth line from the last on page 56, ending with ``$\nu(\sigma_{d-1})=\nu(\overline z_{\overline r+m'})$''.  We must substitute 
$k(c_0,\ldots,c_{\alpha},t_1,\ldots,t_{\beta})$ for $k(c_0,\ldots,c_{\alpha})$,
$\overline r+m'$ for $\overline m$ and $\overline r$ for $s$  in the proof, and for
``(54), (55) for $m<\overline m$'' on line 5 of page 56 \cite{C1} we must substitute
``(A2), (A3) for $m<m'$''.

We now argue as follows.  $\sigma_d$ is a unit. Let
$$
\overline R=k(c_0,\ldots,c_{\alpha+2})(t_1,\ldots,t_{\beta})[[\overline z_1(\alpha+2),\ldots,\overline z_{r+m'}(\alpha+2)]].
$$
 Let $\tilde\alpha_d$ be the residue of $\sigma_d=\overline u_d$
in $k(c_0)(t_1,\ldots,t_{\beta})$, $\tilde\alpha_{d-1}$ be the residue of $\overline u_{d-1}$ in $k(c_0,\ldots,c_{\alpha+2})(t_1,\ldots,t_{\beta})$. 
After replacing $g_0$ with $\frac{1}{\tilde\alpha_d}g_0$, we may assume that $\tilde\alpha_d=1$.
We have 
\begin{equation}\label{eq7}
\lambda_dr\tilde\alpha_d+\lambda_{d-1}\tilde\alpha_{d-1}=0,
\end{equation}
by the argument of page 53 of \cite{C1}, where $\lambda_d,\lambda_{d-1}\in k(c_0,\ldots,c_{\alpha+2})$. Set $\tau=\nu(\sigma_{d-1})=\nu(\overline z_{\overline r+m'})<\infty$. 
Write
$$
\sigma_{d-1}= \sum_{\nu(\overline z_1^{i_1}\cdots\overline z_{\overline r+ m'-1}
^{i_{\overline r+ m'-1}})\le\tau}
g_I\overline z_1^{i_1}\cdots\overline z_{\overline r+ m'-1}
^{i_{\overline r+ m'-1}}
+\sum_{\nu(\overline z_1^{i_1}\cdots\overline z_{\overline r+ m'-1}
^{i_{\overline r+ m'-1}})>\tau}
g_I\overline z_1^{i_1}\cdots\overline z_{\overline r+ m'-1}
^{i_{\overline r+ m'-1}}
$$
with $g_I\in k(c_0)(t_1,\ldots,t_{\beta})$ for all 
$I=(i_1,\ldots,i_{\overline r+m'-1})$.
Set
$$
\Omega=\sum_{\nu(\overline z_1^{i_1}\cdots\overline z_{\overline r+m'-1}
^{i_{\overline r+ m'-1}})\le\tau}
g_I\overline z_1^{i_1}\cdots\overline z_{\overline r+ m'-1}
^{i_{\overline r+ m'-1}}.
$$
We have $\nu(\Omega)=\nu(\sigma_{d-1})$. There exists $0\ne\overline h\in k(c_0)(t_1,\ldots,t_{\beta})$ such that $\overline hg_I\in k(c_0)[t_1,\ldots,t_{\beta}]$
for all $g_I$ in the finite sum $\Omega$.  Thus we have
$$
\overline h\Omega=\sum_{J=(j_1,\ldots,j_{\beta})}\Psi_{J}t_1^{j_1}\cdots t_{\beta}^{j_{\beta}}
$$
with all $\Psi_{J}\in k(c_0)[\overline z_1,\ldots,\overline z_{\overline r+ m'-1}]$.
By Lemma \ref{Lemma4},
$$
\nu(\Omega)=\text{min}\{\nu(\Psi_J)\}.
$$
$\frac{\sigma_{d-1}}{\overline z_{\overline r+m'}}$
has residue $-r\tilde\alpha_d=-r$ in $k(V)(t_1,\ldots,t_{\beta})\subset k(V^*)$ (by the argument of page 53 of \cite{C1}).
Let
$$
\left[\frac{\Psi_J}{\overline z_{\overline r+ m'}}\right]
$$
be the residue of $\frac{\Psi_J}{\overline z_{\overline r+ m'}}$ in $k(V)$.

In $k(V(t_1,\ldots,t_{\beta}))$, we have
$$
-r\overline h=\sum_J\left[\frac{\Psi_J}{\overline z_{\overline r+ m'}}\right]
t_1^{j_1}\cdots t_{\beta}^{j_{\beta}}.
$$

Since $0\ne \overline h\in k(c_0)(t_1,\ldots,t_{\beta})$ and the $t_1^{j_1}\cdots t_{\beta}^{j_{\beta}}$ are linearly independent over
$k(V)$, we have that 
$$
\left[\frac{\Psi_J}{\overline z_{\overline r+m'}}\right]\in k(c_0)\text{ for all }J
$$
and 
$$
\left[\frac{\Psi_{J_0}}{\overline z_{\overline r+m'}}\right]\ne 0\text{ for some }J_0.
$$
Let $c=\left[\frac{\overline z_{\overline r+m'}}{\Psi_{J_0}}\right]$. Set
$$
\overline z_{\overline r+ m'}'=\overline z_{\overline r+ m'}-c\Psi_{J_0}
\in k(c_0)[\overline z_1,\ldots,\overline z_{\overline r+m'}].
$$
We have 
$$
\nu(\overline z_{\overline r+m'}')>\nu(\overline z_{\overline r+m'}).
$$
We further have
$$
\nu(g_0)\ge \nu(\overline z_{\overline r+m'}^r)\ge \nu(\overline z_{\overline r+m'})
$$
since $\overline z_{\overline r+ m'}^r$ is a minimal value term of $g_0$.

Now we finish the proof as on lines 1 - 16 of page 57 of \cite{C1}. On line 7 of page 57 we must replace
``$P_{\overline m,0}\in k(c_0)[\overline z_1,\ldots,\overline z_{\overline m-1}]$''
with ``$P_{m',0}\in k(c_0)[\overline z_1,\ldots,\overline z_{\overline r+ m'-1}]$''. on line 11 we must replace ``$\overline u$ is a unit power series with
coefficients in $k(c_0,\ldots,c_t)$'' with ``$\overline u$ is a unit power series with
coefficients in $k(c_0,\ldots,c_t)(t_1,\ldots,t_{\beta})$''.
(53) on line 12 must be replaced with (A1).
This concludes the proof of Case 1, $\nu(h)<\infty$ of the proof of $B(m')$.

The proof of  $B(m')$, when $\nu(h)=\infty$
 is only a slight modification of the proof of case 2 on page 57 of \cite{C1}.
We must replace  (53) on line 17 with (A1) and replace
``$\sigma_i\in k(\overline T)[[\overline z_1,\ldots,\overline z_{m}]]$'' with
``$\sigma_i\in L[[\overline z_1,\ldots,\overline z_{\overline r+m'}]]$''
on line 21 of page 57. 

The final statements that 
$h\in L[\overline z_1,\ldots,\overline z_{\overline r+\overline m}]$ imply
$h\in L(c_1,\ldots,c_t)[\overline z_1(t),\ldots,\overline z_{\overline r+\overline m}(t)]$
follow since $\overline z_1,\ldots,\overline z_{\overline r+\overline m}\in 
k(c_1,\ldots,c_t)[\overline z_1(t),\ldots,\overline z_{\overline r+\overline m}(t)]$
by the definition of a UTS of form $\overline m$.
\end{pf}

\begin{Lemma}\label{Lemma17}
Suppose that
$$
\overline T\rightarrow \overline T(1)\rightarrow \cdots\rightarrow \overline T(t)
$$
is a UTS of form $\overline m$ as in (\ref{eq52}) of Definition \ref{Definition1}.
Let 
$$
A_i=k(\overline T(i))[[\overline z_1(i),\ldots,\overline z_{\overline r+\overline m}(i)]],
$$
and let
$$
\sigma(i)=\text{ dim }A_i/p_{\overline m}(i)
$$
where $p_{\overline m}(i)$ is defined by (A1) of Theorem \ref{Theorem3}. Then
$$
\sigma(i+1)\le \sigma(i)
$$
for $0\le i\le t-1$.
\end{Lemma}

\begin{pf} There exists an ideal $q\subset A_i$, $0\ne \lambda\in q$ and a maximal ideal
$n$ in $A_i[\frac{q}{n}]$ such that 
$$
A_{i+1}=(A_i[\frac{q}{\lambda}]_n)\sphat.
$$
Let
$$
\tilde p =\cup_{j=1}^{\infty}\left(p_{\overline m}(i)A_i[\frac{q}{\lambda}]_n:
q^jA_i[\frac{q}{\lambda}]_n\right),
$$
the strict transform of $p_{\overline m}(i)$ in $A_i[\frac{q}{\lambda}]_n$.
$q\not\subset p_{\overline m}(i)$ implies $\nu(f)=\infty$ if $f\in \tilde p$.
Thus $\tilde p\subset p_{\overline m}(i+1)$.
$$
A_i/p_{\overline m}(i)\rightarrow A_i[\frac{q}{\lambda}]_{n}/\tilde p
$$
is birational (c.f.  \cite{H1}) and the residue field extension is finite, so
by the dimension formula (c.f. Theorem 15.6 \cite{Ma})
$\sigma(i)=\text{ dim }A_i[\frac{q}{\lambda}]_n/\tilde p$. Thus
$\sigma(i)=\text{ dim }A_{i+1}/\tilde pA_{i+1}$ since completion is flat (c.f. Theorem
8.14 \cite{Ma}) and by Theorem 15.1 \cite{Ma}. We thus have
$$
\sigma(i)\ge\text{ dim }A_{i+1}/p_{\overline m}(i+1)=\sigma(i+1).
$$
\end{pf}

\section{CUTS of form $\overline m$}\label{CUTS}
Let assumptions and notations be as in Section \ref{RR1V} throughout this
section.
Suppose that $(R,\overline T'',\overline T)$
and $(S,\overline U'',\overline U)$ is a CUTS along $\nu^*$, $\overline T''$ has regular 
parameters $(\overline x_1,\ldots,\overline x_m)$ and $\overline U''$ has regular
parameters $(\overline y_1,\ldots,\overline y_n)$ with 
\begin{equation}\label{eq1}
\begin{array}{ll}
\overline x_1&=\overline y_1^{c_{11}}\cdots \overline y_{\overline s}^{c_{1\overline s}}\phi_1\\
&\vdots\\
\overline x_{\overline r}&=\overline y_1^{c_{\overline r1}}\cdots  \overline y_{\overline s}^{c_{\overline r\overline s}}\phi_{\overline r}\\
\overline x_{\overline r+1}&=\overline y_{\overline s+1}\\
&\vdots\\
\overline x_{\overline r+l}&=\overline y_{\overline s+l}\\
\end{array}
\end{equation} 
such that  $\phi_1,\ldots,\phi_{\overline r}\in k(\overline U)$, $\nu(\overline x_1),\ldots,\nu(\overline x_{\overline r})$ are rationally independent, 
\linebreak
$\nu^*(\overline y_1),\ldots,\nu^*(\overline y_{\overline s})$ are rationally independent
and $(c_{ij})$ has rank $\overline r$.

Let $C=(C_1,\ldots,C_{\overline s})$ be the $\overline r\times\overline s$ matrix 
$(c_{ij})$ of
(\ref{eq1}). Multiplication by $C$ defines a linear map
$\Phi:{\bf Q}^{\overline r}\rightarrow {\bf Q}^{\overline s}$, $\Phi(v)=vC$.
$\Phi$ is 1-1 since $C$ has rank $\overline r$.

Suppose that we have a CUTS as in (\ref{eq1}), and that $f\in k(\overline U)[[\overline y_1,
\ldots,\overline y_{\overline s+l}]]$.

For $\Lambda\in{\bf Z}^{\overline s}$, let $[\Lambda]$ denote the class of $\Lambda$
in ${\bf Z}^{\overline s}/({\bf Q}^{\overline r}C)\cap{\bf Z}^{\overline s}$.
$f$ has a unique expression 
\begin{equation}\label{eq2}
f=\sum_{[\Lambda]\in{\bf Z}^{\overline s}/({\bf Q}^{\overline r}C)\cap{\bf Z}^{\overline s}}
h_{[\Lambda]}
\end{equation}
where 
\begin{equation}\label{eq15}
h_{[\Lambda]}=\sum_{\alpha\in{\bf N}^{\overline s}\mid [\alpha]=[\Lambda]}
g_{\alpha}\overline y_1^{\alpha_1}\cdots\overline y_{\overline s}^{\alpha_{\overline s}}
\end{equation}
and $g_{\alpha}\in k(\overline U)[[\overline y_{\overline s+1},\ldots,\overline y_{\overline s+l}]]$.

Set $G=\Phi^{-1}({\bf Z}^{\overline s})$.
For $\Lambda=(\lambda_1,\ldots,\lambda_{\overline s})\in {\bf N}^{\overline s}$, 
define
$$
P_{\Lambda}=\{v\in{\bf Q}^{\overline r}\mid vC_i+\lambda_i\ge 0\text{ for }
1\le i\le \overline s\}.
$$

For $\Lambda\in{\bf N}^{\overline s}$, we have 

\begin{equation}\label{eq3}
h_{[\Lambda]}=\overline y_1^{\lambda_1}\cdots\overline y_{\overline s}^{\lambda_{\overline s}}
\left[\sum_{v=(v_1,\ldots,v_{\overline r})\in G\cap P_{\Lambda}}
\phi_1^{-v_1}\cdots\phi_{\overline r}^{-v_{\overline r}}\overline x_1^{v_1}\cdots\overline x_{\overline r}^{v_{\overline r}}g_{v}\right]
\end{equation}
where each $g_v\in k(\overline U)[[\overline x_{\overline r+1},\ldots,\overline x_{\overline r+l}]]$.
Here we have reindexed the $g_{\alpha}=g_{\Phi(v)+\Lambda}$ in (\ref{eq15}) as $g_v$.
Let 
$$
H=\{v\in{\bf Z}^{\overline r}\mid vC_i\ge 0\text{ for }1\le i\le\overline s\},
$$
$$
I=\{v\in G\mid vC_i\ge 0\text{ for }1\le i\le {\overline s}\}
$$
and for $\Lambda=(\lambda_1,\ldots,\lambda_{\overline s})\in {\bf N}^{\overline s}$,
$$
M_{\Lambda}=\{v\in G\mid vC_i+\lambda_i\ge0\text{ for }1\le i\le\overline s\}.
$$

$P_{\Lambda}$ is a rational polyhedral set in ${\bf Q}^{\overline r}$ whose associated cone
is 
$$
\sigma =\{v\in{\bf Q}^{\overline r}\mid v C_i=0\text{ for }1\le i\le \overline s\}=\{0\}.
$$
Let $W={\bf Q}^{\overline r}$. $G$ is a lattice in $W$.
Thus $P_{\Lambda}$ is strongly convex and $M_{\Lambda}=P_{\Lambda}\cap G$
is a finitely generated module over the semigroup $I$ (c.f. Theorem 7.1 \cite{CHR}).
Let $\overline n=[G:{\bf Z}^{\overline r}]$. We have $\overline nx\in H$ for all $x\in I$.
Gordon's Lemma (c.f. proposition 1, page 12 \cite{Fu}) implies that $H$ and $I$ are
finitely generated semigroups. There exist $w_1,\ldots,w_{\overline m}\in I$ which generate $I$
as a semi-group.  Then the finite set 
$$
\{a_1w_1+\cdots+a_{\overline m}w_{\overline m}\mid a_i\in{\bf N}\text{ and }
0\le a_i<\overline n\text{ for }1\le i\le \overline m\}
$$
generate $I$ as an $H$ module. We have then that $M_{\Lambda}$  is a finitely generated
module over the semigroup $H$.
Thus there exist
$\overline v_1,\ldots. \overline v_{\overline a}\in H$ and $\overline u_1,\ldots,\overline
u_{\overline b}\in M_{\Lambda}$ such that if $v=(v_1,\ldots,v_{\overline r})\in M_{\Lambda}=G\cap P_{\Lambda}$, then
$$
v=\overline u_i+\sum_{j=1}^{\overline a}n_j\overline v_j
$$
for some $1\le i\le \overline b$ and $n_1,\ldots,n_{\overline a}\in{\bf N}$. Thus,
$$
\overline x_1^{v_1}\cdots\overline x_{\overline r}^{v_{\overline r}}
=\overline x_1^{\overline u_{i,1}}\cdots \overline x_{\overline r}^{\overline u_{i,\overline r}}
\prod_{j=1}^{\overline a}(\overline x_1^{\overline v_{j,1}}\cdots
\overline x_{\overline r}^{\overline v_{j,\overline r}})^{n_j}
$$
where $\overline u_i=(\overline u_{i,1},\ldots,\overline u_{i,\overline r})$ and 
$\overline v_j=(\overline v_{j,1},\ldots,\overline v_{j,\overline r})$ for $1\le j\le \overline a$. Since each $\overline v_j\in H$, $\nu(\overline x_1^{\overline v_{j,1}}\cdots
\overline x_{\overline r}^{\overline v_{j,\overline r}})\ge 0$. Thus there exists 
by Lemma 4.2 \cite{C1} and Lemma \ref{Lemma4.3} a CUTS of
type I along $\nu^*$ 
\begin{equation}\label{eq12}
\begin{array}{lll}
\overline U&\rightarrow&\overline U(1)\\
\uparrow&&\uparrow\\
\overline T&\rightarrow&\overline T(1)
\end{array}
\end{equation}
such that
$$
\overline x_1^{\overline v_{j,1}}\cdots\overline x_{\overline r}^{\overline v_{j,\overline r}}
= \overline x_1(1)^{\overline v(1)_{j,1}}\cdots\overline x_{\overline r}(1)^{\overline v(1)_{j,\overline r}}
$$
with $\overline v_j(1)
=(\overline v(1)_{j,1},\ldots,\overline v(1)_{j,\overline r})\in{\bf N}^{\overline r}$ for $1\le j\le\overline a$.

We then have expressions for all $\Lambda=(\lambda_1,\ldots,\lambda_{\overline s})\in{\bf N}^{\overline s}$, where $\overline u_1,\ldots,\overline u_{\overline b}\in{\bf Q}^{\overline r}$ depend on $\Lambda$, 
\begin{equation}\label{eq16}
h_{[\Lambda]}=\overline y_1(1)^{\lambda_1(1)}\cdots
\overline y_{\overline s}(1)^{\lambda_{\overline s}(1)}
\left[\sum_{i=1}^{\overline b}\phi_1^{-\overline u_{i,1}}\cdots\phi_{\overline r}^{-\overline u_{i,\overline r}}\overline x_1(1)^{\overline u_{i,1}(1)}\cdots
\overline x_{\overline r}(1)^{\overline u_{i,\overline r}(1)}g_i\right]
\end{equation}
with $g_i\in k(\overline U)[[\overline x_1(1),\ldots,\overline x_{\overline r+l}(1)]]$,
$$
\overline y_1(1)^{\lambda_1(1)}\cdots
\overline y_{\overline s}(1)^{\lambda_{\overline s}(1)}
=\overline y_1^{\lambda_1}\cdots
\overline y_{\overline s}^{\lambda_{\overline s}}
$$
and
$$
\overline x_1(1)^{\overline u_{i,1}(1)}\cdots
\overline x_{\overline r}(1)^{\overline u_{i,\overline r}(1)}
=\overline x_1^{\overline u_{i,1}}\cdots
\overline x_{\overline r}^{\overline u_{i,\overline r}}
$$
for $1\le i\le \overline b$.

\begin{Remark}\label{Remark3}
If $\Lambda=0$, we have $\nu(\overline x_1^{\overline u_{i,1}}\cdots
\overline x_{\overline r}^{\overline u_{i,\overline r}})\ge 0$ for $1\le i\le\overline b$,
so we can construct our CUTS (\ref{eq12}) so that
$\overline u_i(1)\in{\bf Q}_+^{\overline r}$ for $1\le i\le \overline b$.
Thus if $d$ is a common denominator of the coefficients of the $\overline u_i$,
then $\overline u_i(1)\in\frac{1}{d}{\bf N}^{\overline r}$ for $1 \le i\le\overline b$,
and we have 
\begin{equation}\label{eq5}
h_{[\Lambda]}=h_0\in k(\overline U)[\phi_1^{\frac{1}{d}},\ldots\phi_{\overline r}^{\frac{1}{d}}][[
\overline x_1(1)^{\frac{1}{d}},\ldots,\overline x_{\overline r}(1)^{\frac{1}{d}},
\overline x_{\overline r+1}(1),\ldots,\overline x_{\overline r+l}(1)]].
\end{equation}
\end{Remark}

\begin{Lemma}\label{Lemma12} 
Suppose that $h\in k(\overline U)[[\overline x_1,\ldots,\overline x_m]]\subset \overline U$
is such that $\nu^*(h)<\infty$. Then $\nu^*(h)\in \Gamma_{\nu}\otimes{\bf Q}$.
\end{Lemma}

\begin{pf} Recall that $t_1,\ldots, t_{\overline\alpha}$ is a transcendence basis of $k(\overline U)$ over $k(\overline T)$. Let $A=k(\overline U)[[\overline x_1,\ldots,\overline x_m]]$,
$B=k(\overline T)(t_1,\ldots,t_{\overline\alpha})[[\overline x_1,\ldots,\overline x_m]]$.
We will first assume that $h\in B$. Since $\nu^*(h)<\infty$,
 there exists $s>0$ such that $\nu^*(m(B)^s)>\nu^*(h)$. We have
$$
h=f+g
$$
with $f\in k(\overline T)(t_1,\ldots,t_{\overline\alpha})[\overline x_1,\ldots,\overline x_{m}]$ and $g\in m(B)^s$. $\nu^*(h)=\nu^*(f)$.
There exists $0\ne a\in k(\overline T)[t_1,\ldots, t_{\overline\alpha}]$ such that
we have a finite expansion
$$
af=\sum_{I=(i_1,\ldots,i_{\overline\alpha})\in{\bf N}^{\overline\alpha}}
a_It_1^{i_1}\cdots t_{\overline\alpha}^{i_{\overline\alpha}}
$$
with $a_I\in k(\overline T)[\overline x_1,\ldots,\overline x_{\overline r+l}]$.
By Lemma \ref{Lemma4},
$$
\nu^*(h)=\nu^*(f)=\text{ min }\{\nu(a_I)\mid a_I\ne 0\}\in\Gamma_{\nu}.
$$

Since $A$ is a finite extension of $B$, we have
$\nu^*(h)\in \Gamma_{\nu}\otimes{\bf Q}$ if $h\in B$,
by the Corollary to Lemma 3, Section 11, Chapter VI \cite{ZS}.

\end{pf}

\begin{Lemma}\label{Lemma2} 
\begin{enumerate}
\item Suppose that $\Lambda\in{\bf N}^{\overline s}$ and $\nu^*(h_{[\Lambda]})<\infty$.
 Then  we have 
$$
\nu^*\left(\frac{h_{[\Lambda]}}{\overline y_1^{\lambda_1}\cdots\overline y_{\overline s}^{\lambda_{\overline s}}}\right)\in\Gamma_{\nu}\otimes{\bf Q}.
$$
In particular, 
$\nu^*(h_{[\Lambda]})\in \Gamma_{\nu}\otimes{\bf Q}$ implies $[\Lambda]=0$.
\item
In the expansion (\ref{eq2}), for $\Lambda_1,
\Lambda_2\in{\bf N}^{\overline s}$, suppose that
$$
\nu^*(h_{[\Lambda_1]})=\nu^*(h_{[\Lambda_2]})<\infty. 
$$
Then $[\Lambda_1]=[\Lambda_2]$.
\end{enumerate}
\end{Lemma}

\begin{pf} For $\Lambda\in{\bf N}^{\overline s}$ such that $\nu^*(h_{\Lambda]})<\infty$, consider the expansion (\ref{eq16}) of $h_{[\Lambda]}$.

There exists $w=(w_1,\ldots,w_{\overline r})\in{\bf N}^{\overline r}$
such that $w+\overline u_i\in{\bf Q}_+^{\overline r}$ for $1\le i\le \overline b$.

Let $d$ be a common denominator of the coefficients $\overline u_i$ for $1\le i\le\overline b$.
Let 
$$
A=k(\overline U)[\phi_1^{\frac{1}{d}},\ldots,\phi_{\overline r}^{\frac{1}{d}}]
[[\overline x_1(1)^{\frac{1}{d}},\ldots,\overline x_{\overline r}(1)^{\frac{1}{d}},
\overline x_{\overline r+1}(1),\ldots,\overline x_{\overline r+l}(1)]].
$$
Set 
$$
f=\frac{h_{[\Lambda]}}{\overline y_1^{\lambda_1}\cdots\overline y_{\overline s}^{\lambda_{\overline s}}}\overline x_1^{w_1}\cdots\overline x_{\overline r}^{w_{\overline r}}
\in  A.
$$
If we restrict $\tilde\nu^*$ to $Q(k(\overline U)[[\overline y_1(1),\ldots,\overline y_{\overline s+l}(1)]])$, extend it to the finite extension
$Q(k(\overline U)[\phi_1^{\frac{1}{d}},\ldots,\phi_{\overline r}^{\frac{1}{d}}]
[\overline y_1(1)^{\frac{1}{d}},\ldots,\overline y_{\overline s}(1)^{\frac{1}{d}},
\overline y_{\overline s+1}(1),\ldots,\overline y_{\overline s+l}(1)]])
$ so that it dominates $k(\overline U)[\phi_1^{\frac{1}{d}},\ldots,\phi_{\overline r}^{\frac{1}{d}}]
[\overline y_1(1)^{\frac{1}{d}},\ldots,\overline y_{\overline s}(1)^{\frac{1}{d}},
\overline y_{\overline s+1}(1),\ldots,\overline y_{\overline s+l}(1)]]$
and restrict to $A$, we get a valuation $\overline\nu'$ on $Q(A)$ which
extends $\tilde\nu$ restricted to $Q(k(\overline T)[[\overline x_1(1),\ldots,\overline x_{\overline r+l}(1)]])$.  By Lemma \ref{Lemma12} and the Corollary to Lemma 3, Section 11,
Chapter VI \cite{ZS}, we have 
$
\overline\nu'(f)\in \Gamma_{\nu}\otimes{\bf Q}.
$
Thus, we have
$$
\nu^*(\frac{h_{[\Lambda]}}{\overline y_1^{\lambda_1}\cdots\overline y_{\overline s}^{\lambda_{\overline s}}})\in 
\Gamma_{\nu}\otimes{\bf Q}.
$$
We now compare $\nu^*(h_{[\Lambda_1]})$ and $\nu^*(h_{[\Lambda_2]})$.
$$
\nu^*\left(\frac{h_{[\Lambda_1]}}{h_{[\Lambda_2]}}\right)
=\nu^*\left(\frac{h_{[\Lambda_1]}}{\overline y_1^{\lambda_1^1}\cdots\overline y_{\overline s}^{\lambda_{\overline s}^1}}\right)
-\nu^*\left(\frac{h_{[\Lambda_2]}}{\overline y_1^{\lambda_1^2}\cdots\overline y_{\overline s}^{\lambda_{\overline s}^2}}\right)
+\nu^*\left(\frac{\overline y_1^{\lambda_1^1}\cdots\overline y_{\overline s}^{\lambda_{\overline s}^1}}{\overline y_1^{\lambda_1^2}\cdots\overline y_{\overline s}^{\lambda_{\overline s}^2}}\right)
$$
where $\Lambda_1=(\lambda_1^1,\ldots,\lambda_{\overline s}^1)$, $\Lambda_2=(\lambda_1^2,\ldots,\lambda_{\overline s}^2)$.
Thus $\nu^*(h_{[\Lambda_1]})=\nu^*(h_{[\Lambda_2]})<\infty$
implies 
$$
\nu^*\left(\frac{\overline y_1^{\lambda_1^1}\cdots\overline y_{\overline s}^{\lambda_{\overline s}^1}}{\overline y_1^{\lambda_1^2}\cdots\overline y_{\overline s}^{\lambda_{\overline s}^2}}\right)
\in \Gamma_{\nu}\otimes{\bf Q},
$$
so that
$$
(\lambda_1^1-\lambda_1^2)\nu^*(\overline y_1)+\cdots+(\lambda_{\overline s}^1-\lambda_{\overline s}^2)\nu^*(\overline y_{\overline s})\in\Gamma_{\nu}\otimes{\bf Q}.
$$
Thus, there exists $(a_1,\ldots,a_{\overline r})\in{\bf Q}^{\overline r}$ such that
$$
(\lambda_1^1-\lambda_1^2)\nu^*(\overline y_1)+\cdots+(\lambda_{\overline s}^1-\lambda_{\overline s}^2)\nu^*(\overline y_{\overline s})
=a_1\nu(\overline x_1)+\cdots+ a_{\overline r}\nu(\overline x_{\overline r}).
$$
Substituting from (\ref{eq1}), we get
$$
(a_1,\ldots,a_{\overline r})C=(\lambda_1^1-\lambda_1^2,\ldots,\lambda_{\overline s}^1
-\lambda_{\overline s}^2)
$$
and thus
$$
\Lambda_1-\Lambda_2\in \Phi({\bf Q}^{\overline r})\cap{\bf Z}^{\overline s}.
$$
\end{pf}

\begin{Remark}\label{Remark2} 
In the expansion (\ref{eq2}), Let $\Lambda_0\in{\bf N}^{\overline s}$ be such that
$$
\nu^*(h_{[\Lambda_0]})=\text{min}\{\nu^*(h_{[\Lambda]})\mid \Lambda\in{\bf N}^{\overline s}\}.
$$
This minimum exists since $\overline U$ is Noetherian.
Then, by Lemma \ref{Lemma2}, 
$$
\nu^*(f)=\nu^*(h_{[\Lambda_0]}).
$$
\end{Remark}

\begin{Lemma}\label{Lemma3}
With the notation of (\ref{eq1}), assume that $f\in\overline T''$ and 
$$
f\in k(\overline U)[[\overline y_1,\ldots,
\overline y_{\overline s+l}]].
$$
Then
$$
f\in k(\overline T)[[\overline x_1,\ldots,
\overline x_{\overline r+l}]].
$$
\end{Lemma}

\begin{pf} Let $h_{[\Lambda_0]}$ be the minimum value term of $f$ in the expansion (\ref{eq2}),
so that 
$$
\nu^*(f)=\nu^*(h_{[\Lambda_0]})
$$
 (by Remark \ref{Remark2}). Since $f\in \overline T''$, we have $\nu^*(f)<\infty$.
Since
$$
\nu^*(h_{[\Lambda_0]})=\nu(f)\in\Gamma_{\nu},
$$
 we have $[\Lambda_0]=0$ by Lemma \ref{Lemma2}.
Thus by Remark \ref{Remark3}, there exists a CUTS of type I
$$
\begin{array}{lll}
\overline U&\rightarrow&\overline U(1)\\
\uparrow&&\uparrow\\
\overline T&\rightarrow&\overline T(1)
\end{array}
$$
and $d$, a positive integer, such that 
$$
h_{[\Lambda_0]}=
h_{0}\in A=k(\overline U)[\phi_1^{\frac{1}{d}},\ldots,\phi_{\overline r}^{\frac{1}{d}}]
[[\overline x_1(1)^{\frac{1}{d}},\ldots,\overline x_{\overline r}(1)^{\frac{1}{d}},\overline
x_{\overline r+1}(1),\ldots,\overline x_{\overline r+l}(1)]].
$$

Suppose that $f\ne h_0$. Then there exists $\Lambda_1\not\in({\bf Q}^{\overline r}C)\cap{\bf Z}^{\overline s}$ in the expansion of $f$ in (\ref{eq2}), such that the minimal value
term of $f-h_0$ is $h_{[\Lambda_1]}$. Write $\Lambda_1=(\lambda_1^1,\ldots,\lambda_{\overline s}^1)$.

Consider the $k$-derivation
$$
\partial=\sum_{i=1}^{\overline s}e_i\overline y_i\frac {\partial y_i}{\partial\overline y_i}
$$
on $\overline U''$ (and on $\overline U$) where $e_i\in{\bf Q}$ are chosen so that 
$$
\sum_{i=1}^{\overline s}e_ic_{ji}=0\text{ for }1\le j\le \overline r
$$
and
$$
\sum_{i=1}^{\overline s}e_i\lambda_i^1\ne 0.
$$
We have
$$
\partial(\overline y_1^{b_1}\cdots\overline y_{n}^{b_n})
=(b_1e_1+\cdots+b_{\overline s}e_{\overline s})\overline y_1^{b_1}\cdots\overline y_{n}^{b_n}
$$
for all monomials $\overline y_1^{b_1}\cdots\overline y_{n}^{b_n}$. From (\ref{eq15}),
we see that
$$
\begin{array}{ll}
\partial(h_{[\Lambda]})&=\sum_{\alpha\in{\bf N}^{\overline s}\mid [\alpha]=[\Lambda]} g_{\alpha}\partial(\overline y_1^{\alpha_1}\cdots\overline y_{\overline s}^{\alpha_{\overline s}})\\
&=\sum(\alpha_1e_1+\cdots+\alpha_{\overline s}e_{\overline s})g_{\alpha}\overline y_1^{\alpha_1}
\cdots\overline y_{\overline s}^{\alpha_{\overline s}}\\
&=(\lambda_1 e_1+\cdots+\lambda_{\overline s}e_{\overline s})h_{[\Lambda]}.
\end{array}
$$
In particular, $\partial(h_{[\Lambda_0]})=0$ and 
$$
\partial(h_{[\Lambda_1]})=(\lambda_1^1\alpha_1+\cdots+\lambda_{\overline s}^1\alpha_{\overline s})h_{[\Lambda_1]}\ne 0.
$$
Thus
$$
\partial(f-h_0)=\sum_{[\Lambda]\ne 0}(\lambda_1e_1+\cdots+\lambda_{\overline s}e_{\overline s})h_{[\Lambda]}
$$
has value
$$
\nu^*(\partial(f-h_0))=\nu^*(h_{[\Lambda_1]}).
$$
But $\partial$ is a derivation of $\overline U''$, so that
$\partial(f-h_0)=\partial(f)\in\overline U''$ has finite value, and
$$
\nu^*(h_{[\Lambda_1]})=\nu^*(\partial(f-h_0))<\infty.
$$
Thus $\nu^*(h_{[\Lambda_1]})\not\in \Gamma_{\nu}\otimes{\bf Q}$, by Lemma \ref{Lemma2},
but
$\nu^*(h_{[\Lambda_1]})=\nu^*(f-h_0)\in \Gamma_{\nu}\otimes{\bf Q}$, by Lemma \ref{Lemma12},
since
$$
f-h_0\in k(\overline U)[\phi_1^{\frac{1}{d}},\ldots,\phi_{\overline r}^{\frac{1}{d}}]
[[\overline x_1(1)^{\frac{1}{d}},\ldots,\overline x_{\overline r}(1)^{\frac{1}{d}},
\overline x_{\overline r+1}(1),\ldots,\overline x_m(1)]]
$$
which is a finite extension of $k(\overline U)[[\overline x_1(1),\ldots,\overline x_m(1)]]$,
a contradiction. Thus $f=h_0$.
 We have 
$$
f\in A\cap k(\overline T)[[\overline x_1,\ldots,\overline x_m]]=k(\overline T)[[\overline x_1,\ldots,\overline x_{\overline r+l}]].
$$
\end{pf}

\begin{Theorem}\label{Theorem4.7}
Suppose that $(R,\overline T'',\overline T)$
and $(S,\overline U'',\overline U)$ is a CUTS along $\nu^*$, such that $\overline T''$
contains a subfield isomorphic to $k(c_0)$ for some $c_0\in \overline T''$
and $\overline U''$ contains a subfield isomorphic to $k(\overline U)$.
$\overline T''$ has regular 
parameters $(\overline z_1,\ldots,\overline z_m)$ and $\overline U''$ has regular
parameters $(\overline w_1,\ldots,\overline w_n)$ with
$$
\begin{array}{ll}
\overline z_1&=\overline w_1^{c_{11}}\cdots \overline w_{\overline s}^{c_{1\overline s}}\phi_1\\
&\vdots\\
\overline z_{\overline r}&=\overline w_1^{c_{\overline r1}}\cdots  \overline w_{\overline s}^{c_{\overline r\overline s}}\phi_{\overline r}\\
\overline z_{\overline r+1}&=\overline w_{\overline s+1}\\
&\vdots\\
\overline z_{\overline r+l}&=\overline w_{\overline s+l}
\end{array}
$$

such that $\phi_1,\ldots,\phi_{\overline r}\in k(\overline U)$, $\nu(\overline z_1),\ldots,\nu(\overline z_{\overline r})$ are rationally independent, 
$\nu^*(\overline w_1),\ldots,\nu^*(\overline w_{\overline s})$ are rationally independent
and $(c_{ij})$ has rank $\overline r$. 

Suppose that one of the following three conditions hold.

\begin{equation}\label{eq42}
f\in k(\overline U)[[\overline w_1,\ldots,\overline w_{\overline s+\overline m}]]
\text{ for some $\overline m$ such that $0 \le\overline m\le n-\overline s$ with }
\nu^*(f)<\infty.
\end{equation}

\begin{equation}\label{eq43}
f\in k(\overline U)[[\overline w_1,\ldots,\overline w_{\overline s+\overline m}]]
\text{ for some $\overline m$ such that $0<\overline m\le n-\overline s$ with }
\nu^*(f)=\infty\text{ and $A>0$ is given}.
\end{equation}

\begin{equation}\label{eq44}
f\in \left(k(\overline U)[[\overline w_1,\ldots,\overline w_{\overline s+\overline m}]]
-k(\overline U)[[\overline w_1,\ldots,\overline w_{\overline s+l}]]\right)\cap\overline U''
\text{ for some $\overline m$ such that $l<\overline m\le n-\overline s$}.
\end{equation}
Then there exists a CUTS along $\nu^*$ $(R,\overline T''(t),\overline T(t))$ and
$(S,\overline U''(t),\overline U(t))$ 
\begin{equation}\label{eq45}
\begin{array}{lllllllll}
\overline U&=&\overline U(0)&\rightarrow &\overline U(1)&\rightarrow&\cdots&\rightarrow
&\overline U(t)\\
&&\uparrow&&\uparrow&&&&\uparrow\\
\overline T&=&\overline T(0)&\rightarrow &\overline T(1)&\rightarrow&\cdots&\rightarrow
&\overline T(t)
\end{array}
\end{equation}
such that 
$\overline U''(t)$ has regular
parameters $(\overline w_1(t),\ldots,\overline w_n(t))$.

In case (\ref{eq42}) we have
$$
f=\overline w_1(t)^{d_1}\cdots\overline w_{\overline s}(t)^{d_{\overline s}}u
$$
where $u\in k(\overline U(t))[[\overline w_1(t),\ldots,\overline w_{\overline s+\overline m}(t)]]$
is a unit power series.

In case (\ref{eq43}) we have
 $$
f=\overline w_1(t)^{d_1}\cdots\overline w_{\overline s}(t)^{d_{\overline s}}\Sigma
$$
where $\Sigma\in k(\overline U(t))[[\overline w_1(t),\ldots,\overline w_{\overline s+\overline m}(t)]]$,
$\nu^*(\overline w_1(t)^{d_1}\cdots\overline w_{\overline s}(t)^{d_{\overline s}})>A$.

In case (\ref{eq44}) we have 
$$
f=P+\overline w_1(t)^{d_1}\cdots\overline w_{\overline s}(t)^{ d_{\overline s}}H
$$
for some powerseries $P\in k(\overline U(t))[[\overline w_1(t),\ldots,\overline w_{\overline s+l}(t)]]$,
$$
H=u(\overline w_{\overline s+\overline m}(t)+\overline w_1(t)^{g_1}\cdots\overline w_{\overline s}(t)^{g_{\overline s}}
\Sigma)
$$
where $u\in k(\overline U(t))[[\overline w_1(t),\ldots,\overline w_{\overline s+\overline m}(t)]]$
is a unit, $\Sigma\in k(\overline U(t))[[\overline w_1(t),\ldots,\overline w_{\overline s+\overline m-1}(t)]]$ and $\nu^*(\overline w_{\overline s+\overline m}(t))\le\nu^*(\overline w_1(t)^{g_1}\cdots\overline w_{\overline s}(t)^{g_{\overline s}})$.

(\ref{eq45}) will be such that $\overline T''(\alpha)$ has regular parameters
$$
(\overline z_1(\alpha),\ldots,\overline z_m(\alpha))\text{ and }(\tilde{\overline z}'_1(\alpha),
\ldots,\tilde{\overline z}'_m(\alpha)),
$$
$\overline U''(\alpha)$ has regular parameters
$$
(\overline w_1(\alpha),\ldots,\overline w_n(\alpha))\text{ and }(\tilde{\overline w}'_1(\alpha),
\ldots,\tilde{\overline w}'_n(\alpha))
$$
where $\overline z_i(0)=\overline z_i$ for $1\le i\le m$ and $\overline w_i(0)=\overline w_i$
for $1\le i\le n$. (\ref{eq45}) will consist of three types of CUTS.
\begin{description}
\item[(M1)] $\overline T(\alpha)\rightarrow \overline T(\alpha+1)$ and $\overline U(\alpha)
\rightarrow \overline U(\alpha+1)$ are of type I.
\item[(M2)] $\overline T(\alpha)\rightarrow\overline T(\alpha+1)$ is of type $II_r$,
$1\le r\le \text{ min }\{l,\overline m\}$, and $\overline U(\alpha)\rightarrow \overline U(\alpha+1)$ is a transformation 
of type $II_r$ followed by a transformation of type $I$.
\item[(M3)] $\overline T(\alpha)=\overline T(\alpha+1)$ and $\overline U(\alpha)\rightarrow
\overline U(\alpha+1)$ is of type $II_r$ with $l+1\le r\le \overline m$.
\end{description}
$\overline T''(\alpha)$ contains a subfield isomorphic to $k(c_0,\ldots,c_{\alpha})$ and
$\overline U''(\alpha)$ contains a subfield isomorphic to $k(\overline U(\alpha))$.

We will find polynomials $P_{i,\alpha}$ so that the variables will be related by:
$$
\tilde{\overline z}'_i(\alpha)=\left\{
\begin{array}{ll}
\overline z_i(\alpha)-P_{i,\alpha}&\text{if }\overline r+1\le i\le \overline r+
\text{ min }\{l,\overline m\}\\
\overline z_i(\alpha)&\text{otherwise}
\end{array}\right.
$$
$$
\tilde{\overline w}_i'(\alpha)=\left\{
\begin{array}{ll}
\tilde{\overline z}_{i-\overline s +\overline r}(\alpha)&\text{if }\overline s+1\le i\le 
\overline s+l\\
\overline w_i(\alpha)-P_{i,\alpha}&\text{if }\overline s+l+1\le i\le \overline s+\overline m\\
\overline w_i(\alpha)&\text{otherwise}
\end{array}\right.
$$
We will have $P_{i,\alpha}\in k(c_0,\ldots,c_{\alpha})[\overline z_1(\alpha),\ldots,\overline z_{i-1}(\alpha)]$ if $i\le\overline r+l$,
$$
P_{i,\alpha}\in k(\overline U(\alpha))[\overline w_1(\alpha),\ldots,\overline w_{i-1}(\alpha)]
$$ 
if $i>\overline s+l$.
 For all $\alpha$ we will have

\begin{equation}\label{eq46}
\begin{array}{ll}
\overline z_1(\alpha)&=\overline w_1(\alpha)^{c_{11}(\alpha)}\cdots\overline w_{\overline s}(\alpha)^{c_{1\overline s}(\alpha)}\phi_1(\alpha)\\
&\vdots\\
\overline z_{\overline r}(\alpha)&=\overline w_1(\alpha)^{c_{\overline r 1}(\alpha)}\cdots\overline w_{\overline s}^{c_{\overline r\overline s}}\phi_{\overline r}(\alpha)\\
\overline z_{\overline r+1}(\alpha)&=\overline w_{\overline s+1}(\alpha)\\
&\vdots\\
\overline z_{\overline r+l}(\alpha)&=\overline w_{\overline s+l}(\alpha)
\end{array}
\end{equation}
and

\begin{equation}\label{eq47}
\begin{array}{ll}
\tilde{\overline z}'_1(\alpha)&=\tilde{\overline w}'_1(\alpha)^{c_{11}(\alpha)}\cdots
\tilde{\overline w}'_{\overline s}(\alpha)^{c_{1\overline s}(\alpha)}\phi_1(\alpha)\\
&\vdots\\
\tilde{\overline z}'_{\overline r}(\alpha)&=
\tilde{\overline w}'_1(\alpha)^{c_{\overline r 1}(\alpha)}
\cdots
\tilde{\overline w}'_{\overline s}(\alpha)^{c_{\overline r\overline s}(\alpha)}\phi_{\overline r}(\alpha)\\
\tilde{\overline z}'_{\overline r+1}(\alpha)&=\tilde{\overline w}'_{\overline s+1}(\alpha)\\
&\vdots\\
\tilde{\overline z}'_{\overline r+l}(\alpha)&=\tilde{\overline w}'_{\overline s+l}(\alpha)
\end{array}
\end{equation}
where $\phi_1(\alpha),\ldots,\phi_{\overline r}(\alpha)\in k(\overline U(\alpha))$.
 $\nu(\overline z_1(\alpha)),\ldots,\nu(\overline z_{\overline r}(\alpha))$ are rationally independent, 
$\nu^*(\overline w_1(\alpha)),\ldots,\nu^*(\overline w_{\overline s}(\alpha))$ are rationally independent
and $(c_{ij}(\alpha))$ has rank $\overline r$ for $1\le \alpha\le t$.

In a transformation $\overline T(\alpha)\rightarrow\overline T(\alpha+1)$ of type $I$,
$\overline T''(\alpha+1)$ will have regular parameters $(\overline z_1(\alpha+1),\ldots,
\overline z_m(\alpha+1))$ defined by 
\begin{equation}\label{eq48}
\begin{array}{ll}
\tilde{\overline z}'_1(\alpha)&=\overline z_1(\alpha+1)^{a_{11}(\alpha+1)}\cdots
\overline z_{\overline r}(\alpha+1)^{a_{1\overline r}(\alpha+1)}\\
&\vdots\\
\tilde{\overline z}'_{\overline r}(\alpha)&=\overline z_1(\alpha+1)^{a_{\overline r1}(\alpha+1)}
\cdots\overline z_{\overline r}(\alpha+1)^{a_{\overline r\overline r}(\alpha+1)}.
\end{array}
\end{equation}
and $c_{\alpha+1}$ is defined to be 1. In a transformation $\overline T(\alpha)\rightarrow
\overline T(\alpha+1)$ of type $II_r$ ($1\le r\le \text{ min }\{l,\overline m\}$) $\overline T''(\alpha+1)$
will have regular parameters 
$(\overline z_1(\alpha+1),\ldots,
\overline z_m(\alpha+1))$ defined by 
\begin{equation}\label{eq49}
\begin{array}{ll}
\tilde{\overline z}'_1(\alpha)&=\overline z_1(\alpha+1)^{a_{11}(\alpha+1)}\cdots
\overline z_{\overline r}(\alpha+1)^{a_{1\overline r}(\alpha+1)}
c_{\alpha+1}^{a_{1\overline r+1}(\alpha+1)}\\
&\vdots\\
\tilde{\overline z}'_{\overline r}(\alpha)&=\overline z_1(\alpha+1)^{a_{\overline r1}(\alpha+1)}
\cdots\overline z_{\overline r}(\alpha+1)^{a_{\overline r\overline r}(\alpha+1)}
c_{\alpha+1}^{a_{\overline r,\overline r+1}(\alpha+1)}\\
\tilde{\overline z}'_{\overline r+r}(\alpha)&=
\overline z_1(\alpha+1)^{a_{\overline r+1,1}(\alpha+1)}\cdots\overline z_{\overline r}(\alpha+1)
^{a_{\overline r+1,\overline r}(\alpha+1)}\\
&\,\,\,\,\,\,\,\,\,\,\cdot
(\overline z_{\overline r+r}(\alpha+1)+1)c_{\alpha+1}^{a_{\overline r+1,
\overline r+1}(\alpha+1)}
\end{array}
\end{equation}
In a transformation $\overline U(\alpha)\rightarrow\overline U(\alpha+1)$ of type $I$
$\overline U''(\alpha+1)$ will have regular parameters $(\overline w_1(\alpha+1),\ldots,
\overline w_n(\alpha+1))$ defined by 
\begin{equation}\label{eq50}
\begin{array}{ll}
\tilde{\overline w}'_1(\alpha)&=\overline w_1(\alpha+1)^{b_{11}(\alpha+1)}\cdots
\overline w_{\overline s}(\alpha+1)^{b_{1\overline s}(\alpha+1)}\\
&\vdots\\
\tilde{\overline w}'_{\overline s}(\alpha)&=\overline w_1(\alpha+1)^{b_{\overline s1}(\alpha+1)}
\cdots\overline w_{\overline s}(\alpha+1)^{b_{\overline s\overline s}(\alpha+1)}.
\end{array}
\end{equation}
 In a transformation $\overline U(\alpha)\rightarrow
\overline U(\alpha+1)$ of type $II_r$ ($1\le r\le \overline m$) $\overline U''(\alpha+1)$
will have regular parameters 
$(\overline w_1(\alpha+1),\ldots,
\overline w_n(\alpha+1))$ defined by 
\begin{equation}\label{eq51}
\begin{array}{ll}
\tilde{\overline w}'_1(\alpha)&=\overline w_1(\alpha+1)^{b_{11}(\alpha+1)}\cdots
\overline w_{\overline s}(\alpha+1)^{b_{1\overline s}(\alpha+1)}
d_{\alpha+1}^{b_{1\overline s+1}(\alpha+1)}\\
&\vdots\\
\tilde{\overline w}'_{\overline s}(\alpha)&=\overline w_1(\alpha+1)^{b_{\overline s1}(\alpha+1)}
\cdots\overline w_{\overline s}(\alpha+1)^{b_{\overline s\overline s}(\alpha+1)}
d_{\alpha+1}^{b_{\overline s,\overline s+1}(\alpha+1)}\\
\tilde{\overline w}'_{\overline s+r}(\alpha)&=
\overline w_1(\alpha+1)^{b_{\overline s+1,1}(\alpha+1)}\cdots\overline w_{\overline s}(\alpha+1)
^{b_{\overline s+1,\overline s}(\alpha+1)}\\
&\,\,\,\,\,\,\,\,\,\,\cdot
(\overline w_{\overline s+r}(\alpha+1)+1)d_{\alpha+1}^{b_{\overline s+1,
\overline s+1}(\alpha+1)}
\end{array}
\end{equation}
We will call a CUTS as in (\ref{eq45}) a CUTS of form $\overline m$.
Observe that the UTS $\overline T\rightarrow \overline T(t)$ is a UTS of form
$\text{ min }\{l,\overline m\}$.
\end{Theorem}

\begin{pf}  
We will first assume that $f$ satisfies (\ref{eq42}) or (\ref{eq43}) with
$0\le \overline m\le l$.

By (A1) of Theorem \ref{Theorem3}, after performing a CUTS of form $\overline m$, we may assume that 
\begin{equation}\label{eq20}
p_{\overline m}=(\overline z_{r(1)}-Q_{r(1)}(\overline z_1,\ldots,\overline z_{r(1)-1}),\ldots,
\overline z_{r(\tilde m)}-Q_{r(\tilde m)}(\overline z_1,\ldots,\overline z_{r(\tilde m)-1})
\end{equation} 
where the coefficients of $Q_{r(i)}$ are in $k(c_0)$. We have that $f\in k(\overline U)[[\overline w_1,\ldots,\overline w_{\overline s+\overline m}]]$.

Given a CUTS (\ref{eq45}), define $\sigma(i)$ as in Lemma \ref{Lemma17} for the UTS
$$
\overline T\rightarrow\cdots\rightarrow\overline T(t).
$$

If $\sigma(i)$ drops during the course of the proof, we can start the corresponding
algorithm again with this smaller value of $\sigma(i)$. Eventually $\sigma(i)$ must stablize,
so we may assume that $\sigma(i)$ is constant throughout the proof. 

We have the expansion
$$
f=\sum_{[\Lambda]\in {\bf Z}^{\overline s}/({\bf Q}^{\overline r}C)\cap{\bf Z}^{\overline s}}h_{[\Lambda]}
$$
of (\ref{eq2}). Let $\Lambda_0\in{\bf N}^{\overline s}$ be such that 
$$
\nu^*(h_{[\Lambda_0]})=\text{min}\{\nu^*(h_{[\Lambda]})\mid h_{[\Lambda]}\ne 0\}.
$$
By Remark \ref{Remark2}, $\nu^*(f)=\nu^*(h_{[\Lambda_0]})$. We either have $\nu^*(f)<\infty$ and
$\nu^*(h_{[\Lambda]})>\nu^*(h_{[\Lambda_0]})$ if $[\Lambda]\ne [\Lambda_0]$ (by Lemma \ref{Lemma2}) or $\nu^*(f)=\infty$ and $\nu^*(h_{[\Lambda]})=\infty$ for all $[\Lambda]$.
Let $I$ be the ideal 
$$
I=(h_{[\Lambda]}\mid [\Lambda]\ne[\Lambda_0])\subset k(\overline U)[[\overline w_1,\ldots,
\overline w_{\overline s+\overline m}]].
$$
Let $h_{[\Lambda_1]},\ldots,h_{[\Lambda_{\beta}]}$ be generators of $I$. We will construct
a CUTS along $\nu^*$ of form $\overline m$
$$
\begin{array}{lll}
\overline U&\rightarrow&\overline U(\alpha)\\
\uparrow&&\uparrow\\
\overline T&\rightarrow&\overline T(\alpha)
\end{array}
$$
such that
$$
h_{[\Lambda_i]}=\overline w_1(\alpha)^{b_1^i}\cdots\overline w_{\overline s}(\alpha)^{b_{\overline s}^i}\psi_i
$$
for $0\le i\le \beta$, with $\psi_i\in \overline U(\alpha)$.
If $\nu^*(h_{[\Lambda_i]})<\infty$, $\psi_i$ will be a unit.

 If $\nu^*(f)<\infty$, so that
$f$ satisfies the conditions of (\ref{eq42}), then set 
$$
B=\nu^*(f)=\nu^*(h_{[\Lambda_0]}).
$$
 If $\nu^*(h_{[\Lambda_i]})
=\infty$, we will have $\nu^*(\overline w_1(\alpha)^{b_1^i}\cdots\overline w_{\overline s}(\alpha)^{b_{\overline s}^i})>B$.

 If $\nu^*(f)=\infty$, so that $f$ satisfies the conditions of 
(\ref{eq43}), we will have 
$$
\nu^*(\overline w_1(\alpha)^{b_1^i}\cdots\overline w_{\overline s}(\alpha)^{b_{\overline s}^i})>A
$$
 for all $i$. In this case $\nu^*(h_{[\Lambda_i]})=\infty$
for $0\le i\le \beta$.

Assume that the above CUTS has been constructed.
There exists (by Lemma 4.2 \cite{C1}) a CUTS of type (M1) along $\nu^*$
$$
\begin{array}{lll}
\overline U(\alpha)&\rightarrow&\overline U(\alpha+1)\\
\uparrow&&\uparrow\\
\overline T(\alpha)&=&\overline T(\alpha+1)
\end{array}
$$
such that if $\nu^*(f)<\infty$, then for $1\le i\le \beta$, $h_{[\Lambda_0]}$ divides 
$h_{[\Lambda_i]}$ in $\overline U(\alpha+1)$.  If $\nu^*(f)=\infty$, then there exists
$\overline w_1(\alpha+1)^{a_1}\cdots\overline w_{\overline s}(\alpha+1)^{a_{\overline s}}$
such that 
$$
\nu^*(\overline w_1(\alpha+1)^{a_1}\cdots\overline w_{\overline s}(\alpha+1)^{a_{\overline s}})>A
$$
 and
$\overline w_1(\alpha+1)^{a_1}\cdots\overline w_{\overline s}(\alpha+1)^{a_{\overline s}}$
divides $h_{[\Lambda_i]}$ for $0\le i\le \beta$ in $\overline U(\alpha+1)$.
Thus the conclusions of the theorem hold for $f$ satisfying the conditions (\ref{eq42}) or (\ref{eq43}) with $0\le \overline m\le l$.

We are thus reduced to proving the theorem (with our assumption that $f$ satisfies
(\ref{eq42}) or (\ref{eq43}) and $0\le \overline m\le l$) when $f=h_{[\Lambda]}$ for some
$\Lambda\in{\bf N}^{\overline s}$. Assume that $f$ has this form. There exists a CUTS
along $\nu^*$ (using Lemma \ref{Lemma4.3})
$$
\begin{array}{lll}
\overline U&\rightarrow&\overline U(1)\\
\uparrow&&\uparrow\\
\overline T&\rightarrow&\overline T(1)
\end{array}
$$
of type (M1) such that there is an expression of the form of (\ref{eq16}),
$$
h_{[\Lambda]}=\overline w_1(1)^{\lambda_1(1)}\cdots\overline w_{\overline s}(1)^{\lambda_{\overline s}(1)}\left(\sum_{i=1}^{\overline b}\overline z_1(1)^{\overline u_{i,1}}\cdots\overline z_{\overline r}(1)^{\overline u_{i,\overline r}}g_i\right),
$$
$d\in{\bf N}$
with $\overline u_i=(\overline u_{i,1},\ldots,\overline u_{i,\overline r})\in\frac{1}{d}{\bf Z}^{\overline r}$
for $1\le i\le \overline b$ , $g_i\in k(\overline U)[\phi_1^{\frac{1}{d}},\ldots,\phi_{\overline r}^{\frac{1}{d}}][[\overline z_1(1),\ldots,\overline z_{\overline r+\overline m}(1)]]$, and
$$
\overline w_1^{\lambda_1}\cdots\overline w_{\overline s}^{\lambda_{\overline s}}
=\overline w_1^{\lambda_1(1)}\cdots\overline w_{\overline s}(1)^{\lambda_{\overline s}(1)}.
$$
Set $L'=k(\overline U)[\phi_1^{\frac{1}{d}},\ldots,\phi_{\overline r}^{\frac{1}{d}}]$.
There exists $a_1,\ldots,a_{\overline r}\in{\bf N}$ such that
$$
a_j+\overline u_{i,j}\ge 0
$$
for $1\le j\le \overline r$ and $1\le i\le \overline b$. Set
$$
\delta=\frac {h_{[\Lambda]}\overline z_1(1)^{a_1}\cdots\overline z_{\overline r}(1)^{a_{\overline r}}}{\overline w_1(1)^{\lambda_1(1)}\cdots\overline w_{\overline s}(1)^{\lambda_{\overline s}(1)}}
\in L'[[\overline z_1(1)^{\frac{1}{d}},\ldots,\overline z_{\overline r}(1)^{\frac{1}{d}},\overline z_{\overline r+1}(1),\ldots,\overline z_{\overline r+\overline m}(1)]]
$$
and let
$$
\sigma=\nu(\overline z_1(1)^{a_1}\cdots\overline z_{\overline r}(1)^{a_{\overline r}}).
$$
Let $\omega$ be a primitive $d$-th root of unity (in an algebraic closure of $L'$).
Set
$$
\begin{array}{ll}
\delta_{j_1,\ldots,j_{\overline r}}&=\delta(\omega^{j_1}\overline z_1(1)^{\frac{1}{d}},
\ldots,\omega^{j_{\overline r}}\overline z_{\overline r}(1)^{\frac{1}{d}},\overline z_{\overline r+1}(1),\ldots,\overline z_{\overline r+\overline m}(1))\\
&\in B_1=L'[\omega][[\overline z_1(1)^{\frac{1}{d}},\ldots,\overline z_{\overline r}(1)^{\frac{1}{d}},
\overline z_{\overline r+1}(1),\ldots,\overline z_{\overline r+\overline m}(1)]].
\end{array}
$$
Set
$$
\epsilon= \prod_{j_1,\ldots,j_{\overline r}=1}^d \delta_{j_1,\ldots,j_{\overline r}}\in A_1=L'[[\overline z_1(1),\ldots,\overline z_{\overline r+\overline m}(1)]].
$$
Identify $\nu^*$ with an extension of (a restriction of $\nu^*$ to $Q(F)$) which dominates $B_1$.

We will now prove that
$$
\nu^*(h_{[\Lambda]})=\infty\Leftrightarrow \nu^*(\delta)=\infty\Leftrightarrow
\nu^*(\epsilon)=\infty.
$$
We certainly have that 
$$
\nu^*(h_{[\Lambda]})=\infty\Leftrightarrow \nu^*(\delta)=\infty.
$$
$\nu^*(\delta)=\infty\Rightarrow \nu^*(\epsilon)=\infty$ since $\delta\mid\epsilon$
in $B_1$.

Suppose that $\nu^*(\epsilon)=\infty$ and $\nu^*(\delta)<\infty$. We will derive a
contradiction. 
$p_{\overline m}$ has the form of (\ref{eq20}) and $\overline T\rightarrow \overline T(1)$
of type $I$ implies 
$q=p_{\overline m}k(\overline T(1))[[\overline z_1(1),\ldots,\overline z_{\overline r+\overline m}(1)]]$
is a prime which is a complete intersection of height $\tilde m=m-\sigma(0)$.
As $q\subset p_{\overline m}(1)$ and since, by assumption, $\sigma(1)=\sigma(0)$, so that $p_{\overline m}(1)$ is a prime
ideal of the same height, we must have $q=p_{\overline m}(1)$.

By (\ref{eq6}),
$$
p_{A_1}=\{a\in A_1\mid\nu^*(a)=\infty\}=p_{\overline m}(1)A_1.
$$
Since $A_1\rightarrow B_1$ is finite and $p_{\overline m}(1)B_1$ is a prime ideal,
$$
p_{B_1}=\{a\in B_1\mid \nu^*(a)=\infty\}=p_{\overline m}(1)B_1.
$$
Since $B_1$ is Galois over $A_1$, the automorphisms of $B_1$ over $A_1$ fix $p_{B_1}$
and $\delta\mid \epsilon$ in $B_1$, so that some conjugate of $\delta$ is in $p_{B_1}$,
we have $\delta\in p_{B_1}$. Thus $\nu^*(\delta)=\infty$, a contradiction.
We have completed the verification that
$$
\nu^*(h_{[\Lambda]})=\infty\Leftrightarrow \nu^*(\delta)\Leftrightarrow \nu^*(\epsilon)=\infty.
$$

We now continue with our proof of the theorem for $h_{[\Lambda]}$ satisfying (\ref{eq42})
or (\ref{eq43}) when $\overline m\le l$.

By (A2) and (A3) of Theorem \ref{Theorem3} and Lemma 4.2 of \cite{C1} there exists a UTS
$\overline T(1)\rightarrow \overline T(2)$ of form $\overline m$ along $\nu$ such that
$$
\epsilon=\overline z_1(2)^{g_1}\cdots\overline z_{\overline r}(2)^{g_{\overline r}}\Sigma
$$
where
$$
\Sigma\in C=L'*k(\overline T(2))[[\overline z_1(2),\ldots,\overline z_{\overline r+\overline m}(2)]]
$$
with $\Sigma$ a unit in $C$ if $\nu^*(h_{[\Lambda}])<\infty$ and
$$
\nu(\overline z_1(2)^{g_1}\cdots\overline z_{\overline r}(2)^{g_{\overline r}})>d^{\overline r}
(A+\sigma)
$$
if $\nu^*(h_{[\Lambda]})=\infty$.

We can further assume that $\overline z_i(2)$ does not divide $\Sigma$ if $1\le i\le \overline r$.
We have expressions
$$
\begin{array}{ll}
\overline z_1(1)&=\overline z_1(2)^{a_{11}}\cdots\overline z_{\overline r}(2)^{a_{1\overline r}}
b_1\\
&\vdots\\
\overline z_{\overline r}(1)&=\overline z_1(2)^{a_{\overline r1}}\cdots\overline z_{\overline r}
(2)^{a_{\overline r\overline r}}b_{\overline r}
\end{array}
$$
with $b_1,\ldots,b_{\overline r}\in k(\overline T(2))$ and there exist polynomials
$$
a_i\in k(\overline T(2))[\overline z_1(2),\ldots,\overline z_{\overline r+\overline m}(2)]
\text{ for }\overline r+1\le i\le \overline r+\overline m
$$
such that $\overline z_i(1)=a_i$ for $\overline r+1\le i\le\overline r+\overline m$.
Thus there exists a series  in indeterminates $x_1,\ldots,x_{2\overline r+\overline  m}$
$$
\delta'\in L'*k(\overline T(2))[[x_1,\ldots,x_{2\overline r+\overline m}]]
$$
 such that
$$
\delta_{j_1\cdots j_{\overline r}}=\delta'(\omega^{j_1}\overline z_1(2)^{\frac{a_{11}}{d}}
\cdots \overline z_{\overline r}
(2)^{\frac{a_{1\overline r}}{d}}b_1^{\frac{1}{d}},\cdots,
\omega^{j_{\overline r}}\overline z_1(2)^{\frac{a_{\overline r1}}{d}}
\cdots \overline z_{\overline r}
(2)^{\frac{a_{\overline r\overline r}}{d}}b_{\overline r}^{\frac{1}{d}},\overline z_{1}(2),\ldots,\overline z_{\overline r+\overline m}(2))
$$
for $1\le j_1,\ldots,j_{\overline r}\le d$. Let
$$
D=L'*k(\overline T(2))[\omega,b_1^{\frac{1}{d^{\overline r}}},\ldots,b_{\overline r}^{\frac{1}{d^{\overline r}}}]
[[\overline z_1(2)^{\frac{1}{d^{\overline r}}},\ldots,\overline z_{\overline r}(2)^{\frac{1}{d^{\overline r}}},
\overline z_{\overline r+1}(2),\ldots,\overline z_{\overline r+\overline m}(2)]].
$$
We have $\delta_{j_1\cdots j_{\overline r}}\in D$ for all $j_1,\ldots, j_{\overline r}$.
Since for  any natural numbers
$a_1,\ldots, a_{\overline r}$ we have 
$$
\overline z_1(2)^{\frac{a_1}{d^{\overline r}}}\cdots\overline z_{\overline r}(2)^{\frac{a_{\overline r}}{d^{\overline r}}}
\mid \delta_{j_1\cdots j_{\overline r}}\Leftrightarrow \overline z_1(2)^{\frac{a_1}{d^{\overline r}}}
\cdots\overline z_{\overline r}(2)^{\frac{a_{\overline r}}{d^{\overline r}}}\mid\delta
$$
in $D$, we have
$$
\overline z_1(2)^{\frac{g_1}{d^{\overline r}}}\cdots\overline z_{\overline r}(2)^{\frac{g_{\overline r}}{d^{\overline r}}}\mid \delta
$$
in $D$, so that we have a factorization
$$
\delta=\overline z_1(2)^{\frac{g_1}{d^{\overline r}}}\cdots\overline z_{\overline r}(2)^{\frac
{g_{\overline r}}{d^{\overline r}}}\delta''
$$
where $\delta''\in D$ is such that $\overline z_i(2)^{\frac{1}{d^{\overline r}}}$ does not divide
$\delta''$ for $1\le i\le \overline r$ and if $\nu^*(h_{[\Lambda]})<\infty$,
$\delta''$ is a unit.

$\overline T(1)\rightarrow \overline T(2)$ extends to a CUTS of form $\overline m$
$$
\begin{array}{lll}
\overline U(1)&\rightarrow&\overline U(2)\\
\uparrow&&\uparrow\\
\overline T(1)&\rightarrow&\overline T(2)
\end{array}
$$
by Lemma \ref{Lemma4.3} and \ref{Lemma4.4}.
$$
\begin{array}{ll}
\overline z_1(1)^{a_1}\cdots\overline z_{\overline r}(1)^{a_{\overline r}}h_{[\Lambda]}
&=\overline w_1(1)^{\lambda_1(1)}\cdots\overline w_{\overline s}(1)^{\lambda_{\overline s}(1)}
\delta\\
&=\overline w_1(1)^{\lambda_1(1)}\cdots\overline w_{\overline s}(1)^{\lambda_{\overline s}(1)}
\overline z_1(2)^{\frac{g_1}{d^{\overline r}}}\cdots\overline z_{\overline r}(2)^{
\frac{g_{\overline r}}{d^{\overline r}}}
\delta''
\end{array}
$$
in 
$$
E=L'*k(\overline U(2))[\omega,b_1^{\frac{1}{d^{\overline r}}},\ldots,b_{\overline r}^{\frac{1}{d^{\overline r}}},
\phi_1(1)^{\frac{1}{d^{\overline r}}},\ldots,\phi_{\overline r}(1)^{\frac{1}{d^{\overline r}}}]
[[\overline w_1(2)^{\frac{1}{d^{\overline r}}},\ldots,\overline w_{\overline r}(2)^{\frac{1}{d^{\overline r}}},
\overline w_{\overline r+1}(2),\ldots,\overline w_{\overline r+\overline m}(2)]].
$$
Since we necessarily have that
$$
\overline z_1(1)^{a_1}\cdots\overline z_{\overline r}(1)^{a_{\overline r}}\mid
\overline w_1(1)^{\lambda_1(1)}\cdots\overline w_{\overline s}(1)^{\lambda_{\overline s}(1)}
\overline z_1(2)^{\frac{g_1}{d^{\overline r}}}\cdots\overline z_{\overline r}(2)^{\frac{g_{\overline r}}{d^{\overline r}}}
$$

in $E$, there exist $\phi\in k(\overline U(2))[\phi_1(1)^{\frac{1}{d^{\overline r}}},\ldots,\phi_{\overline r}(1)
^{\frac{1}{d^{\overline r}}}]$ and $e_1,\ldots,e_{\overline r}\in{\bf N}$ such that

$$
\frac{\overline w_1(1)^{\lambda_1(1)}\cdots\overline w_{\overline s}(1)^{\lambda_{\overline s}(1)}
\overline z_1(2)^{\frac{g_1}{d^{\overline r}}}\cdots\overline z_{\overline r}(2)^{\frac{g_{\overline r}}{d^{\overline r}}}}
{\overline z_1(1)^{a_1}\cdots\overline z_{\overline r}(1)^{a_{\overline r}}}
=\overline w_1(2)^{\frac{e_1}{d^{\overline r}}}\cdots\overline w_{\overline r}(2)^{\frac{e_{\overline r}}{d^{\overline r}}}
\phi.
$$

We thus have that 
$$
f_1=\frac{e_1}{d^{\overline r}},\ldots,f_{\overline r}=\frac{e_{\overline r}}{d^{\overline r}}\in{\bf N}
$$
and $\delta'''=\phi\delta''\in\overline U(2)$ since $h_{[\Lambda]}\in\overline U(2)$.
Thus 
$$
h_{[\Lambda]}=\overline w_1(2)^{f_1}\cdots\overline w_{\overline r}(2)^{f_{\overline r}}
\delta'''
$$
in $\overline U(2)$. If $\nu^*(h_{[\Lambda]})<\infty$, we have that $\delta'''$ is a unit,
and if $\nu^*(h_{[\Lambda]})=\infty$, we have that
$$
\begin{array}{ll}
\nu^*(\overline w_1(2)^{f_1}\cdots\overline w_{\overline r}(2)^{f_{\overline r}})
&=\nu^*(\overline w_1(1)^{\lambda_1(1)}\cdots\overline w_{\overline s}(1)^{\lambda_{\overline s}(1)})-\nu^*(\overline z_1(1)^{a_1}\cdots\overline z_{\overline r}(1)^{a_{\overline r}})\\
&+\nu^*(\overline z_1(2)^{\frac{g_1}{d^{\overline r}}}\cdots\overline z_{\overline r}(2)^{\frac{g_{\overline r}}{d^{\overline r}}})\\
&>\frac{1}{d^{\overline r}}[d^{\overline r}(A+\sigma)]-\sigma=A.
\end{array}
$$ 
This concludes the proof of the analysis of $f$ satisfying (\ref{eq42})
 or (\ref{eq43}) when $\overline m\le l$.

The proof of the theorem when $f$ satisfies (\ref{eq42}) or (\ref{eq43}) with $\overline m> l$ is, with some obvious notational changes,  the same as
Case 2 of pages 59 -61 of \cite{C1}. The induction on line 9 of page 60 \cite{C1}
is now on $\overline m$ in the conclusions of Theorem \ref{Theorem4.7} of this paper.

The proof of the theorem when $f$ satisfies (\ref{eq44})  is the same, with obvious notational changes, and after replacing references to (42), (43) and (44) of \cite{C1} with (\ref{eq42}) and
(\ref{eq43}) and (\ref{eq44}) of this theorem,  as
``the proof when (44) holds'' on pages 61-65 of \cite{C1}.
\end{pf}

\begin{Theorem}\label{Theorem1}
Suppose that $(R,\overline T'',\overline T)$
and $(S,\overline U'',\overline U)$ is a CUTS along $\nu^*$, such that $\overline T''$
contains the subfield $k(c_0)$ for some $c_0\in\overline T''$
and $\overline U''$ contains a subfield isomorphic to $k(\overline U)$.
$\overline T''$ has regular 
parameters $(\overline z_1,\ldots,\overline z_m)$ and $\overline U''$ has regular
parameters $(\overline w_1,\ldots,\overline w_n)$ with
$$
\begin{array}{ll}
\overline z_1&=\overline w_1^{c_{11}}\cdots \overline w_{\overline s}^{c_{1\overline s}}\phi_1\\
&\vdots\\
\overline z_{\overline r}&=\overline w_1^{c_{\overline r1}}\cdots  \overline w_{\overline s}^{c_{\overline r\overline s}}\phi_{\overline r}\\
\overline z_{\overline r+1}&=\overline w_{\overline s+1}\\
&\vdots\\
\overline z_{\overline r+l}&=\overline w_{\overline s+l}
\end{array}
$$

such that $\phi_1,\ldots,\phi_{\overline r}\in k(\overline U)$, $\nu(\overline z_1),\ldots,\nu(\overline z_{\overline r})$ are rationally independent, 
$\nu^*(\overline w_1),\ldots,\nu^*(\overline w_{\overline s})$ are rationally independent
and $(c_{ij})$ has rank $\overline r$. 

Suppose that  $\overline m>l$ and $f\in\overline T''$ is such that

\begin{equation}\label{eq9}
f=P+\overline w_1^{d_1}\cdots\overline w_{\overline s}^{d_{\overline s}}H
\end{equation}
for some powerseries $P\in k(\overline U)[[\overline w_1,\ldots,\overline w_{\overline s+l}]]$,
$$
H=u(\overline w_{\overline s+\overline m}+\overline w_1^{g_1}\cdots\overline w_{\overline s}^{g_{\overline s}}
\Sigma)
$$
where $u\in k(\overline U)[[\overline w_1,\ldots,\overline w_{\overline s+\overline m}]]$
is a unit, $\Sigma\in k(\overline U)[[\overline w_1,\ldots,\overline w_{\overline s+\overline m-1}]]$ and 
$$
\nu^*(\overline w_1^{g_1}\cdots\overline w_{\overline s}^{g_{\overline s}})\ge\nu^*(\overline w_{\overline s+\overline m}).
$$

Then there exists a CUTS of form $\overline m$ along $\nu^*$ $(R,\overline T''(t),\overline T(t))$ and
$(S,\overline U''(t),\overline U(t))$ 
\begin{equation}\label{eq11}
\begin{array}{lllllllll}
\overline U&=&\overline U(0)&\rightarrow &\overline U(1)&\rightarrow&\cdots&\rightarrow
&\overline U(t)\\
&&\uparrow&&\uparrow&&&&\uparrow\\
\overline T&=&\overline T(0)&\rightarrow &\overline T(1)&\rightarrow&\cdots&\rightarrow
&\overline T(t)
\end{array}
\end{equation}

such that 
$\overline T''(i)$
contains a subfield isomorphic to $k(c_0,\ldots,c_i)$
and $\overline U''(i)$ contains a subfield isomorphic to $k(\overline U(i))$.
$\overline T''(i)$ has regular 
parameters $(\overline z_1(i),\ldots,\overline z_m(i))$ and $\overline U''$ has regular
parameters $(\overline w_1(i),\ldots,\overline w_n(i))$ with

$$
\begin{array}{ll}
\overline z_1(i)&=\overline w_1(i)^{c_{11}(i)}\cdots \overline w_{\overline s}(i)^{c_{1\overline s}(i)}\phi_1(i)\\
&\vdots\\
\overline z_{\overline r}(i)&=\overline w_1(i)^{c_{\overline r1}(i)}\cdots  \overline w_{\overline s}(i)^{c_{\overline r\overline s}(i)}\phi_{\overline r}(i)\\
\overline z_{\overline r+1}(i)&=\overline w_{\overline s+1}(i)\\
&\vdots\\
\overline z_{\overline r+l}(i)&=\overline w_{\overline s+l}(i)
\end{array}
$$
such that $\phi_1(i),\ldots,\phi_{\overline r}(i)\in k(\overline U(i))$, $\nu(\overline z_1(i)),\ldots,\nu(\overline z_{\overline r}(i))$ are rationally independent, 
$\nu^*(\overline w_1(i)),\ldots,\nu^*(\overline w_{\overline s}(i))$ are rationally independent
and $(c_{ij}(i))$ has rank $\overline r$ for $1\le i\le t$. 

We further have that
$$
f=\overline P+\overline w_1(t)^{\overline d_1}\cdots\overline w_{\overline s}(t)^{\overline d_{\overline s}}\overline H
$$
with $\overline P\in k(\overline U(t))[\overline w_1(t),\ldots,\overline w_{\overline s+l}(t)]$,
and there exists a finite extension $L$ of the algebraic closure of $k(\overline T(t))$
in $k(\overline U(t))$, and a positive integer $d$,  such that 
$$
\overline P\in L[\overline z_1(t)^{\frac{1}{d}},\ldots,\overline z_{\overline r}(1)^{\frac{1}{d}},
\overline z_{\overline r+1}(1),\ldots,\overline z_{\overline r+l}(1)],
$$
$$
\overline H=\overline u(\overline w_{\overline s+\overline m}(t)+\overline w_1(t)^{\overline g_1}\cdots\overline w_{\overline s}(t)^{\overline g_{\overline s}}
\overline \Sigma)
$$
where $\overline u\in k(\overline U(t))[[\overline w_1(t),\ldots,\overline w_{\overline s+\overline m}(t)]]$
is a unit, $\overline \Sigma\in k(\overline U(t))[[\overline w_1(t),\ldots,\overline w_{\overline s+\overline m-1}(t)]]$ and $\nu^*(\overline w_1(t)^{\overline g_1}\cdots\overline w_{\overline s}(t)^{\overline g_{\overline s}})\ge \nu^*(\overline w_{\overline s+\overline m}(t))$.
\end{Theorem}

\begin{pf} Set $\rho=\nu^*(\overline w_1^{d_1}\cdots\overline w_{\overline s}^{\overline d_{\overline s}})$. There is an expression 
\begin{equation}\label{eq14}
P=\sum_{\Lambda\in{\bf Z}^{\overline s}/({\bf Q}^{\overline r}C)\cap{\bf Z}^{\overline s}}
h_{[\Lambda]}
\end{equation}
of the form (\ref{eq2}). We can rewrite as
$$
P=\sum_{\nu^*(h_{[\Lambda]})\le\rho}h_{[\Lambda]}+\Omega
$$
with $\nu^*(\Omega)>\rho$.  For each of the finitely many $\Lambda$ 
with $\nu^*(h_{[\Lambda]})\le\rho$ (c.f. Lemma 2.3 \cite{C1}), we can further
write
$$
h_{[\Lambda]}=\overline h_{[\Lambda]}+\Omega_{[\Lambda]}
$$
with $\nu^*(\Omega_{[\Lambda]})>\rho$ and $\overline h_{[\Lambda]}\in k(\overline U)
[\overline w_1,\ldots,\overline w_{\overline s+l}]$ of the form (\ref{eq15}). Set
$$
\overline P_1=\sum \overline h_{[\Lambda]},
$$
$$
\overline P_2=\sum \Omega_{[\Lambda]}+\Omega,
$$
$$
P=\overline P_1+\overline P_2.
$$
Observe that $\nu^*(\overline h_{[\Lambda]})=\nu^*(h_{[\Lambda]})\le \rho$ for all 
$\Lambda$ in the (finite) sum $\overline P_1$, and $\nu^*(\overline P_2)>\rho$.

By Theorem \ref{Theorem4.7} applied to $\overline P_2$ in  equations
(\ref{eq43}) or (\ref{eq44}) and by Lemma 4.2 \cite{C1},  there exists a CUTS 
$(R,\overline T''(1),\overline T(1))$ and $(S,\overline U''(1),\overline U(1))$ along
$\nu^*$ of form $l$ 
$$
\begin{array}{lll}
\overline U&\rightarrow&\overline U(1)\\
\uparrow&&\uparrow\\
\overline T&\rightarrow&\overline T(1)
\end{array}
$$
such that $\overline P_2=\overline w_1(1)^{e_1}\cdots\overline w_{\overline s}(1)^{e_{\overline s}}\Phi$ with $\Phi\in k(\overline U(1))[[\overline w_1(1),\ldots,\overline w_{\overline s+l}(1)]]$,
and
$$
\overline w_1^{d_1}\cdots\overline w_{\overline s}^{d_{\overline s}}=
\overline w_1(1)^{d_1(1)}\cdots\overline w_{\overline s}(1)^{d_{\overline s}(1)}
$$
with $e_i>d_i(1)$ for all $i$. Thus
$$
f=\overline P_1+\overline w_1(1)^{d_1(1)}\cdots\overline w_{\overline s}(1)^{d_{\overline s}(1)}
\Omega'
$$
where 
$$
\Omega'\in m(k(\overline U(1))[[\overline w_1(1),\ldots,\overline w_{\overline s+\overline m}(1)]]
$$
 and 
$$
\frac{\partial \Omega'}{\partial w_{\overline s+\overline m}}\not\in m\left(k(\overline U(1))
[[\overline w_1(1),\ldots,\overline w_{\overline s+\overline m}(1)]]\right).
$$
 By the implicit function theorem (the case s=1 of the Weierstrass Preparation Theorem,
Corollary 1 to Theorem 5, Section 1, Chapter VII \cite{ZS}),
$$
\Omega'=u'(\overline w_{\overline s+\overline m}(1)+\overline\Sigma)
$$
with $u'\in k(\overline U(1))[[\overline w_1(1),\ldots,\overline w_{\overline s+\overline m}(1)]]$ a unit, $\overline\Sigma\in k(\overline U(1))[[\overline w_1(1),\ldots,\overline 
w_{\overline s+\overline m-1}(1)]]$.

After replacing $\overline w_{\overline s+\overline m}(1)=\overline w_{\overline s+\overline m}$ with
$\overline w_{\overline s+\overline m}+\Psi$ with 
$$
\Psi\in k(\overline U(1))[\overline w_1(1),\ldots,\overline w_{\overline s+\overline m-1}(1)]
\subset \overline U''(1),
$$
we may assume that
$$
\overline\Sigma\in (\overline w_1(1),\ldots,\overline w_{\overline s+\overline m-1}(1))^B
$$
where $B$ is arbitrarily large. If $\nu^*(\Omega')<\infty$, we can choose $B$ so large that $$
\nu^*(\Omega')=\nu^*(\overline w_{\overline s+\overline m}(1))<\nu^*(\overline\Sigma).
$$
If $\nu^*(\Omega')=\infty$, we have $\nu^*(\overline\Sigma)=\nu^*(\overline w_{\overline s+\overline m}(1))<\infty$.
Then by cases (\ref{eq42}) and (\ref{eq43}) of Theorem \ref{Theorem4.7}, we can perform a CUTS
of the form $\overline m-1$, 
$$
\begin{array}{lll}
\overline U(1)&\rightarrow&\overline U(2)\\
\uparrow&&\uparrow\\
\overline T(1)&\rightarrow&\overline T(2)
\end{array}
$$
to get
$$
f=\overline P_1+\overline w_1(2)^{d_1(2)}\cdots\overline w_{\overline s}^{d_{\overline s}(2)}H'
$$
where
$$
H'=
u'(\overline w_{\overline s+\overline m}(2)+\overline w_1(2)^{g_1(2)}\cdots\overline w_{\overline s}(1)^{g_{\overline s}(2)}\Sigma')
$$
with $u'$ a unit and $\nu^*(\overline w_1(1)^{g_1(2)}\cdots\overline w_{\overline s}(1)^{g_{\overline s}(2)})\ge\nu^*(\overline w_{\overline s+\overline m}(2))$.

Thus we may assume that in (\ref{eq9}), $P\in k(\overline U)[\overline w_1,\ldots,\overline w_{\overline s+l}]$, and 
$$
P=\sum_{[\Lambda]\in {\bf Z}^{\overline s}/{\bf Q}^{\overline r}C\cap {\bf Z}^{\overline s}}h_{[\Lambda]}
$$
with $\nu^*(h_{[\Lambda]})\le \rho$ for all $\Lambda$.

Recall that $\nu^*(f)<\infty$ since $f\in\overline T''$.  If $P= 0$  we have
proven the Theorem. Thus we may assume that $P\ne0$.  There exists $\Lambda_0$ such that $\nu^*(f)=\nu^*(P)=\nu^*(h_{[\Lambda_0]})$ where $\Lambda_0$ is such that
$$
\nu^*(h_{[\Lambda_0]})=\text{ min }\{\nu^*(h_{[\Lambda]})\}
$$
by Remark
\ref{Remark2}. Thus
$[\Lambda_0]=0$ by Lemma \ref{Lemma2}. By Remark \ref{Remark3}, there exists a CUTS of type
(M1),
$$
\begin{array}{lll}
\overline U&\rightarrow&\overline U(1)\\
\uparrow&&\uparrow\\
\overline T&\rightarrow&\overline T(1)
\end{array}
$$
and $d\in{\bf N}$, such that
$$
h_0=h_{[\Lambda_0]}\in k(\overline U)[\phi_1^{\frac{1}{d}},\ldots,\phi_{\overline r}^{\frac{1}{d}}]
[\overline z_1(1)^{\frac{1}{d}},\ldots,\overline z_{\overline r}(1)^{\frac{1}{d}},
\overline z_{\overline r+1}(1),\ldots,\overline z_{\overline r+l}(1)].
$$
Let 
$$
A_1=k(\overline U)[\phi_1^{\frac{1}{d}},\ldots,\phi_{\overline r}^{\frac{1}{d}}]
[[\overline w_1(1)^{\frac{1}{d}},\ldots,\overline w_{\overline s}(1)^{\frac{1}{d}},
\overline w_{\overline s+1}(1),\ldots,\overline w_{n}(1)]].
$$
Extend $\tilde\nu^*$ to the finite extension $Q(A_1)$ so that it
dominates $A_1$. Let 
$$
B_1=\overline T(1)''\otimes_kk(t_1,\ldots,t_{\overline \alpha}).
$$
Let $\nu'=\tilde\nu^*\mid Q(B_1)$. By Lemma \ref{Lemma4}, $\nu'$ is a rank one valuation 
with value group $\Gamma_{\nu}$, since it extends the rank 1 valuation $\nu\mid Q(\overline T(1)'')$.  

Let $k'$ be the algebraic closure of $k$ in $\overline T(1)''$, and let
$$
C=\overline T(1)''\otimes_{k'}k(\overline U(1))[\phi_1^{\frac{1}{d}},\ldots,\phi_{\overline r}^{\frac{1}{d}},\overline z_1^{\frac{1}{d}},\ldots,\overline z_{\overline r}^{\frac{1}{d}}].
$$
Let $\hat \nu$ be the restriction to $Q(C)$ of an extension of $\tilde\nu^*$ to $Q(A_1)$. 
$\hat\nu$ dominates $C$. Since $C$ is finite over $B_1$, $\hat\nu$ has rank 1 and $\Gamma_{\hat\nu}
\subset \Gamma_{\nu}\otimes{\bf Q}$.

If $P\ne h_0$, then there exists $h_{[\Lambda_1]}$ such that $\Lambda_1\not\in {\bf Q}^{\overline r}C\cap{\bf Z}^{\overline s}$ and
$$
\nu^*(f-h_0)=\nu^*(h_{[\Lambda_1]}).
$$
 $f-h_0\in C$ implies 
\begin{equation}\label{eq23}
\nu^*(f-h_0)\in \Gamma_{\nu}\otimes{\bf Q}.
\end{equation}
But 
$$
\nu^*(h_{[\Lambda_1]})\not\in \Gamma_{\nu}\otimes{\bf Q}
$$
by Lemma \ref{Lemma2}, since
 $\Lambda_1\not \in{\bf Q}^{\overline r}C\cap {\bf Z}^{\overline s}$, a contradiction.

Thus we may assume that in (\ref{eq9}), $P\in k(\overline U)[\overline w_1,\ldots,\overline w_{\overline s+l}]$ is such that 
\begin{equation}\label{eq100}
P\in k(\overline U)[\phi_1^{\frac{1}{d}},\ldots,\phi_{\overline r}^{\frac{1}{d}}]
[\overline z_1^{\frac{1}{d}},\ldots,\overline z_{\overline r}^{\frac{1}{d}},
\overline z_{\overline r+1},\ldots,\overline z_{\overline r+l}]
\end{equation} 
and $\nu^*(P)\le\rho$.

Recall that $t_1,\ldots,t_{\overline \alpha}$ is a transcendence basis of $k(\overline U)$ over
$k(\overline T)$.  Write 
\begin{equation}\label{eq101}
P=\sum a_I\overline w_1^{b_1(I)}\cdots\overline w_{\overline s+l}^{b_{\overline s+l}(I)}
\end{equation}
with $a_I\in k(\overline U)$. Let $a_{I_1},\ldots, a_{I_{\gamma}}$ be the (finitely many) nonzero terms in the sum.
Let $A$ be the integral closure of
$$
k(\overline T)[t_1,\ldots,t_{\overline \alpha},a_{I_1},\ldots,a_{I_{\gamma}},
\phi_1,\frac{1}{\phi_1},
\ldots,\phi_{\overline r},\frac{1}{\phi_{\overline r}}]
$$ in $k(\overline U)$. There
exists an algebraic regular local ring $B$ of $k(\overline U)$ such that $B$ dominates $A$,
and the residue field of $B$ is finite over $k(\overline T)$ (c.f. Theorem 2.9 \cite{C1}).

Let $(v_1,\ldots,v_{\overline \alpha})$ be a regular system of parameters in $B$.
We have an inclusion
$$
B\rightarrow \hat B=L_1[[v_1,\ldots,v_{\overline \alpha}]]
$$
where $L_1=k(B)$ is a finite extension of $k(\overline T)$.

After reindexing $\overline w_1,\ldots,\overline w_{\overline s}$, we may assume that the
matrix
$$
\overline C=\left(\begin{array}{lll}
c_{11}&\cdots&c_{1\overline r}\\
\vdots&&\vdots\\
c_{\overline r1}&\cdots&c_{\overline r\overline r}
\end{array}\right)
$$
has positive determinant $e=\text{det}(\overline C)$. Let
 $$
\overline B=(b_{ij})=\frac{1}{de}\text{adj}\overline C=\frac{1}{d}\overline C^{-1}.
$$
Let 
$$
\psi_i=\phi_1^{b_{i1}}\cdots\phi_{\overline r}^{b_{i\overline r}}
$$
for $1\le i\le \overline r$, and let
$$
\tilde w_i=\left\{\begin{array}{ll}\psi_i^d\overline w_i&1\le i\le \overline r\\
\overline w_i&\overline r<i\le \overline s.
\end{array}\right.
$$
We then have equations
$$
\begin{array}{ll}
\overline z_1&=\tilde w_1^{c_{11}}\cdots\tilde w_{\overline s}^{c_{1\overline s}}\\
&\vdots\\
\overline z_{\overline r}&=\tilde w_1^{c_{\overline r1}}\cdots\tilde w_{\overline s}^{c_{\overline r\overline s}}\\
\overline z_{\overline r+1}&=\overline w_{\overline s+1}\\
&\vdots\\
\overline z_{\overline r+l}&=\overline w_{\overline s+l}.
\end{array}
$$
$\psi_i^{de}\in k(\overline U)$ for $1\le i\le \overline r$.
$\psi_i^{de}\in \hat B$ has residue $0\ne\lambda_i\in L_1$.
Let $L''=L_1[\lambda_1^{\frac{1}{de}},\ldots,\lambda_{\overline r}^{\frac{1}{de}}]$,
$E=L''[[v_1,\ldots,v_{\overline \alpha}]]$. 
 Since
$$
\frac{\psi_i^{de}}{\lambda_i}
$$
has residue 1 in $E$, there exists a $de$-th root $\sigma_i$ of 
$$
\frac{\psi_i^{de}}{\lambda_i}
$$
in $E$, with residue 1 in $L''$ for $1\le i\le \overline r$. Thus $\psi_i=\lambda_i^{\frac{1}{de}}\sigma_i\in E$ for $1\le i\le \overline r$.
Note that 
$$
\psi_1^{c_{i1}}\cdots\psi_{\overline r}^{c_{i\overline r}}=\phi_i^{\frac{1}{d}}
$$
for $1\le i\le \overline r$, so that we have natural inclusions
$$
k(\overline U)\subset k(\overline U)[\phi_1^{\frac{1}{d}},\ldots,\phi_{\overline r}^{\frac{1}{d}}]
\subset k(\overline U)[\psi_1,\ldots,\psi_{\overline r}]\subset Q(E).
$$

Let $\tau_1,\ldots,\tau_{\overline \alpha}\in{\bf R}_+$ be rationally independent. Let
$\overline\nu$ be the $L''$-valuation on the quotient field of $E$ defined by $\overline \nu(v_i)=\tau_i$ for $1\le i\le\overline \alpha$. Let 
$$
A_2=Q(E)[\overline w_1^{\frac{1}{d}},\ldots,
\overline w_{\overline s}^{\frac{1}{d}},\overline w_{\overline s+1},\cdots,\overline w_n]
_{(\overline w_1^{\frac{1}{d}},\ldots,
\overline w_{\overline s}^{\frac{1}{d}},\overline w_{\overline s+1},\cdots,\overline w_n)}.
$$
Let $w$ be an extension of $\nu^*$ to $Q(A_2)$ which dominates $A_2$.
Any element of $A_2$ is a quotient of elements of the form $g=\frac{\overline x}{\overline a}$
with $0\ne \overline a\in L''[[v_1,\ldots, v_{\overline\alpha}]]$ and
$$
\begin{array}{ll}
\overline x&\in L''[[v_1,\ldots,v_{\overline\alpha}]][\overline w_1^{\frac{1}{d}},\ldots,
\overline w_{\overline s}^{\frac{1}{d}},\overline w_{\overline s+1},\ldots,\overline w_n]\\
&=L''[\overline w_1^{\frac{1}{d}},\ldots,
\overline w_{\overline s}^{\frac{1}{d}},\overline w_{\overline s+1},\ldots,\overline w_n]
[[v_1,\ldots,v_{\overline\alpha}]].
\end{array}
$$
We have that the value group of $w\mid Q(L''[\overline w_1^{\frac{1}{d}},\ldots,
\overline w_{\overline s}^{\frac{1}{d}},\overline w_{\overline s+1},\ldots,\overline w_n])$
is contained in $\Gamma_{\nu^*}\otimes{\bf Q}$ and thus $w(g)=w(\overline x)\in
\Gamma_{\nu^*}\otimes{\bf Q}$ by Lemma \ref{Lemma25}. We further have that $k(V_w)$ is an algebraic extension of $Q(E)$.
We thus have $\Gamma_w\subset \Gamma_{\nu^*}\otimes{\bf Q}$. Identify  $\overline\nu$ with an extension of $\overline\nu$ to $k(V_w)$. 
Let $\tilde \nu'$ be the composite $w\circ\overline \nu$ of $w$ and $\overline\nu$.
$$
\tilde\nu'(v_1),\ldots,\tilde\nu'(v_{\overline \alpha}),\tilde\nu'(\tilde w_1),\ldots,\tilde \nu'
(\tilde w_{\overline s})
$$
is  a rational basis of $\Gamma_{\tilde\nu'}\otimes{\bf Q}$.  We have 
an equality
$$
Q(E)[\overline w_1^{\frac{1}{d}},\ldots,
\overline w_{\overline s}^{\frac{1}{d}},\overline w_{\overline s+1},\ldots,\overline w_n]
=Q(E)[\tilde w_1^{\frac{1}{d}},\ldots,
\tilde w_{\overline s}^{\frac{1}{d}},\overline w_{\overline s+1},\ldots,\overline w_n].
$$

Let $F=k(\overline T)[\overline z_1,\ldots,\overline z_{m}]$,
$G=L''[\overline z_1^{\frac{1}{d}},\ldots,\overline z_{\overline r}^{\frac{1}{d}},
\overline z_{\overline r+1},\ldots,\overline z_{m}]$. The restriction of $w$
to $Q(F)$ is the restriction of $\tilde\nu$ to $Q(F)$ so $w\mid Q(F)$ has residue field which
is an algebraic extension of $k(\overline T)$.

Now $Q(G)$ is finite over $Q(F)$, so that the restriction of $w$ to $Q(G)$ has a residue field
which is an algebraic extension of $k(\overline T)$.

Thus if $h\in Q(G)$ is such that $w(h)=0$, then if $[h]$ is the residue class of $h$ in
$k(V_w)$, we have that $[h]$ is contained in the algebraic closure $M$ of $k(\overline T)$
in $k(V_w)$, and thus by Lemma 1, Section 11, Chapter VI \cite{ZS}, $\overline \nu([h])=0$,
since $\overline\nu$ is a $k(\overline T)$ valuation.

For $I=(i_1,\ldots,i_{\overline\alpha})\in{\bf N}^{\overline\alpha}$, let $v^I$ denote
$v_1^{i_1}\cdots v_{\overline\alpha}^{i_{\overline\alpha}}$. By (\ref{eq100}) and (\ref{eq101}), there is a series expansion
$$
P=\sum_Ig_Iv^I
$$
with each 
$$
g_I\in L''[\overline z_1^{\frac{1}{d}},\ldots,\overline z_{\overline r}^{\frac{1}{d}},
\overline z_{\overline r+1},\ldots,\overline z_{\overline r+l}].
$$

Suppose that $\tilde\nu'(g_Iv^I)=\tilde\nu'(g_Jv^J)$. Then 
$$
\tilde\nu'(\frac{g_I}{g_J})
=\tilde\nu'(\frac{v^J}{v^I})
$$
so that $w(\frac{g_I}{g_J})=0$. Thus
$$
\tilde\nu'(\frac{g_I}{g_J})=\overline \nu(\frac{g_I}{g_J})=0
$$
since $\frac{g_I}{g_J}\in Q(G)$. We have $\overline\nu(\frac{v^J}{v^I})=0$,
so that $I=J$. Thus 
\begin{equation}\label{eq102}
\tilde\nu'(g_Iv^I)=\tilde\nu'(g_Jv^J)\Leftrightarrow I=J.
\end{equation}

Let $N=k(\overline T)[\overline z_1,\ldots,\overline z_m]$.
We will now establish that if $0\ne h\in Q(N)$, then 
\begin{equation}\label{eq21}
\tilde\nu'(h)\in {\bf Q}\tilde\nu'(\overline z_1)+\cdots+{\bf Q}\tilde\nu'(\overline z_{\overline r}).
\end{equation}

To establish (\ref{eq21}), we first observe that since $\nu(h)<\infty$, there is a UTS $\overline T\rightarrow
\overline T_1$ along $\nu$,  such that 
$\overline T_1$ has regular parameters $(\overline z_1(1),
\ldots,\overline z_m(1))$ and
$$
h=\overline z_1(1)^{a_1}\cdots\overline z_{\overline r}(1)^{a_{\overline r}}u
$$
where $a_1,\ldots,a_{\overline r}\in{\bf Z}$, $u\in \overline T_1$ is a unit, and
$$
\begin{array}{ll}
\overline z_1&=\overline z_1(1)^{b_{11}}\cdots\overline z_{\overline r}(1)^{b_{1\overline r}}c_1\\
&\vdots\\
 \overline z_{\overline r}&=\overline z_1(1)^{b_{\overline r1}}\cdots\overline z_{\overline r}(1)^{b_{\overline r}}c_{\overline r}
\end{array}
$$
with $b_{ij}\in{\bf N}$ for all $i,j$, $\text{Det }(b_{ij})\ne 0$ and $c_1,\ldots,c_{\overline r}\in k(\overline T_1)$.

We identify $\tilde\nu'$ with an extension of $\tilde\nu'$ to $Q(\hat A_2)$ which dominates
$\hat A_2$.
Now $\overline T\rightarrow \overline T_1$ is also a UTS along $\tilde\nu'$ since the
center of $\tilde\nu'$ on a UTS of $\overline T$ must be the center of $\nu$.

We thus have that 
$$
\begin{array}{ll}
\tilde\nu'(h)&=\tilde\nu'(\overline z_1(1)^{a_1}\cdots\overline z_{\overline r}(1)^{a_{\overline r}})+\overline\nu(u)\\
&=\tilde\nu'(\overline z_1(1)^{a_1}\cdots\overline z_{\overline r}(1)^{a_{\overline r}})
\in {\bf Q}\tilde\nu'(\overline z_1)+\cdots+{\bf Q}\tilde\nu'(\overline z_{\overline r})
\end{array}
$$
Here $\overline\nu(u)=0$ since $\nu(u)=0$ and $k(V_{\nu})$ is algebraic over $k$.
We have thus established (\ref{eq21}).

Since $Q(G)$ is a finite extension of $Q(N)$,  we have that
$$
\tilde\nu'(h)\in {\bf Q}\tilde\nu'(\overline z_1)+\cdots+{\bf Q}\tilde\nu'(\overline z_{\overline r})\text{ if }0\ne h\in Q(G).
$$  

Since $L''[\overline z_1^{\frac{1}{d}},\ldots,\overline z_{\overline r}^{\frac{1}{d}},
\overline z_{\overline r+1},\ldots,\overline z_{\overline r+l}][[v_1,\ldots,v_{\overline\alpha}]]$ is a Noetherian ring, there exists $I_0$ such that
$$
\tilde\nu'(g_{I_0}v^{I_0})=\text{ min }\{\tilde\nu'(g_Iv^I)\mid g_I\ne 0\}.
$$
We necessarily have that $\tilde\nu'(P)=\tilde\nu'(g_{I_0}v^{I_0})$
by (\ref{eq102}). Thus 
$$
\tilde\nu'(P)
\in {\bf Q}\tilde \nu'(\overline z_1)+\cdots+{\bf Q}\tilde\nu'(\overline z_{\overline r})
$$
if and only if $I_0=0$.

Write $P=P_1+P_2$ where
$$
P_1=\sum_{w(g_I)\le\rho}g_Iv^I,
$$
$$
P_2=\sum_{w(f_{g_I})>\rho}g_Iv^I
$$
(each sum is possibly infinite). If $P_1\ne0$, let $I_0$ be such that
$\tilde\nu'(t^{I_0}g_{I_0})=\tilde\nu'(P_1)$. Then
$\tilde\nu'(v^{I_0}g_{I_0})=\nu(f)\in\Gamma_{\nu}\otimes{\bf Q}$ implies
$$
\tilde\nu'(v^{I_0})\in{\bf Q}\tilde\nu'(\overline z_1)+\cdots+{\bf Q}
\tilde\nu'(\overline z_{\overline r})\subset \tilde\nu'(\tilde w_1){\bf Q}+\cdots+
\tilde\nu'(\tilde w_s){\bf Q}
$$
 which implies $I_0=0$.
Thus $v^{I_0}g_{I_0}=g_{I_0}\in L''[\overline z_1^{\frac{1}{d}},\ldots,\overline z_{\overline r}^{\frac{1}{d}},\overline z_{\overline r+1},\ldots,\overline z_{\overline r+l}]$. Now
$\nu(f-g_{I_0})\in \Gamma_{\nu}\otimes{\bf Q}$ implies $P_1=g_{I_0}$.

In the sum $P_2$, let $g_1,\ldots,g_{\beta}$ be generators of
the ideal
$$
(g_I\mid P_2=\sum g_Iv^I)\subset L''[\overline z_1^{\frac{1}{d}},
\ldots,\overline z_{\overline r}^{\frac{1}{d}},\overline z_{\overline z+1},\ldots,\overline z_{\overline r+l}].
$$

Let $\omega$ be a primitive $d$th root of unity. Let
$$
d_j=\prod_{i_1,\ldots,i_{\overline r}=1}^d g_j( \omega^{i_1}\overline z_1^{\frac{1}{d}},\ldots, \omega^{i_{\overline r}}\overline z_{\overline r}^{\frac{1}{d}},\overline z_{\overline r+1},\ldots,\overline z_{\overline r+l})\in L''[\overline z_1,\ldots,\overline z_{\overline r+l}]
$$
for $1\le j\le\beta$. The $d_j$ are of the form of  (A2) of Theorem 
\ref{Theorem3} with $\overline m=l$.

Now apply Theorem \ref{Theorem3} to $d_j$ for $1\le j\le\beta$ (and Lemmas
\ref{Lemma4.3} and \ref{Lemma4.4}) to construct a CUTS 
of form $l$ along $\nu^*$,

\begin{equation}\label{eq22}
\begin{array}{lll}
\overline T(1)&\rightarrow &\overline U(1)\\
\uparrow&&\uparrow\\
\overline T&\rightarrow&\overline U
\end{array}
\end{equation}
so that for $1\le j\le \beta$,
$$
d_j=\overline z_1(1)^{e_1(j)}\cdots\overline z_{\overline r}(1)^{e_{\overline r}(j)}u_j
$$
with $u_j$ a unit in $k(T(1))*L''[\overline z_1(1),\ldots,\overline z_{\overline r+l}(1)]$. 
We have 
$$
\begin{array}{ll}
\overline z_1&=\overline z_1(1)^{a_{11}}\cdots\overline z_{\overline r}(1)^{a_{\overline r1}}c_1\\
&\vdots\\
\overline z_r&=\overline z_1(1)^{a_{11}}\cdots\overline z_{\overline r}(1)^{a_{\overline r\overline r}}c_{\overline r}
\end{array}
$$
for some $c_1,\ldots,c_{\overline r}\in k(\overline T(1))$.

We have 
$$
P_1\in D=L''*k(\overline T(1))[c_1^{\frac{1}{d}},\ldots,c_{\overline r}^{\frac{1}{d}}]
[\overline z_1(1)^{\frac{1}{d}},\ldots,\overline z_{\overline r}(1)^{\frac{1}{d}},
\overline z_{\overline r+1}(1),\ldots,\overline z_{\overline r+l}(1)].
$$
$g_{j}$ divides $d_j$ in $D$ implies $g_{j}=\overline w_1(1)^{\overline e_1(j)}\cdots
\overline w_{\overline s}(1)^{\overline e_{\overline s}(j)}f_j
$ 
in 
$$
M=L''*k(\overline U(1))[c_1^{\frac{1}{d}},\ldots,c_{\overline r}^{\frac{1}{d}},
\phi_1(1)^{\frac{1}{d}},\ldots,\phi_{\overline r}(1)^{\frac{1}{d}}]
[\overline w_1(1)^{\frac{1}{d}},\ldots,\overline w_{\overline s}(1)^{\frac{1}{d}},
\overline w_{\overline s+1}(1),\ldots,\overline w_{\overline s+l}(1)]
$$
 for $1\le j\le d$
where $f_j$ is a unit, $\overline e_i(j)\in\frac{1}{d}{\bf N}$, and
$\nu^*(\overline w_1(1)^{\overline e_1(j)}\cdots\overline w_{\overline s}(1)^{\overline e_{\overline s}(j)})>\rho$ for all $1\le j\le\beta$, where $\nu^*$ is identified with an
extension of $\nu^*$ to $Q(M)$ which dominates $M$.

There exists a CUTS of type (M1) along $\nu^*$
$$
\begin{array}{lll}
\overline U(1)&\rightarrow&\overline U(2)\\
\uparrow&&\uparrow\\
\overline T(1)&=&\overline T(2)
\end{array}
$$
so that if
$$
\overline w_1^{d_1}\cdots\overline w_{\overline s}^{d_{\overline s}}=\overline w_1(2)^{d_1(2)}\cdots\overline w_{\overline s}(2)^{d_{\overline s}(2)},
$$
$\overline w_1(2)^{d_1(2)+1}\cdots\overline w_{\overline s}(2)^{d_{\overline s}(2)+1}$
divides $P_2$ in 
$$
F=L''*k(\overline U(2))[c_1^{\frac{1}{d}},\ldots,c_{\overline r}^{\frac{1}{d}},
\phi_1(1)^{\frac{1}{d}},\ldots,\phi_{\overline r}(1)^{\frac{1}{d}}]
[\overline w_1(2)^{\frac{1}{d}},\ldots,\overline w_{\overline s}(2)^{\frac{1}{d}},
\overline w_{\overline s+1}(2),\ldots,\overline w_n(2)]. 
$$
$$
P=P_1+\overline w_1(2)^{d_1(2)+1}\cdots\overline w_{\overline s}(2)^{d_{\overline s}(2)+1}\Phi
$$
with 
$$
P_1\in F'=L''*k(\overline T(2))[c_1^{\frac{1}{d}},\ldots,c_{\overline r}^{\frac{1}{d}}]
[\overline z_1(2)^{\frac{1}{d}},\ldots,\overline z_{\overline r}(2)^{\frac{1}{d}},
\overline z_{\overline r+1}(2),\ldots,\overline z_{\overline r+l}(2)].
$$
We can thus rewrite this equation to get
$$
P=\tilde P_1+\overline w_1(2)^{d_1(2)+1}\cdots\overline w_{\overline s}(2)^{d_{\overline s}(2)+1}
\tilde\Phi
$$
with $\tilde P_1\in F',\tilde \Phi\in F$
and where $\overline w_1(2)^{d_1(2)+1}\cdots\overline w_{\overline s}(2)^{d_{\overline s}(2)+1}$
does not divide any monomial in the expansion of $\tilde P_1$ in $F$.

Comparing with the extension of the expansion of $P$ in 
$k(\overline U(2))[[\overline w_1(2),\ldots,\overline w_{\overline s+l}(2)]]$
to the expansion in  $\hat F$, we see that 
$$
\tilde \Phi\in k(\overline U(2))[[\overline w_1(2),\ldots,\overline w_{\overline s+l}(2)]]\cap F
=k(\overline U(2))[\overline w_1(2),\ldots,\overline w_{\overline s+l}(2)]
$$
and
$$
\tilde P_1\in k(\overline U(2))[[\overline w_1(2),\ldots,\overline w_{\overline s+l}(2)]]\cap F'
\subset L[\overline z_1(2)^{\frac{1}{d}},\ldots,\overline z_{\overline r}(2)^{\frac{1}{d}},
\overline z_{\overline r+1}(2),\ldots,\overline z_{\overline r+l}(2)]
$$
where $L$ is a finite extension of the algebraic closure of $k(\overline T(2))$ in 
$k(\overline U(2))$.
As in the first part of this proof,
we can make a change of variables in $\overline w_{\overline s+\overline m}(2)$ to get the conclusions of the theorem.
\end{pf} 

\section{Conclusion of the proof for rank 1 valuations}
In this section, assumptions and notations will be as in Section \ref{RR1V}.

\begin{Theorem}\label{Theorem4.8}
 Suppose that $T''(0)\subset \hat R$
is a regular local ring, essentially of finite type over $R$ such that the quotient
field of $T''(0)$ is finite over $K$, $U''(0)\subset\hat S$ is a regular local ring,
essentially of finite type over $S$ such that the quotient field of $U''(0)$ is finite over
$K^*$, $U''(0)$ dominates $T''(0)$, $T''(0)$ contains a subfield isomorphic to $k(c_0)$,
for some $c_0\in k(T''(0))$,
$U''(0)$ contains a subfield isomorphic to $k(U''(0))$. Suppose that $R$ has regular
parameters $(x_1,\ldots,x_m)$, $S$ has regular parameters $(y_1,\ldots,y_n)$,
$T''(0)$ has regular parameters $(\tilde{\overline x}_1,\ldots,\tilde{\overline x}_m)$ and
$U''(0)$ has regular parameters $(\tilde{\overline y}_1,\ldots,\tilde{\overline y}_n)$ such
that
$$
\begin{array}{ll}
\tilde{\overline x}_1&=\tilde{\overline y}_1^{c_{11}}\cdots\tilde{\overline y}_{\overline s}^{c_{1\overline s}}\phi_{1}\\
&\vdots\\
\tilde{\overline x}_{\overline r}&=\tilde{\overline y}_1^{c_{\overline r1}}\cdots
\tilde{\overline y}_{\overline s}^{c_{\overline r\overline s}}\phi_{\overline r}\\
\tilde{\overline x}_{\overline r+1}&=\tilde{\overline y}_{\overline s+1}\\
&\vdots\\
\tilde{\overline x}_{\overline r+l}&=\tilde{\overline y}_{\overline s+l}
\end{array}
$$
where $\phi_1,\ldots,\phi_{\overline r}\in k(U''(0))$,  $\nu(\tilde{\overline x}_1),\ldots,
\nu(\tilde{\overline x}_{\overline r})$  are rationally independent, $\nu^*(\tilde{\overline y}_1),\ldots,
\nu^*(\tilde{\overline y}_{\overline s})$ are rationally independent and $(c_{ij})$ has rank $\overline r$. 

Suppose that there exists an algebraic regular local ring $\tilde R\subset R$ such that 
$(x_1,\ldots,x_{\overline r+l})$ are regular parameters in $\tilde R$, $k(\tilde R)\cong
k(c_0)$ and
$$
x_i=\left\{\begin{array}{ll}
\gamma_i\tilde{\overline x}_i&1\le i\le \overline r+l\\
\overline x_i&\overline r+l<i\le m 
\end{array}\right.
$$ 
with $\gamma_i\in k(c_0)[[x_1,\ldots,x_{\overline r+l}]\cap T''(0)$ for 
$1\le i\le \overline r+l$ and $\gamma_i\equiv 1\text{ mod }(x_1,\ldots,x_{\overline r+l})$,
there exist $\gamma_i^y\in U''(0)$ such that $y_i=\gamma_i^y\tilde{\overline y}_i$, $\gamma_i^y\equiv 1
\text{ mod }m(U''(0))$ for $1\le i\le n$.

Suppose that one of the following three conditions holds 

\begin{equation}\label{eq60}
f\in k(U''(0))[[\tilde{\overline y}_1,\ldots,\tilde{\overline y}_{\overline s+\overline m}]]
\text{ for some }\overline m\text{ with }l\le\overline m\le n-\overline
 s\text{ and }\nu^*(f)<\infty.
\end{equation}

\begin{equation}\label{eq61}
f\in k(U''(0))[[\tilde{\overline y}_1,\ldots,\tilde{\overline y}_{\overline s+\overline m}]]
\text{ for some }\overline m\text{ with }l<\overline m\le n-\overline s, \nu^*(f)=\infty,\text{ and }A\in{\bf N}
\text{ is given}.
\end{equation}

\begin{equation}\label{eq62}
f\in U''(0)-k(U''(0))[[\tilde{\overline y}_1,\ldots,\tilde{\overline y}_{\overline s+\overline l}]].
\end{equation}

Then there exists a positive integer $N_0$ such that for $N\ge N_0$, we can construct a
CRUTS along $\nu^*$ $(R,T''(t),T(t))$ and $(S,U''(t),U(t))$ with associated MTSs
$$
\begin{array}{lll}
S&\rightarrow&S(t)\\
\uparrow&&\uparrow\\
R&\rightarrow&R(t)
\end{array}
$$
such that the following holds. $T''(t)$ contains a subfield isomorphic to $k(c_0,\ldots,c_t)$,
$U''(t)$ contains a subfield isomorphic to $k(U(t))$,
$R(t)$ has regular parameters $(x_1(t),\ldots,x_m(t))$, $T''(t)$ has regular parameters
$(\tilde{\overline x}_1(t),\ldots,\tilde{\overline x}_m(t))$, $S(t)$ has regular parameters
$(y_1(1),\ldots,y_n(t))$, $U''(t)$ has regular parameters $(\tilde{\overline y}_1(t),\ldots
\tilde{\overline y}_n(t))$ such that
$$
x_i(t)=\left\{\begin{array}{ll}
\gamma_i(t)\tilde{\overline x}_i(t)&1\le i\le \overline r+l\\
\tilde{\overline x}_i(t)&\overline r+l<i\le m
\end{array}\right.
$$
where $\gamma_i(t)\in k(c_0,\ldots,c_t)[[x_1(t),\ldots,x_{\overline r+l}(t)]]$ are units such that
$$
\gamma_i(t)\equiv 1\text{ mod }(x_1(t),\ldots,x_{\overline r+l}(t)).
$$
In particular,
$$
k(c_0,\ldots,c_t)[[x_1(t),\ldots,x_{\overline r+l}(t)]]=k(c_0,\ldots,c_t)[[\tilde{\overline x}_1(t),\ldots,
\tilde{\overline x}_{\overline r+l}(t)]].
$$
For $1\le i\le n$ there exists $\gamma_i^y(t)\in U''(t)$ such that $y_i(t_=\gamma_i^y(t)
\tilde{\overline y}_i(t)$,
$$
\gamma_i^y(t)\equiv 1\text{ mod } m(U''(t)).
$$

\begin{equation}\label{eq63}
\begin{array}{ll}
\tilde{\overline x}_1(t)&=\tilde{\overline y}_1(t)^{c_{11}(t)}\cdots
\tilde{\overline y}_{\overline s}(t)^{c_{1\overline s}(t)}\phi_1(t)\\
&\vdots\\
\tilde{\overline x}_{\overline r}(t)&=\tilde{\overline y}_1(t)^{c_{\overline r1}}
\cdots\tilde{\overline y}_{\overline s}^{c_{\overline r\overline s}}\phi_{\overline r}(t)\\
\tilde{\overline x}_{\overline r+1}(t)&=\tilde{\overline y}_{\overline s+1}(t)\\
&\vdots\\
\tilde{\overline x}_{\overline r+l}&=\tilde{\overline y}_{\overline s+l}(t)
\end{array}
\end{equation}

$\phi_1(t),\ldots\phi_{\overline r}(t)\in k(U(t))$, $\nu(\tilde{\overline x}_1(t)),\ldots,
\nu(\tilde{\overline x}_{\overline r}(t))$ are rationally independent, 
\linebreak $\nu^*(\tilde{\overline y}_1(t)),\ldots,\nu^*(\tilde{\overline y}_{\overline s}(t))$
are rationally independent and $(c_{ij}(t))$ has rank $\overline r$. There exists an
algebraic regular
local ring $\tilde R(t)\subset R(t)$ such that $(x_1(t),\ldots,x_{\overline r+l})$ are regular
parameters in $\tilde R(t)$ and $k(\tilde R(t))\cong k(c_0,\ldots,c_t)$. Furthermore,
$x_i(t)=x_i$ for $\overline r+l+1\le i\le m$, $y_i(t)=y_i$ for $\overline s+\overline m+1\le
i\le n$, so that the CRUTS is of the form $\overline m$  where 
$\overline s+\overline m=n$ in case (\ref{eq62}). Set
$n_{t,l}=m(k(U(t))[[\tilde{\overline y}_1(t),\ldots,\tilde{\overline y}_{\overline s+l}(t)]])$.

In case (\ref{eq60}) we have 
\begin{equation}\label{eq64}
f\equiv \tilde{\overline y}_1(t)^{d_1}\cdots\tilde{\overline y}_{\overline s}(t)^{ d_{\overline s}}u\text{ mod }m(U(t))^N
\end{equation}
where $u\in k(U(t))[[\tilde{\overline y}_1(t),\ldots,\tilde{\overline y}_{\overline s+\overline m}(t)]]$ is a unit power series. Further, if $f\in k(U)[[\tilde{\overline y}_1,\ldots,
\tilde{\overline y}_{\overline s+l}]]$, then
$$
f\equiv \tilde{\overline y}_1(t)^{d_1}\cdots\tilde{\overline y}_{\overline s}(t)^{ d_{\overline s}}u\text{ mod }n_{t,l}^N
$$
where $u\in k(U(t))[[\tilde{\overline y}_1(t),\ldots,\tilde{\overline y}_{\overline s+\overline l}(t)]]$ is a unit power series.

In case (\ref{eq61}) we have 
\begin{equation}\label{eq65}
f\equiv \tilde{\overline y}_1(t)^{d_1}\cdots\tilde{\overline y}_{\overline s}(t)^{ d_{\overline s}}\Sigma\text{ mod }m(U(t))^N
\end{equation}
where $\Sigma\in k(U(t))[[\tilde{\overline y}_1(t),\ldots,\tilde{\overline y}_{\overline s+\overline m}(t)]]$ and $\nu^*(\tilde{\overline y}_1(t)^{d_1}\cdots\tilde{\overline y}_{\overline s}^{d_{\overline s}}(t))>A$. 
Further, if $f\in k(U)[[\tilde{\overline y}_1,\ldots,
\tilde{\overline y}_{\overline s+l}]]$, then
$$
f\equiv \tilde{\overline y}_1(t)^{d_1}\cdots\tilde{\overline y}_{\overline s}(t)^{ d_{\overline s}}\Sigma\text{ mod }n_{t,l}^N
$$
where $\Sigma\in k(U(t))[[\tilde{\overline y}_1(t),\ldots,\tilde{\overline y}_{\overline s+\overline l}(t)]]$ and $\nu^*(\tilde{\overline y}_1(t)^{d_1}\cdots\tilde{\overline y}_{\overline s}^{d_{\overline s}}(t))>A$.

In case (\ref{eq62}) we have 
\begin{equation}\label{eq66}
f\equiv P+\tilde{\overline y}_1(t)^{d_1}\cdots\tilde{\overline y}_{\overline s}(t)^{d_{\overline s}}H\text{ mod }m(U(t))^N
\end{equation}
where $P\in k(U(t))[\tilde{\overline y}_1(t),\ldots\tilde{\overline y}_{\overline s+l}(t)]$
and there exists a finite extension $L$ of the algebraic closure of $k(T(t))$ in $k( U(t))$ and a positive integer $d$ such that
$$
P\in L[\tilde{\overline x}_1(t)^{\frac{1}{d}},\ldots,\tilde{\overline x}_{\overline r}(t)^{\frac{1}{d}},\tilde{\overline x}_{\overline r+1}(t),\ldots,\tilde{\overline x}_{\overline r+l}(t)],
$$ 
$$
H=u(\tilde{\overline y}_{\overline s+l+1}(t)+\tilde{\overline y}_1(t)^{g_1}\cdots
\tilde{\overline y}_{\overline s}(t)^{g_{\overline s}}\Sigma)
$$
where $u\in U(t)$ is a unit series, $\Sigma\in k(U(t))[[\tilde{\overline y}_1(t),\ldots
\tilde{\overline y}_{\overline s+l}(t),\tilde{\overline y}_{\overline s+l+2}(t),\ldots,
\tilde{\overline y}_n(t)]]$ and
$\nu^*(\tilde{\overline y}_{\overline s+l+1}(t))\le\nu^*(\tilde{\overline y}_1(t)^{g_1}\cdots
\tilde{\overline y}_{\overline s}(t)^{g_{\overline s}})$.

We further have that $\nu^*(m(U(t))$ is a constant which is independent of $N$ for $N\ge N_0$. 
\end{Theorem}

\begin{pf} The proof is essentially the same as the proof of Theorem 4.8 \cite{C1}, with
some modification of notation. On page 67, line 2 of the proof, ``By Theorem 4.7 there is
a CUTS'' should be replaced with ``By Theorems \ref{Theorem4.7} and \ref{Theorem1}
there is a CUTS''. On page 67, line 7 of the proof in \cite{C1}, ``notation of Theorem 4.7''
should be ``notation of Theorems \ref{Theorem4.7} and \ref{Theorem1}''.
On page 67, line 18 of the proof, replace ``$P\in k(\overline U(t))[[\overline w_1(t),\ldots,\overline w_l(t)]]$'' with ``$P\in k(\overline U(t))[\overline w_1(t),\ldots,\overline w_{\overline s+l}(t)]$ and there exists a finite extension $L$ of the algebraic closure of $k(\overline T(t))$ in $k(\overline U(t))$, and a positive integer $d$ such that
$$
P\in L[\overline z_1(t)^{\frac{1}{d}},\ldots,\overline z_{\overline r}(t)^{\frac{1}{d}},
\overline z_{\overline r+1}(1),\ldots,\overline z_{\overline r+l}(1)].\text{''}
$$

  All later references in this proof to Theorem 4.7 and to equations
(46), (48), (49), (50) and (51) should be replaced with references to Theorem
\ref{Theorem4.7}, (\ref{eq46}),
(\ref{eq48}), (\ref{eq49}), (\ref{eq50}) and (\ref{eq51}) of this paper. References to Lemma 4.4 should
be replaced with references to Lemma \ref{Lemma4.4} of this paper.

The independence of $\nu^*(m(U(t))=\text{min}\{\nu^*(f)\mid f\in m(U(t))\}$ of $N$ follows
from (A3) of page 83 of the proof of Theorem 4.8 in \cite{C1}.
\end{pf}

\begin{Theorem}\label{Theorem4.9}
 Suppose that $T''(0)\subset \hat R$
is a regular local ring, essentially of finite type over $R$ such that the quotient
field of $T''(0)$ is finite over $K$, $U''(0)\subset\hat S$ is a regular local ring,
essentially of finite type over $S$ such that the quotient field of $U''(0)$ is finite over
$K^*$, $U''(0)$ dominates $T''(0)$, $T''(0)$ contains a subfield isomorphic to $k(c_0)$
for some $c_0\in k(T''(0))$,
$U''(0)$ contains a subfield isomorphic to $k(U''(0))$. Suppose that $R$ has regular
parameters $(x_1,\ldots,x_m)$, $S$ has regular parameters $(y_1,\ldots,y_n)$,
$T''(0)$ has regular parameters $({\overline x}_1,\ldots,{\overline x}_m)$ and
$U''(0)$ has regular parameters $({\overline y}_1,\ldots,{\overline y}_n)$ such
that
$$
\begin{array}{ll}
{\overline x}_1&={\overline y}_1^{c_{11}}\cdots{\overline y}_{\overline s}^{c_{1\overline s}}\phi_{ 1}\\
&\vdots\\
{\overline x}_{\overline r}&={\overline y}_1^{c_{\overline r1}}\cdots
{\overline y}_{\overline s}^{c_{\overline r\overline s}}\phi_{\overline r}\\
{\overline x}_{\overline r+1}&={\overline y}_{\overline s+1}\\
&\vdots\\
{\overline x}_{\overline r+l}&={\overline y}_{\overline s+l}
\end{array}
$$
where $\phi_1,\ldots,\phi_{\overline r}\in k(U''(0))$,  $\nu({\overline x}_1),\ldots,
\nu({\overline x}_{\overline r})$  are rationally independent, $\nu^*({\overline y}_1),\ldots,
\nu^*({\overline y}_{\overline s})$ are rationally independent and $(c_{ij})$ has rank $\overline r$. 

Suppose that there exists an algebraic regular local ring $\tilde R\subset R$ such that 
$(x_1,\ldots,x_{\overline r+l})$ are regular parameters in $\tilde R$, $k(\tilde R)\cong
k(c_0)$ and
$$
x_i=\left\{\begin{array}{ll}
\gamma_i{\overline x}_i&1\le i\le \overline r+l\\
\overline x_i&\overline r+l<i\le m
\end{array}\right.
$$ 
with $\gamma_i\in k(c_0)[[x_1,\ldots,x_{\overline r+l}]\cap T''(0)$ for 
$1\le i\le \overline r+l$ and $\gamma_i\equiv 1\text{ mod }(x_1,\ldots,x_{\overline r+l})$,
there exist $\gamma_i^y\in U''(0)$ such that $y_i=\gamma_i^y{\overline y}_i$, $\gamma_i^y\equiv 1
\text{ mod }m(U''(0))$ for $1\le i\le n$.

Further, suppose that
$$
x_{\overline r+l+1}=\overline P
+\overline y_1^{d_1}\cdots \overline y_{\overline s}^{d_{\overline s}}\overline H+\Omega
$$
where $\overline P\in k(U''(0))[\overline y_1,\ldots,\overline y_{\overline s+l}]$
and $\overline P\in L[\overline x_1^{\frac{1}{d}},\ldots,\overline x_{\overline r}^{\frac{1}{d}},
\overline x_{\overline r+1},\ldots,\overline x_{\overline s+l}]$, where 
$d$ is a positive integer, $L$ is a finite extension of the algebraic closure of $k(T''(0))$ in $k(U''(0))$,
$$
\overline H=\overline u(\overline y_{\overline s+l+1}+\overline y_1^{\overline g_1}\cdots
\overline y_{\overline s}^{\overline g_{\overline s}}\overline\Sigma)
$$
where $\overline u\in U''(0)\sphat$ is a unit,
$$
\overline\Sigma\in k(U''(0))[[\overline y_1,\ldots,\overline y_{\overline s+l},
\overline y_{\overline s+l+2},\ldots,\overline y_n]],
$$
$$
\nu^*(\overline y_{\overline s+l+1})\le\nu(\overline y_1^{\overline g_1}\cdots\overline y_{\overline s}^{\overline g_{\overline s}})\text{ and}
$$
$$
\Omega\in m(U''(0))^N\text{ with }N\nu^*(m(U''(0)))>\nu^*(\overline y_1^{d_1}\cdots
\overline y_{\overline s}^{d_{\overline s}}\overline y_{\overline s+l+1}).
$$

Then there exists a 
CRUTS along $\nu$ $(R,T''(t),T(t))$ and $(S,U''(t),U(t))$ with associated MTSs
$$
\begin{array}{lll}
S&\rightarrow&S(t')\\
\uparrow&&\uparrow\\
R&\rightarrow&R(t')
\end{array}
$$
such that the following holds. $T''(t')$ contains a subfield isomorphic to $k(c_0,\ldots,c_{t'})$,
$U''(t')$ contains a subfield isomorphic to $k(U(t'))$,
$R(t')$ has regular parameters $(x_1(t'),\ldots,x_m(t'))$, $T''(t')$ has regular parameters
$(\tilde{\overline x}_1(t'),\ldots,\tilde{\overline x}_m(t'))$, $S(t')$ has regular parameters
$(y_1(t'),\ldots,y_n(t'))$, $U''(t')$ has regular parameters $(\tilde{\overline y}_1(t'),\ldots
\tilde{\overline y}_n(t'))$  where

$$
\begin{array}{ll}
\tilde{\overline x}_1(t')&=\tilde{\overline y}_1(t')^{c_{11}(t')}\cdots
\tilde{\overline y}_{\overline s}(t')^{c_{1\overline s}(t')}\phi_1(t')\\
&\vdots\\
\tilde{\overline x}_{\overline r}(t')&=\tilde{\overline y}_1(t')^{c_{\overline r1}(t')}
\cdots\tilde{\overline y}_{\overline s}(t')^{c_{\overline r\overline s}}\phi_{\overline r}(t')\\
\tilde{\overline x}_{\overline r+1}(t')&=\tilde{\overline y}_{\overline s+1}(t')\\
&\vdots\\
\tilde{\overline x}_{\overline r+l}(t')&=\tilde{\overline y}_{\overline s+l}(t')\\
x_{\overline r+l+1}(t')=\tilde{\overline x}_{\overline r+l+1}(t')&=P+
\tilde{\overline y}_1(t')^{d_1(t')}\cdots
\tilde{\overline y}_{\overline s}(t')^{d_{\overline s}(t')}H
\end{array}
$$
where 
$P\in k(U(t'))[\tilde{\overline y}_1(t'),\ldots,\tilde{\overline y}_{\overline s+l}(t')]$
and 
$$
\overline P\in L'[\tilde{\overline x}_1(t')^{\frac{1}{d}},\ldots,\tilde{\overline x}_{\overline r}(t')^{\frac{1}{d}},
\tilde{\overline x}_{\overline r+1}(t'),\ldots,\tilde{\overline x}_{\overline r+l}(t')],
$$
 where $L'$ is a finite extension of the algebraic closure of $k(T(t'))$ in $k(U(t'))$,
$H\in k(U(t'))[[\tilde{\overline y}_1(t'),\ldots,\tilde{\overline y}_{n}(t')]]$ is such that
$$
\text{mult }H(0,\ldots,0,\tilde{\overline y}_{\overline s+l+1}(t'),0,\ldots,0)=1,
$$
$\phi_1(t'),\ldots\phi_{\overline r}(t')\in k(U(t'))$, $\nu(\tilde{\overline x}_1(t')),\ldots,
\nu(\tilde{\overline x}_{\overline r}(t'))$ are rationally independent,
\linebreak 
$\nu^*(\tilde{\overline y}_1(t')),\ldots,\nu^*(\tilde{\overline y}_{\overline s}(t'))$
are rationally independent and $(c_{ij}(t'))$ has rank $\overline r$. There exists an
algebraic regular
local ring $\tilde R(t')\subset R(t')$ such that $(x_1(t'),\ldots,x_{\overline r+l}(t'))$ are regular
parameters in $\tilde R(t')$ and $k(\tilde R(t'))\cong k(c_0,\ldots,c_{t'})$.
$$
x_i(t')=\left\{\begin{array}{ll}
\gamma_i(t')\tilde{\overline x}_i(t')&1\le i\le \overline r+l\\
\tilde{\overline x}_i(t')&\overline r+l<i\le m
\end{array}\right.
$$
with $\gamma_i(t')\in k(c_0,\ldots,c_{t'})[[x_1(t'),\ldots,x_{\overline r+l}(t')]]\cap T''(t')$
units for $1\le i\le\overline r+l$, such that 
$$
\gamma_i(t')\equiv 1\text{ mod }(x_1(t'),\ldots,x_{\overline r+l}(t'))
$$
and for $1\le i\le n$ there exists $\gamma_i^y(t')\in U''(t')$ such that $y_i(t')=\gamma_i^y(t')
\tilde{\overline y}_i(t')$,
$$
\gamma_i^y(t')\equiv 1\text{ mod } m(U''(t')).
$$
\end{Theorem}

The proof of Theorem \ref{Theorem4.9} is similar to the proof of Theorem 4.9 of \cite{C1}.

\begin{Theorem}\label{Theorem4.10} Let $n_{0,l}=m(k(U''(0))[[\overline y_1,\ldots,
\overline y_{\overline s+l}]]$ in the assumptions of Theorem \ref{Theorem4.9}.
\begin{enumerate}
\item If $\Omega\in n_{0,l}^N$ in the assumptions of Theorem \ref{Theorem4.9}, then a sequence
of MTSs of type (M2) and a MTS of type (M1) (so that the CRUTS along $\nu^*$ is of form $l$)
 are sufficient to transform $x_{\overline r+l+1}$ into the form of the
conclusions of Theorem \ref{Theorem4.9}.
\item Suppose that
$$
g=\overline y_1^{d_1}\cdots\overline y_{\overline s}^{d_{\overline s}}u+\Omega
$$
where $u\in k(U''(0))[[\overline y_1,\ldots,\overline y_{\overline s+l}]]$ is a unit
power series and $\Omega\in n_{0,l}^N$ with $N\nu^*(n_{0,l})>\nu(\overline y_1^{d_1}\cdots\overline y_{\overline s}^{d_{\overline s}})$. Then a sequence of MTSs
of type (M2) and a MTS of type (M1) (so that the CRUTS along $\nu^*$ is of form $l$)
 are sufficient to transform $g$ into the form
$$
g=\overline y_1(t')^{d_1(t')}\cdots\overline y_{\overline s}(t')^{d_{\overline s}(t')}
\overline u
$$
where $\overline u\in k(U''(t'))[[\overline y_1(t'),\ldots\overline y_{\overline s+l}(t')]]$ is a unit power series.
\item  Suppose that
$$
g=\overline y_1^{d_1}\cdots\overline y_{\overline s}^{d_{\overline s}}\Sigma+\Omega
$$
where $\Sigma\in k(U''(0))[[\overline y_1,\ldots,\overline y_{\overline s+l}]]$,
$\nu^*(\overline y_1^{d_1}\cdots\overline y_{\overline s}^{d_{\overline s}})>A$ and
$\Omega\in n_{0,l}^N$ with $N\nu^*(n_{0,l})>\nu(\overline y_1^{d_1}\cdots\overline y_{\overline s}^{d_{\overline s}})$. Then a sequence of MTSs
of type (M2) and a MTS of type (M1) (so that the CRUTS along $\nu^*$ is of form $l$)
 are sufficient to transform $g$ into the form
$$
g=\overline y_1(t')^{d_1(t')}\cdots\overline y_{\overline s}(t')^{d_{\overline s}(t')}
\overline \Sigma
$$
where $\overline \Sigma\in k(U''(t'))[[\overline y_1(t'),\ldots\overline y_{\overline s+l}]]$ and
$\nu^*(\overline y_1(t')^{d_1(t')}\cdots\overline y_{\overline s}(t')^{d_{\overline s}(t')})>A$
\item Suppose that
$$
g=\overline y_1^{d_1}\cdots\overline y_{\overline s}^{d_{\overline s}}u+\Omega
$$
where $u\in k(U''(0))[[\overline y_1,\ldots,\overline y_{\overline s+l}]]$ is a unit
power series and $\Omega\in m(U(0))^N$ with $N\nu^*(m(U''(0)))>\nu^*(\overline y_1^{d_1}\cdots\overline y_{\overline s}^{d_{\overline s}})$. Then  there exists a CRUTS along $\nu^*$ as in the conclusions of Theorem \ref{Theorem4.9} such that 
$$
g=\overline y_1(t')^{d_1(t')}\cdots\overline y_{\overline s}(t')^{d_{\overline s}(t')}
\overline u
$$
where $\overline u\in k(U(t'))[[\overline y_1(t'),\ldots\overline y_{\overline s+l}]]$ is a unit power series.
\item Suppose that
$$
g=\overline y_1^{d_1}\cdots\overline y_{\overline s}^{d_{\overline s}}\Sigma+\Omega
$$
where $\Sigma\in k(U''(0))[[\overline y_1,\ldots,\overline y_{\overline s+l}]]$, $\nu^*(\overline y_1^{d_1}\cdots\overline y_{\overline s}^{d_{\overline s}})>A$  and $\Omega\in m(U(0))^N$ with $N\nu^*(m(U''(0)))>\nu^*(\overline y_1^{d_1}\cdots\overline y_{\overline s}^{d_{\overline s}})$. Then  there exists a CRUTS along $\nu^*$ as in the conclusions of Theorem \ref{Theorem4.9} such that 
$$
g=\overline y_1(t')^{d_1(t')}\cdots\overline y_{\overline s}(t')^{d_{\overline s}(t')}
\overline \Sigma
$$
where $\overline \Sigma\in k(U(t'))[[\overline y_1(t'),\ldots\overline y_{\overline s+l}]]$ and
$\nu^*(\overline y_1(t')^{d_1(t')}\cdots\overline y_{\overline s}(t')^{d_{\overline s}(t')})>A$.
\end{enumerate}
\end{Theorem}

Theorem \ref{Theorem4.10} and its proof is a modification of the statement and proof of
Theorem 4.10 of \cite{C1}.

\begin{Theorem}\label{Theorem4.11}
Suppose that $T''(0)\subset \hat R$
is a regular local ring, essentially of finite type over $R$ such that the quotient
field of $T''(0)$ is finite over $K$, $U''(0)\subset\hat S$ is a regular local ring,
essentially of finite type over $S$ such that the quotient field of $U''(0)$ is finite over
$K^*$, $U''(0)$ dominates $T''(0)$, $T''(0)$ contains a subfield isomorphic to $k(c_0)$
for some $c_0\in k(T''(0))$,
$U''(0)$ contains a subfield isomorphic to $k(U''(0))$. Suppose that $R$ has regular
parameters $(x_1,\ldots,x_m)$, $S$ has regular parameters $(y_1,\ldots,y_n)$,
$T''(0)$ has regular parameters $({\overline x}_1,\ldots,{\overline x}_m)$ and
$U''(0)$ has regular parameters $({\overline y}_1,\ldots,{\overline y}_n)$ such
that
$$
\begin{array}{ll}
{\overline x}_1&=\tilde{\overline y}_1^{c_{11}}\cdots{\overline y}_{\overline s}^{c_{1\overline s}}\phi_{ 1}\\
&\vdots\\
{\overline x}_{\overline r}&={\overline y}_1^{c_{\overline r1}}\cdots
{\overline y}_{\overline s}^{c_{\overline r\overline s}}\phi_{\overline r}\\
{\overline x}_{\overline r+1}&={\overline y}_{\overline s+1}\\
&\vdots\\
{\overline x}_{\overline r+l}&={\overline y}_{\overline s+l}
\end{array}
$$
where $\phi_1,\ldots,\phi_{\overline r}\in k(U''(0))$,  $\nu({\overline x}_1),\ldots,
\nu({\overline x}_{\overline r})$  are rationally independent, $\nu^*({\overline y}_1),\ldots,
\nu^*({\overline y}_{\overline s})$ are rationally independent and $(c_{ij})$ has rank $\overline r$. Further suppose that $l<m-\overline r$.

Suppose that there exists an algebraic regular local ring $\tilde R\subset R$ such that 
$(x_1,\ldots,x_{\overline r+l})$ are regular parameters in $\tilde R$, $k(\tilde R)\cong
k(c_0)$ and
$$
x_i=\left\{\begin{array}{ll}
\gamma_i{\overline x}_i&1\le i\le \overline r+l\\
\overline x_i&\overline r+l<i\le m 
\end{array}\right.
$$ 
with $\gamma_i\in k(c_0)[[x_1,\ldots,x_{\overline r+l}]]\cap T''(0)$ for 
$1\le i\le \overline r+l$ and $\gamma_i\equiv 1\text{ mod }(x_1,\ldots,x_{\overline r+l})$,
there exist $\gamma_i^y\in U''(0)$ such that $y_i=\gamma_i^y{\overline y}_i$, $\gamma_i^y\equiv 1
\text{ mod }m(U''(0))$ for $1\le i\le n$.

Then there exists a 
CRUTS along $\nu^*$ $(R,T''(t),T(t))$ and $(S,U''(t),U(t))$ with associated MTSs
$$
\begin{array}{lll}
S&\rightarrow&S(t)\\
\uparrow&&\uparrow\\
R&\rightarrow&R(t)
\end{array}
$$
such that the following holds. $T''(t)$ contains a subfield isomorphic to $k(c_0,\ldots,c_t)$,
$U''(t)$ contains a subfield isomorphic to $k(U(t))$,
$R(t)$ has regular parameters $(x_1(t),\ldots,x_m(t))$, 
$S(t)$ has regular parameters $(y_1(t),\ldots,y_n(t))$, $T''(t)$ has regular parameters
$(\tilde{\overline x}_1(t),\ldots,\tilde{\overline x}_m(t))$,
 $U''(t)$ has regular parameters $(\tilde{\overline y}_1(t),\ldots
\tilde{\overline y}_n(t))$  where

$$
\begin{array}{ll}
\tilde{\overline x}_1(t)&=\tilde{\overline y}_1(t)^{c_{11}(t)}\cdots
\tilde{\overline y}_{\overline s}(t)^{c_{1\overline s}(t)}\phi_1(t)\\
&\vdots\\
\tilde{\overline x}_{\overline r}(t)&=\tilde{\overline y}_1(t)^{c_{\overline r1}(t)}
\cdots\tilde{\overline y}_{\overline s}(t)^{c_{\overline r\overline s}(t)}\phi_{\overline r}(t)\\
\tilde{\overline x}_{\overline r+l}(t)&=\tilde{\overline y}_{\overline s+1}(t)\\
&\vdots\\
\tilde{\overline x}_{\overline r+l+1}(t)&=\tilde{\overline y}_{\overline s+l+1}(t)
\end{array}
$$
such that $\phi_1(t),\ldots\phi_{\overline r}(t)\in k(U(t))$, $\nu(\tilde{\overline x}_1(t)),\ldots,
\nu(\tilde{\overline x}_{\overline r}(t))$ are rationally independent, 
$\nu^*(\tilde{\overline y}_1(t)),\ldots,\nu^*(\tilde{\overline y}_{\overline s}(t))$
are rationally independent and $(c_{ij}(t))$ has rank $\overline r$.
$$
x_i(t)=\tilde{\overline x}_i(t)\text{ for }1\le i\le m.
$$
For $1\le i\le n$ there exists $\gamma_i^y(t)\in U''(t)$ such that $y_i(t)=\gamma_i^y(t)
\tilde{\overline y}_i(t)$,
$$
\gamma_i^y(t)\equiv 1\text{ mod } m(U''(t)).
$$
\end{Theorem}

\begin{pf} The proof of this theorem is similar to that of Theorem 4.11 \cite{C1}.
References to Theorems 4.8, 4.9 and 4.10 must be replaced with references to
Theorems \ref{Theorem4.8}, \ref{Theorem4.9}, \ref{Theorem4.10} of this paper.
The argument of Lines 11-13 of page 95 \cite{C1}, ``By Theorem 2.12 $\ldots\ldots$ $\text{mult }
\Sigma(0,\ldots,0,\overline y_{l+1},0\ldots,0)=1$'', must be replaced with:

``By Lemma \ref{Lemma3} and Theorem \ref{Theorem4.8}
(with $f=x_{\overline r+l+1}$ in (\ref{eq62})) and Theorem \ref{Theorem4.9},
we may assume that
$$
x_{\overline r+l+1}=\overline x_{\overline r+l+1}
=P+\overline y_1^{d_1}\cdots\overline y_{\overline s}^{d_{\overline s}}\Sigma_0
$$
where $P\in k(U''(0))[\overline y_1,\ldots,\overline y_{\overline r+l}]$ and
$P\in L[\overline x_1^{\frac{1}{d}},\ldots,\overline x_{\overline r}^{\frac{1}{d}},
\overline x_{\overline r+1},\ldots,\overline x_{\overline r+l}]$, with $L$ a finite
extension of the algebraic closure of $k(T''(0))$ in $k(U''(0))$, $d$ a natural number.
$\Sigma_0$ is a series with $\text{mult }
\Sigma(0,\ldots,0,\overline y_{\overline r+l+1},0\ldots,0)=1$''.

On lines 21-31 of page 95 \cite{C1}, ``Suppose that $\nu(P)<\infty$ $\ldots\ldots$
\linebreak $g\in k(c_0)[[\overline x_1(1),\ldots,\overline x_l(1)]][x_{l+1}]$'' must be replaced with:

``Suppose that $\nu^*(P)<\infty$. Let $\omega$ be a primitive $d$th root of unity in an
algebraic closure  of $L$. Set
$$
g'=\prod_{i_1,\ldots,i_{\overline r}=1}^d(x_{\overline r+l+1}-P(\omega^{i_1}\overline x_1^{\frac{1}{d}},\ldots,\omega^{i_{\overline r}}\overline x_{\overline r}^{\frac{1}{d}},
\overline x_{\overline r+1},\ldots,\overline x_{\overline r+l})).
$$
$$
g'\in L[\overline x_1,\ldots,\overline x_{\overline r+l},x_{\overline r+l+1}].
$$
Let $G$ be the Galois group of a Galois closure of $L$ over $k(c_0)$. 
We can define
$$
g=\prod_{\tau\in G}\tau(g')
$$
where $G$ acts on the coefficients of $g'$.
$$
g\in k(c_0)[\overline x_1,\ldots,\overline x_{\overline r+l},x_{\overline r+l+1}].''
$$

On page 97, lines 8-20 of \cite{C1}, ``Set $(e_{ij})=\,\,\,\ldots\ldots\,\,\,
x_l(\alpha+1)=\tilde{\overline y}_l(\alpha+1)$'' should be replaced with:
``After possibly interchanging $\tilde{\overline y}_1(\alpha+1),\ldots,\tilde{\overline y}_{\overline s}(\alpha+1)$, we may assume that $\text{Det}(\tilde C)\ne 0$ where
$$
\tilde C=\left(\begin{array}{lll}
c_{11}(\alpha+1)&\cdots&c_{1\overline r}(\alpha+1)\\
\vdots&&\vdots\\
c_{\overline r1}&\cdots&c_{\overline r\overline r}(\alpha+1)
\end{array}
\right)
$$
Set $(e_{ij})=\tilde C^{-1}$, $d=\text{Det}(\tilde C)$. We can replace $\tilde{\overline y}_i(\alpha+1)$ with 
$$
\tilde{\overline y}_i(\alpha+1)\gamma_1(\alpha+1)^{e_{i1}}\cdots\gamma_{\overline r}(\alpha+1)^{e_{i\overline r}}
$$
for $1\le i\le \overline r$, $\tilde{\overline y}_i(\alpha+1)$ with $\tilde{\overline y}_i(\alpha+1)\gamma_i(\alpha+1)$ for $\overline s+1\le i\le \overline s+l$ and replace
$U''(\alpha+1)$ with $U''(\alpha+1)[\gamma_1(\alpha+1)^{\frac{1}{d}},\ldots,
\gamma_{\overline r}(\alpha+1)^{\frac{1}{d}}]_q$ where 
$$
q=m(U(\alpha+1))\cap
\left(U''(\alpha+1)[\gamma_1(\alpha+1)^{\frac{1}{d}},\ldots,\gamma_{\overline r}(\alpha+1)^{\frac{1}{d}}]\right)
$$
to get
$$
\begin{array}{ll}
x_1(\alpha+1)&=\tilde{\overline y}_1(\alpha+1)^{c_{11}(\alpha+1)}\cdots
\tilde{\overline y}_{\overline s}(\alpha+1)^{c_{1\overline s}(\alpha+1)}\phi_1(\alpha+1)\\
&\vdots\\
x_{\overline r}(\alpha+1)&=\tilde{\overline y}_1(\alpha+1)^{c_{\overline r1}(\alpha+1)}\cdots
\tilde{\overline y}_{\overline s}(\alpha+1)^{c_{\overline r\overline s}(\alpha+1)}\phi_{\overline r}(\alpha+1)\\
&\vdots\\
x_{\overline r+1}(\alpha+1)&=\tilde{\overline y}_{\overline s+1}(\alpha+1)\\
x_{\overline r+l}(\alpha+1)&=\tilde{\overline y}_{\overline s+l}(\alpha+1)\text{''}
\end{array}
$$

On page 98, line 18 to page 99 line 8, ``By construction, $\ldots\ldots$ 
$x_s^*(\alpha+2)=\hat y_1(\alpha+2)^{c_{s1}(\alpha+2)}\cdots \hat y_s(\alpha+2)
^{c_{ss}(\alpha+2)}\psi_s$'' should be replaced with ``By construction there are positive
integers $f_{ij}$ such that
$$
\begin{array}{ll}
x_1^*(\alpha+2)&=\overline y_1(\alpha+2)^{f_{11}}\cdots\overline y_{\overline s}(\alpha+2)^{f_{1\overline s}}
\gamma^{e_{1,\overline r+1}}\tau^{e_{1,\overline r+1}}\\
&\,\,\,\,\,\cdot\phi_1(\alpha+1)^{e_{11}}\cdots
\phi_{\overline r}(\alpha+1)^{e_{1\overline r}}\\
&\vdots\\
x^*_{\overline r}(\alpha+2)&=\overline y_1(\alpha+2)^{f_{\overline r1}}\cdots
\overline y_{\overline s}(\alpha+2)^{f_{\overline r\overline s}}\gamma^{e_{\overline r,\overline r+1}}\tau^{e_{\overline r,\overline r+1}}\\
&\,\,\,\,\,\cdot\phi_1(\alpha+1)^{e_{\overline r1}}\cdots
 \phi_{\overline r}(\alpha+1)^{e_{\overline r\overline r}}\\
x_{\overline r+l+1}^*(\alpha+2)+c_{\alpha+2}&=
\overline y_1(\alpha+2)^{f_{\overline r+1,1}}\cdots
\overline y_{\overline s}(\alpha+2)^{f_{\overline r+1,\overline s}}\gamma^{e_{\overline r+1,\overline r+1}}\tau^{e_{\overline r+1,\overline r+1}}\\
&\,\,\,\,\,\cdot\phi_1(\alpha+1)^{e_{\overline r+1,1}}\cdots
 \phi_{\overline r}(\alpha+1)^{e_{\overline r+1,\overline r}}
\end{array}
$$
in $S(\alpha+2)\sphat$\,\,. $\nu(x^*_{\overline r+1}(\alpha+2)+c_{\alpha+2})=0$ implies
$$
f_{\overline r+1,1}=\cdots=f_{\overline r+1,\overline s}=0.
$$
Set
$$
\tilde\omega=
\phi_1(\alpha+1)^{e_{\overline r+1,1}}\cdots\phi_{\overline r}(\alpha+1)^{e_{\overline r+1,\overline r}}\tau^{e_{\overline r+1,\overline r+1}}\in k(U(\alpha+1)).
$$
Substituting 
$$
\gamma=\overline P_{\alpha+1}+\tilde{\overline y}_1(\alpha+1)^{\epsilon_1(\alpha+1)}\cdots
\tilde{\overline y}_{\overline s}(\alpha+1)^{\epsilon_{\overline s}(\alpha+1)}\Sigma_{\alpha+1}
$$
we have
$$
\begin{array}{ll}
x_{\overline r+l+1}^*(\alpha+2)+c_{\alpha+2}&=\tilde\omega
(\overline P_{\alpha+1}+\tilde{\overline y}_1(\alpha+1)^{\epsilon_1(\alpha+1)}\cdots
\tilde{\overline y}_{\overline s}(\alpha+1)^{\epsilon_{\overline s}(\alpha+1)}\Sigma_{\alpha+1})
^{e_{\overline r+1,\overline r+1}}\\
&=Q_0(\overline y_1(\alpha+2),\ldots,\overline y_{\overline s+l}(\alpha+2))\\
&\,\,\,\,\,
+\overline y_1(\alpha+2)^{\epsilon_1(\alpha+2)}\cdots\overline y_{\overline s}(\alpha+2)
^{\epsilon_{\overline s}(\alpha+2)}\Lambda_0(\overline y_1(\alpha+2),\ldots,\overline y_n(\alpha+2))
\end{array}
$$
where $Q_0$ is a unit series and 
$$
\text{mult }\Lambda_0(0,\ldots,0,\overline y_{\overline s+l+1}(\alpha+2),0,\ldots,0)=1.
$$
The $\overline r\times\overline s$ matrix $(f_{ij})$ with $1\le i\le \overline r$,
$1\le j\le \overline s$, has rank $\overline r$, so after
possibly reindexing $\overline y_1(\alpha+2),\ldots,\overline y_{\overline s}(\alpha+2)$,
we may assume that
$$
\left(\begin{array}{lll}
f_{11}&\cdots& f_{1\overline r}\\
\vdots&&\vdots\\
f_{\overline r1}&\cdots&f_{\overline r\overline r}
\end{array}\right)
$$
has rank $\overline r$. Define $\alpha_i\in{\bf Q}$ by
$$
\left(\begin{array}{l}
\alpha_1\\ \vdots\\ \alpha_{\overline r}\end{array}\right)=
\left(\begin{array}{lll}
f_{11}&\cdots& f_{1\overline r}\\
\vdots&&\vdots\\
f_{\overline r1}&\cdots&f_{\overline r\overline r}
\end{array}\right)^{-1}
\left(\begin{array}{l}
-e_{1,\overline r+1}\\ \vdots\\ -e_{\overline r,\overline r+1}
\end{array}\right)
$$
and set
$$
\hat y_i(\alpha+2)=\left\{
\begin{array}{ll}
\gamma^{-\alpha_i}\overline y_i(\alpha+2)&\text{ for }1\le i\le \overline r\\
\overline y_i(\alpha+2)&\text{ for }\overline r<i
\end{array}
\right.
$$
to get
$$
\begin{array}{ll}
x_1^*(\alpha+2)&=\hat y_1(\alpha+2)^{c_{11}(\alpha+2)}\cdots\hat y_{\overline s}(\alpha+2)^{c_{1\overline s}(\alpha+2)}\psi_1\\
&\vdots\\
 x_{\overline r}^*(\alpha+2)&=\hat y_1(\alpha+2)^{c_{\overline r1}(\alpha+2)}\cdots\hat y_{\overline s}(\alpha+2)^{c_{\overline r\overline s}(\alpha+2)}\psi_{\overline r}
\end{array}
$$
where $c_{ij}(\alpha+2)=f_{ij}$, $(\hat y_1(\alpha+2),\ldots,\hat y_i(\alpha+2))$ are
regular parameters in 
$S(\alpha+2)\sphat$, $\psi_1,\ldots,\psi_{\overline r}\in k(S(\alpha+2))$.''

On page 107,  lines 3-12 substitute for ``$x_1(t)=\ldots\ldots i=l+1$'' the following:
$$
\begin{array}{ll}
``x_1(t')&=y_1(t')^{c_{11}(t')}\ldots y_{\overline s}(t')^{c_{1\overline s}(t')}\tau_1(t')\\
&\vdots\\
x_{\overline r}(t')&=y_1(t')^{c_{\overline r1}(t')}\ldots y_{\overline s}(t')^{c_{\overline r\overline s}(t')}\tau_{\overline r}(t')\\
x_{\overline r+1}(t')&=y_{\overline s+1}(t')\\
&\vdots\\
x_{\overline r+l}(t')&=y_{\overline s+l}(t')\\
x_{\overline r+l+1}(t')&=y_1(t')^{d_1(t')}\cdots y_{\overline s}(t')^{d_{\overline s}(t')}
\overline y_{\overline s+l+1}(t')
\end{array}
$$
Let $\phi_i(t')$ be the residue of $\tau_i(t')$ in $k(S(t'))$, 
$$
\overline \tau_i=\frac{\tau_i(t')}{\phi_i(t')}.
$$
After possibly reindexing $y_1(t'),\ldots, y_{\overline s}(t')$, we may assume that
$$
\tilde C=
\left(\begin{array}{lll}
c_{11}(t')&\cdots&c_{1\overline r}(t')\\
\vdots&&\vdots\\
c_{\overline r1}(t')&\cdots&c_{\overline r\overline r}(t')
\end{array}\right)
$$
has rank $\overline r$. Let $(e_{ij})=\tilde C^{-1}$. Define
$$
\overline y_i(t')=\left\{\begin{array}{ll}
 \overline \tau_1^{e_{i1}}\cdots\overline \tau_{\overline r}^{e_{i\overline r}}y_i(t')&\text{ if }1\le i\le \overline r\\
y_i(t')&\text{ if }\overline r<i,i\ne \overline s+l+1\\
\overline \tau_1^{-e_{11}d_1(t')-\cdots-e_{\overline r1}d_{\overline r}(t')}\cdots
\overline \tau_{\overline r}^{-e_{1\overline r}d_1(t')-\cdots-e_{\overline r\overline r}d_{\overline r}(t')}
y_{\overline s+l+1}(t')&\text{ if }i=\overline r+l+1''
\end{array}
\right.
$$
\end{pf}

\section{Monomialization}

\begin{Theorem}\label{Theorem5.1} Suppose that $k$ is a field of
characteristic zero, $K\rightarrow K^*$ is a (possibly transcendental) extension
of algebraic function fields over $k$.  Suppose that $\nu^*$ is a rank 1 valuation of
$K^*$ which is trivial on $k$. Suppose that $R$ is an algebraic local ring of $K$,
$S$ is an algebraic local ring of $K^*$ such that $S$ dominates $R$ and $\nu^*$ dominates
$S$.  
Let $\nu=\nu^*\mid K$,
$$
\overline s=\text{ratrank }\nu^*\ge\overline r=\text{ratrank }\nu.
$$
Then there exist sequences of monoidal transforms $R\rightarrow R'$ and $S\rightarrow S'$
along $\nu^*$ such that $R'$ and $S'$ are regular local rings, $S'$ dominates $R'$,
there exist regular parameters
$(y_1',\ldots,y_n')$ in $S'$, $(x_1',\ldots,x_m')$ in $R'$, where 
$$
n=\text{trdeg}_kK^*-\text{trdeg}_kk(V^*),
$$
$$
m=\text{trdeg}_kK-\text{trdeg}_kk(V),
$$
$\nu(x_1'),\ldots,\nu(x_{\overline r}')$ is a rational basis of $\Gamma_{\nu}\otimes{\bf Q}$,
$\nu^*(y_1'),\ldots,\nu^*(y_{\overline s}')$ is a rational basis of $\Gamma_{\nu^*}\otimes{\bf Q}$, there are
units $\delta_1,\ldots,\delta_{\overline r}\in S'$ and an $\overline r\times\overline s$
matrix $(c_{ij})$ of nonnegative integers such that $(c_{ij})$ has rank $\overline r$, and
$$
\begin{array}{ll}
x_1'&=(y_1')^{c_{11}}\cdots(y_{\overline s}')^{c_{1\overline s}}\delta_1\\
&\vdots\\
x_{\overline r}'&=(y_1')^{c_{\overline r1}}\cdots(y_{\overline s}')^{c_{\overline r\overline s}}
\delta_{\overline r}\\
x_{\overline r+1}'&=y_{\overline s+1}'\\
&\vdots\\
x_m'&=y_{\overline s+m-\overline r}'.
\end{array}
$$
\end{Theorem}

\begin{pf} $k(V)$ and $k(V^*)$ have finite transcendence degree over $k$ by Theorem 1 \cite{Ab2}
or Appendix 2 \cite{ZS}. We have $\text{ rank }\nu\le \text{ rank }\nu^*=1$.
By Hironaka's theorems on resolution, resolution of singularities Theorem $I_2^{m,n}$ \cite{H1}
(c.f. Theorem 2.9 \cite{C1}) and resolution of indeterminancy (c.f. Theorem 2.6 [C1],
the statement and proof of Theorem 2.6 are valid if $R$ is not regular) we can assume that $R$ and $S$ are
regular local rings.

By  resolution of indeterminancy (c.f. Theorem 2.6 \cite{C1}) and Theorem 2.7
\cite{C1}, applied to a lift to $V$ of a transcendence basis of $k(V)$ over $k$, and Theorem 2.7
\cite{C1} applied to a lift to $V^*$ of a transcendence basis of $k(V^*)$
over $k$, there exist
MTSs along $\nu^*$ $R\rightarrow R(1)$ and $S\rightarrow S(1)$ such that $R(1)$ and $S(1)$
are regular local rings, $V^*$ dominates $S(1)$, $S(1)$ dominates $R(1)$ and
$$
\text{trdeg}_{k(R(1))}k(V)=0,
$$
$$
\text{trdeg}_{k(S(1))}k(V^*)=0.
$$

First assume that $\text{rank }\nu=1$. 
Let $\{t_1,\ldots,t_{\beta}\}$ be a lift of a transcendence basis of $k(R(1))$ over $k$ to
$R(1)$. Let $L=k(t_1,\ldots,t_{\beta})\subset R(1)$. By replacing $k$ with $L$,
we may assume that 
$$
\text{trdeg}_kk(R(1))=\text{trdeg}_kk(V)=0.
$$
There exist $f_1,\ldots,f_{\overline r}\in K$ such that $\nu(f_1),\ldots,
\nu(f_{\overline r})$ are positive and rationally independent. By Theorem
2.7 \cite{C1}, there exists a MTS $R(1)\rightarrow R(2)$ along $\nu$ such
that $f_1,\ldots,f_{\overline r}\in R(2)$. By Theorem 2.5 \cite{C1},
there exists a MTS $R(2)\rightarrow R(3)$ along $\nu$ such that 
$f_1\cdots f_{\overline r}$ is a SNC divisor in $R(3)$. Thus $R(3)$ has
regular parameters $(x_1(3),\ldots,x_m(3))$ such that
$\nu(x_1(3)),\ldots,\nu(x_{\overline r}(3))$ are  a rational basis of $\Gamma_{\nu}\otimes{\bf Q}$.

By Theorem 2.6, there exists a MTS $S(1)\rightarrow S(2)$ along $\nu^*$
such that $S(2)$ dominates $R(3)$.

As in the construct of $R(1)\rightarrow R(3)$ there exists a MTS
$S(2)\rightarrow S(3)$ along $\nu^*$ such that $S(3)$ has regular
parameters $(y_1(3),\ldots,y_n(3))$ such that 
$\nu^*(y_1(3)),\ldots,\nu^*(y_{\overline s}(3))$ are a rational basis of $\Gamma_{\nu^*}\otimes{\bf Q}$.

 By (\ref{eq60}) of Theorem \ref{Theorem4.8}, with the $R$, $S$, $f$, $\overline m$, $l$
 of the hypothesis of that theorem set to $R=S(3)$ and $S=S(3)$, $f=x_1(3)\cdots x_{\overline r}(3)$, $\overline m=n-\overline s$, $l=0$,
and by (4) of Theorem \ref{Theorem4.10}, there exists a MTS $S(3)\rightarrow S(4)$
along $\nu^*$ such that
$$
x_i(3)=y_1(4)^{c_{i1}}\cdots y_{\overline s}(4)^{c_{i\overline s}}\psi_i
$$
where $\psi_i\in S(4)$ are units for $1\le i\le \overline r$, 
$\nu^*(y_1(4)),\ldots,\nu^*(y_{\overline s}(4))$ are a rational basis of
$\Gamma_{\nu^*}\otimes{\bf Q}$, and $\text{rank }(c_{ij})=\overline r$.

After possibly permuting the first $\overline s$ variables $y_i(4)$, we may assume
that the matrix
$$
\tilde C=\left(\begin{array}{lll}
c_{11}&\cdots&c_{1\overline r}\\
\vdots&&\vdots\\
c_{\overline r1}&\cdots&c_{\overline r\overline r}
\end{array}
\right)
$$
has nonzero determinant.

Let $\phi_i$ be the residue of $\psi_i$ in $k(S(4))$. For $1\le i\le \overline r$, set
$$
\epsilon_i=\prod_{j=1}^{\overline r}\left(\frac{\psi_j}{\phi_j}\right)^{e_{ij}}
$$
where  $(e_{ij})=\tilde C^{-1}$, a matrix with rational coefficients. We have
$\epsilon_j\in S(4)\sphat$ for $1\le j\le \overline r$.

Set
$$
\overline y_j(4)=\left\{
\begin{array}{ll} \epsilon_jy_j(4)& 1\le j\le\overline r\\
y_j(4)&\overline r<j
\end{array}\right.
$$
we have 
$$
x_i(3)=\prod_{j=1}^{\overline s}\overline y_j(4)^{c_{ij}}\phi_i
$$
for $1\le i\le\overline r$.

In the notation of Theorem \ref{Theorem4.11}, set $R=R(3)$, $T''(0)=R(3)$, $x_i=x_i(3)$ for
$1\le i\le m$, $\overline x_i=x_i(3)$ for $1\le i\le m$, $c_0=1$,
$\tilde R=k[x_1(3),\ldots,x_{\overline r}(3)]_q$ where $q=m(R(3))\cap k[x_1(3),\ldots,x_{\overline r}(3)]$, $\gamma_i=1$ for $1\le i\le \overline r$. Set $S=S(4)$,
$U''(0)=S(4)[d_0,\epsilon_1,\ldots,\epsilon_{\overline r}]_p$ where $k(d_0)\cong k(S(4))$,
$p=m(S(4)\sphat)\cap S(4)[d_0,\epsilon_1,\ldots,\epsilon_{\overline r}]$,
$y_i=y_i(4)$ for $1\le i\le n$, $\overline y_i=\overline y_i(4)$ for $1\le i\le n$.

Then the assumptions of Theorem \ref{Theorem4.11} are satisfied with $l=0$.
 By induction on $l$ in Theorem \ref{Theorem4.11},
we construct the desired MTSs, and finish the proof of the Theorem when $\text{rank }\nu=1$.

If $\text{rank }\nu=0$, then $\nu$ is trivial, so that $V=K$, $R=K$ and $\overline r=0$.
We can then construct a MTS $S\rightarrow S'$ along $\nu^*$ as in the first part of the proof,
so that $S'$ has a regular system of parameters $(y_1',\ldots, y_n')$ such that
$\nu^*(y_1'),\ldots,\nu^*(y_{\overline s}')$ is a rational basis of $\Gamma_{\nu^*}\otimes{\bf Q}$ and reach the conclusions of the theorem (see Remark \ref{Remark1} below).
\end{pf}

\begin{Remark}\label{Remark1} The degenerate case of Theorem \ref{Theorem5.1} is when $\nu$ is trivial
on $K$. This case can only occur if $K^*$ is  transcendental over $K$. When $\nu$ is trivial,
$V=K$, so we must have that $R=K$.
The conclusions of Theorem \ref{Theorem5.1} are in this case that there exists a sequence of monoidal transforms $S\rightarrow S'$
along $\nu^*$ such that $S'$ is a  regular local ring (which contains $K$), and there exist regular parameters $(y_1',\ldots,y_n')$ in $S'$, where 
$$
n=\text{trdeg}_kK^*-\text{trdeg}_kk(V^*),
$$
such that $\nu^*(y_1'),\ldots,\nu^*(y_n')$ is a rational basis of $\Gamma_{\nu^*}\otimes{\bf Q}$.
\end{Remark}

We now introduce notation that will be used in the proof of Theorem \ref{Theorem5.3}.
Suppose that $k$ is a field of
characteristic zero, $K\rightarrow K^*$ is a (possibly transcendental) extension
of algebraic function fields over $k$.  Suppose that $\nu^*$ is a  valuation of
$K^*$ of arbitrary (but necessarily finite) rank which is trivial on $k$. Suppose that $R$ is an algebraic local ring of $K$,
$S$ is an algebraic local ring of $K^*$ such that $S$ dominates $R$ and $\nu^*$ dominates
$S$.  Let $\nu=\nu^*\mid K$. Let $V^*$ be the valuation rings of $\nu^*$ and $V$ be the
valuation ring of $\nu$.

\begin{Lemma}\label{Lemma10}
Suppose that $p$ is a prime ideal of $V$. Then there exists a prime ideal $q$ of $V^*$
such that $q\cap V=p$.
\end{Lemma}

\begin{pf} By Theorem 15 of Section 10, Chapter VI \cite{ZS}, there exists an isolated
subgroup $\Delta$ of $\Gamma_{\nu}$ such that
$$
p=\{a\in K\mid\nu(a)=\beta\text{ for some }\beta\in \Gamma_{\nu}-\Delta\text{ with }\beta\ge 0\}\cup\{0\}.
$$
Set 
$$
\Delta^*=\{\beta\in\Gamma_{\nu^*}\mid \mid\beta\mid\le\mid\alpha\mid\text{ for some }
\alpha\in\Delta\}.
$$
$\Delta^*$ is an isolated subgroup of $\Gamma_{\nu^*}$. Theorem 15 
of Section 10, Chapter VI \cite{ZS} implies that 
$$
q\in \{a\in K\mid \nu(a)=\beta\text{ for some }\beta\in \Gamma_{\nu^*}-\Delta^*
\text{ with }\beta\ge 0\}\cup\{0\}
$$
is a prime ideal of $V^*$. 
$$
(\Gamma_{\nu^*}-\Delta^*)\cap \Gamma_{\nu}=\Gamma_{\nu}-\Delta
$$
implies $q\cap V=p$.
\end{pf}

Let $\beta=\text{rank }V$.  The primes of
$V$ are a finite chain
$$
0=p_0\subset\cdots\subset p_{\beta}\subset V.
$$
Note that if $\beta=0$ then $V=K$. The primes of $V^*$ are a finite chain
$$
0=q_{0,1}\subset\cdots\subset q_{0,\sigma(0)}\subset q_{1,1}\subset\cdots \subset q_{\beta,\sigma(\beta)}
$$
where  $p_i=q_{i,j}\cap V$ for $1\le j\le \sigma(i)$.

The isolated subgroups of $\Gamma_{\nu}$ are a chain
$$
0=\Gamma_{\beta}\subset\cdots\subset\Gamma_0=\Gamma_{\nu}.
$$
There is a corresponding chain of isolated subgroups of $\Gamma_{\nu^*}$
$$
0=\Gamma_{\beta,\sigma(\beta)}\subset \cdots\subset\Gamma_{0,1}=\Gamma_{\nu^*}.
$$
For $i\le j$, $\nu$ induces a valuation on the field $(V/p_i)_{p_i}$
with valuation ring $(V/p_i)_{p_j}$ and value group $\Gamma_i/\Gamma_j$.
If $j=i+1$ then $\Gamma_i/\Gamma_j$ has rank 1.
For $i< j$, $1\le a\le\sigma(i)$ and $1\le b\le\sigma(j)$, $\nu^*$ induces a valuation on
the field $(V^*/q_{i,a})_{q_{i,a}}$ with valuation ring   $(V^*/q_{i,a})_{q_{j,b}}$
and value group $\Gamma_{i,a}/\Gamma_{j,b}$. If $i=j$ and $1\le a\le b\le\sigma(i)$,
$\nu^*$ induces a valuation on
the field $(V^*/q_{i,a})_{q_{i,a}}$ with valuation ring   $(V^*/q_{i,a})_{q_{i,b}}$
and value group $\Gamma_{i,a}/\Gamma_{i,b}$.

We have dominant inclusions  of valuation rings
$$
(V/p_i)_{p_j}\rightarrow (V^*/q_{ia})_{q_{jb}}
$$
if $i< j$, $1\le a\le\sigma(i)$, $1\le b\le\sigma(j)$ which induce inclusions of
valuation groups
$$
\Gamma_i/\Gamma_j\rightarrow\Gamma_{i,a}/\Gamma_{j,b}.
$$

We also have dominant inclusions  of valuation rings
$$
(V/p_i)_{p_i}\rightarrow (V^*/q_{ia})_{q_{ib}}
$$
if $i= j$, $1\le a\le b\le\sigma(i)$.
Note that the value group of the field $(V/p_i)_{p_i}$ is $\Gamma_i/\Gamma_i=0$.

\begin{Lemma}\label{Lemma5}
There exist MTSs along $\nu^*$
$$
\begin{array}{lll}
R'&\rightarrow&S'\\
\uparrow&&\uparrow\\
R&\rightarrow&S
\end{array}
$$
such that  $R'$ and $S'$ are regular, 
$$
\text{trdeg}_{k(R'_{p_i\cap R'})}k(V_{p_i})=0
$$
for all $i$ and
$$
\text{trdeg}_{k(S'_{q_{ij}\cap S'})}k(V^*_{q_{ij}})=0
$$
for all $i,j$.
\end{Lemma}

\begin{pf} By Hironaka's theorem on resolution of singularities (Theorem $I_2^{m,n}$\cite{H1}
or Theorem 2.9 \cite{C1}) and resolution of indeterminancy (c.f. Theorem 2.6 \cite{C1},
the statement and proof are valid if $R$ is not regular) we can assume
that $R$ and $S$ are regular local rings.

For all $i$, $V_{p_i}$ is a valuation ring of $K$ dominating $R_{p_i\cap R}$. Thus
$$
\text{trdeg}_{(R/p_i\cap R)_{p_i\cap R}}(V/p_i)_{p_i}<\infty
$$
by Theorem 1 \cite{Ab2} or Appendix 2 \cite{ZS}. We can lift transcendence bases of
$(V/p_i)_{p_i}$ over $(R/p_i\cap R)_{p_i\cap R}$ for $1\le i\le\beta$
to $t_1,\ldots,t_a\in V$. After possibly replacing the $t_i$ with $\frac{1}{t_i}$,
we have $\nu(t_i)\ge 0$ for all $t_i$.

By Theorem 2.7 \cite{C1}, there exists a MTS $R\rightarrow R'$ along $\nu$ such that
$t_i\in R'$ for all $i$. Let $p_i'=R'\cap p_i$. Then
$$
\text{trdeg}_{(R'/p_i')_{p_i'}}(V/p_i)_{p_i}=0\text{ for }1\le i\le \beta.
$$
By Theorem 2.6 \cite{C1}, there exists a MTS $S\rightarrow S''$ along $\nu^*$ such that
$S''$ dominates $R'$. As argued above for $R$, there exists a MTS $S''\rightarrow S'$ along
$\nu^*$ such that if $q_{ij}'=S'\cap q_{ij}$, then
$$
\text{trdeg}_{(S'/q_{ij}')_{q_{ij}'}}(V^*/q_{ij})_{q_{ij}}=0\text{ for all }i,j.
$$
\end{pf}

\begin{Theorem}\label{Theorem5.3} Let notation be as above.
Suppose that  $R$ and $S$ are regular, 
$$
\text{trdeg}_{k(R_{p_i\cap R})}k(V_{p_i})=0
$$
for all $i$ and
$$
\text{trdeg}_{k(S_{q_{ij}\cap S})}k(V^*_{q_{ij}})=0
$$
for all $i,j$.
Suppose that the rank 1 valuation groups $\Gamma_{i-1}/\Gamma_i$ has rational rank $\overline r_i$ for $1\le i\le\beta$,
$\Gamma_{i-1,\sigma(i-1)}/\Gamma_{i,1}$ have rational rank $\overline s_i=\overline s_{i1}$
for $1\le i\le \beta$ and $\Gamma_{i,a-1}/\Gamma_{i,a}$ have rational rank $\overline s_{ia}$
for $1\le i\le \beta$ and $2\le a\le\sigma(i)$.

Set $t_i=\text{dim}(R/p_{i-1}\cap R)_{p_i\cap R}$ for $1\le i\le\beta$, so that 
$$
m=\text{dim }R=t_1+\cdots +t_{\beta}.
$$
 Set
$$
\overline t_{ij}=\left\{\begin{array}{ll}
0&\text{ if }i=0, j=1\\
\text{dim }(S/q_{i-1,\sigma(i-1)}\cap S)_{q_{i,1}\cap S}&\text{ if } 1\le i\le\beta, j=1\\
\text{dim }(S/q_{i,j-1}\cap S)_{q_{i,j}\cap S}&\text{ if } 0\le i\le\beta, 2\le j\le\sigma(i).
\end{array}\right.
$$
For $1\le i\le \beta$ set
$$
\overline t_i=\overline t_{i,1}+\cdots+\overline t_{i,\sigma(i)}=\text{dim }(S/q_{i-1,\sigma(i-1)}\cap S)_{q_{i\sigma(i)}\cap S},
$$
and set
$$
\overline t_0=\overline t_{0,1}+\cdots+\overline t_{0,\sigma(0)}=\text{ dim }
S_{q_{0,\sigma(0)}\cap S}
$$
so that $n=\text{ dim }S=\overline t_0+\cdots+\overline t_{\beta}$.

Then there exist sequences of monoidal transforms $R\rightarrow R'$ and $S\rightarrow S'$ 
 such that $V^*$ dominates $S'$,
$S'$ dominates $R'$, $R'$ has regular parameters $(z_1,\dots,z_m)$, $S'$ has regular 
parameters $(w_1,\ldots,w_n)$ and there are units $\delta_j\in S'$  such that
$$
p_i\cap R'=(z_1,\ldots,z_{t_1+\cdots +t_i})
$$
 for $1\le i\le \beta$
and
$$
q_{ij}\cap S'=(w_1,\ldots,w_{\overline t_0+\cdots+\overline t_{i-1}+\overline t_{i,1}+
\cdots+\overline t_{i,j}})
$$
for $0\le i\le \beta$ and $1\le j\le \sigma(i)$.
$$
\nu(z_{t_1+\cdots+t_{i-1}+1}),\cdots,\nu(z_{t_1+\cdots+t_{i-1}+\overline r_i})
$$
is a rational basis of $\Gamma_{i-1}/\Gamma_i\otimes{\bf Q}$ for $1\le i\le \beta$,
$$
\nu^*(w_{\overline t_0+\cdots+\overline t_{i-1}+1}),\ldots,\nu^*(w_{\overline t_0+
\cdots+\overline t_{i-1}+\overline s_i})
$$
is a rational basis of $\Gamma_{i-1,\sigma(i-1)}/\Gamma_{i,1}\otimes{\bf Q}$ for $1\le i\le \beta$
and
$$
\nu^*(w_{\overline t_0+\cdots+\overline t_{i-1}+\overline t_{i,1}+\cdots+\overline t_{i,j-1}+1}),\ldots,\nu^*(w_{\overline t_0+
\cdots+\overline t_{i-1}+\overline t_{i,1}+\cdots+\overline t_{i,j}+\overline s_{i,j}})
$$
is a rational basis of $\Gamma_{i,j-1}/\Gamma_{i,j}\otimes{\bf Q}$ for $0\le i\le \beta$
and $2\le j\le \sigma(i)$. Furthermore

$$
\begin{array}{ll}
z_1&=w_{\overline t_0+1}^{g_{11}(1)}\cdots w_{\overline t_0+\overline s_{1}}^{g_{1\overline s_{1}}(1)}
w_{\overline t_0+\overline t_{1,1}+1}^{h_{1,\overline t_0+\overline t_{1,1}+1}(1)}\cdots w_n^{h_{1,n}(1)}\delta_1\\
&\vdots\\
z_{\overline r_1}&=w_{\overline t_0+1}^{g_{\overline r_1,1}(1)}\cdots w_{\overline t_0+\overline s_1}^{g_{\overline r_1,\overline s_1}(1)}
w_{\overline t_0+\overline t_{1,1}+1}^{h_{t_1,\overline t_0+\overline t_{1,1}+1}(1)}\cdots w_n^{h_{\overline r_1,n}(1)}\delta_{\overline r_1}\\
z_{\overline r_1+1}&=w_{\overline t_0+\overline s_1+1}w_{\overline t_0+\overline t_{1,1}+1}^{h_{\overline r_1+1,\overline t_0+\overline t_{1,1}+1}(1)}\cdots w_n^{h_{\overline r_1+1,n}(1)}\delta_{\overline r_1+1}\\
&\vdots\\
z_{t_1}&=w_{\overline t_0+\overline s_1+t_1-\overline r_1} w_{\overline t_0+\overline t_{1,1}+1}^{h_{t_1,\overline t_0+\overline t_{1,1}+1}(1)}\cdots
w_n^{h_{t_1,n}(1)}\delta_{t_1}\\
z_{t_1+1}&=w_{\overline t_0+\overline t_1+1}^{g_{1,1}(2)}\cdots w_{\overline t_0+\overline t_1+\overline s_2}^{g_{1,\overline s_2}(2)}
w_{\overline t_0+\overline t_1+\overline t_{2,1}+1}^{h_{1,\overline t_0+\overline t_1+\overline t_{2,1}+1}(2)}
\cdots w_n^{h_{1,n}(2)}\delta_{t_1+1}\\
&\vdots\\
z_{t_1+\overline r_2}&=w_{\overline t_0+\overline t_1+1}^{g_{\overline r_2,1}(2)}\cdots w_{\overline t_0+\overline t_1+\overline s_2}^{g_{\overline r_2,\overline s_2}(2)}
w_{\overline t_0+\overline t_1+\overline t_{2,1}+1}
^{h_{\overline r_2,\overline t_0+\overline t_1+\overline t_{2,1}+1}(2)}\cdots
w_n^{h_{\overline r_2,n}(2)}\delta_{t_1+\overline r_2}\\
z_{t_1+\overline r_2+1}&=w_{\overline t_0+\overline t_1+\overline s_2+1} w_{\overline t_0+\overline t_1+\overline t_{2,1}+1}
^{h_{\overline r_2+1,\overline t_0+\overline t_1+\overline t_{2,1}+1}(2)}\cdots w_n^{h_{\overline r_2+1,n}(2)}\delta_{t_1+\overline r_2+1}\\
&\vdots\\
z_{t_1+t_2}&=w_{\overline t_0+\overline t_0+t_2+\overline s_2-\overline r_2} w_{\overline t_0+\overline t_1+\overline t_{2,1}+1}
^{h_{t_2,\overline t_0+\overline t_1+\overline t_{2,1}+1}(2)}\cdots w_n^{h_{t_2,n}(2)}\delta_{t_1+t_2}\\
&\vdots\\
z_{t_0+\cdots+t_{\beta-1}+1}&=w_{\overline t_0+\cdots+\overline t_{\beta-1}+1}
^{g_{11}(\beta)}\cdots w_{\overline t_0+\cdots+\overline t_{\beta-1}+\overline s_{\beta}}
^{g_{1,\overline s_{\beta}}(\beta)}
w_{\overline t_0+\cdots+\overline t_{\beta,1}+1}^{h_{1,\overline t_0+\cdots+\overline t_{\beta,1}+1}(\beta)}\cdots w_n^{h_{1,n}(\beta)}\delta_{t_1+\cdots+t_{\beta-1}+1}\\
&\vdots\\
z_{t_1+\cdots+t_{\beta-1}+\overline r_{\beta}}&=w_{\overline t_0+\cdots+\overline t_{\beta-1}+1}
^{g_{\overline r_{\beta},1}(\beta)}\cdots w_{\overline t_0+\cdots+\overline t_{\beta-1}+\overline s_{\beta}}
^{g_{\overline r_{\beta},\overline s_{\beta}}(\beta)}
w_{\overline t_0+\cdots+\overline t_{\beta,1}+1}^{h_{\overline r_{\beta},\overline t_0+\cdots+t_{\beta,1}+1}(\beta)}\cdots w_n^{h_{\overline r_{\beta},n}(\beta)}\delta_{t_1+\cdots+t_{\beta-1}+\overline r_{\beta}}\\
z_{t_1+\cdots+t_{\beta-1}+\overline r_{\beta}+1}&=
w_{\overline t_0+\cdots+\overline t_{\beta-1}+\overline s_{\beta}+1}
w_{\overline t_0+\cdots+\overline t_{\beta,1}+1}^{h_{\overline r_{\beta}+1,\overline t_0+\cdots+\overline t_{\beta,1}+1}(\beta)}\cdots w_n^{h_{\overline r_{\beta}+1,n}(\beta)}\delta_{t_1+\cdots+t_{\beta-1}+\overline r_{\beta}+1}\\
z_{t_1+\cdots+t_{\beta}}&=
w_{\overline t_0+\cdots+\overline t_{\beta-1}+t_{\beta}+\overline s_{\beta}-\overline r_{\beta}}
w_{\overline t_0+\cdots+\overline t_{\beta,1}+1}^{h_{t_{\beta},\overline t_0+\cdots+\overline t_{\beta,1}+1}(\beta)}\cdots w_n^{h_{t_{\beta},n}(\beta)}\delta_{t_1+\cdots+t_{\beta}}
\end{array}
$$
Where 
$$
\text{rank }
\left(\begin{array}{lll}
g_{11}(i)&\cdots&g_{1\overline s_i}(i)\\
\vdots&&\vdots\\
g_{\overline r_i1}(i)&\vdots&g_{\overline r_i\overline s_i}(i)
\end{array}\right)
=\overline r_i
$$
for $1\le i\le \beta$.
\end{Theorem}

\begin{pf} 
We prove the theorem by induction on $\text{rank }V^*$. If $\text{rank }V^*=1$, then the
theorem is immediate from Theorem \ref{Theorem5.1}.

By induction on $\gamma=\text{ rank }V^*$, we may assume that the theorem is true whenever 
$\text{ rank }V^*=\gamma-1$. We are reduced to proving the theorem in the following two 
cases.
\begin{description}
\item[Case 1] $\sigma(\beta)=1$
\item[Case 2] $\sigma(\beta)>1$.
\end{description}

Suppose that we are in Case 1. Then $V^*/q_{\beta-1,\sigma(\beta-1)}$ is a rank 1 valuation
ring which dominates the rank 1 valuation ring $V/p_{\beta-1}$. 
$V^*/q_{\beta-1,\sigma(\beta-1)}$ has rational rank $\overline s_{\beta}$ and
$V/p_{\beta-1}$ has rational rank $\overline r_{\beta}$.

The proof of Theorem \ref{Theorem5.3} in case 1 follows from the proof of Theorem
5.3 \cite{C1} with some changes in notation and references to
supporting lemmas and theorems. We must replace $r$ with $\beta$, $p_{r-1}(i)$
with $p_{\beta-1}\cap R(i)$ and $q_{r-1}(i)$ with $q_{\beta-1,\sigma(\beta-1)}\cap S(i)$
throughout the proof.  Then $R(i)_{p_{\beta-1}(i)}$ has a system of
$\lambda=t_1+\cdots+t_{\beta-1}$ regular parameters, while
$S(i)_{q_{\beta-1}}(i)$ has a system of $\overline\lambda=\overline t_0+\cdots+
\overline t_{\beta-1}\ge\lambda$ regular parameters.

References to Theorems 4.8, 4.10 and 5.1 of \cite{C1}
must be replaced with references to Theorems  \ref{Theorem4.8}, \ref{Theorem4.10} and \ref{Theorem5.1} of this paper.

Now suppose that we are Case 2, $\sigma(\beta)>1$. Then
$V^*/q_{\beta,\sigma(\beta)-1}$ is a rank 1, rational rank 
$\overline s_{\beta,\sigma(\beta)}$ valuation ring which dominates the 
rank 0 valuation ring $V/p_{\beta}$. That is, $V/p_{\beta}$ is a field.
The proof in Case 2 is thus a substantial simplification of the proof in Case 1.
\end{pf}

The proof of Theorem \ref{Theorem2} is now an immediate corollary.

\ \\ \\
\noindent
Steven Dale Cutkosky, Department of Mathematics, University of
Missouri\\
Columbia, MO 65211, USA\\
cutkoskys@@missouri.edu

\end{document}